\documentclass[11pt, letterpaper, leqno]{article}
\usepackage{amsmath,amsthm,amscd,amssymb,amsbsy,cite,amstext,mathtools}
\usepackage{pgf}
\usepackage{color}
\usepackage{enumerate}
\usepackage{mathrsfs}
\usepackage[symbol]{footmisc}
\usepackage[colorlinks=true, pdfstartview=FitV, linkcolor=black, citecolor=orange, urlcolor=blue]{hyperref}
\usepackage{cleveref}

\usepackage{amsfonts}
\usepackage[T1]{fontenc}
\usepackage{euler}

\theoremstyle{plain}
\newtheorem{theorem}{Theorem}
\newtheorem{lemma}{Lemma}[section]
\newtheorem{corollary}{Corollary}[section]

\newtheorem{claim}{Claim}[section]
\newtheorem{problem}{Problem}

\newcommand{\namedtheorem}[3][2]{
	\theoremstyle{plain}
	\newtheorem*{#1}{#1}
	\begin{#1}[#3]
		#2
	\end{#1}
}

\makeatletter
\newenvironment{subtheorem}[1]{%
	\def\subtheoremcounter{#1}%
	\refstepcounter{#1}%
	\protected@edef\theparentnumber{\csname the#1\endcsname}%
	\setcounter{parentnumber}{\value{#1}}%
	\setcounter{#1}{0}%
	\expandafter\def\csname the#1\endcsname{\theparentnumber.\Alph{#1}}%
	\expandafter\def\csname theH#1\endcsname{theorem.\theparentnumber.\Alph{#1}}%
	\unskip\ignorespaces
}{%
	\setcounter{\subtheoremcounter}{\value{parentnumber}}%
	\ignorespacesafterend
}
\makeatother
\newcounter{parentnumber}

\theoremstyle{remark}
\newtheorem{remark}{Remark}[section]
\newtheorem*{acknowledgment}{Acknowledgment}

\theoremstyle{definition}
\newtheorem{algorithm}{Algorithm}
\newtheorem*{algorithm*}{Algorithm}
\newtheorem*{CZalg}{CZ Algorithm}
\newtheorem{definition}{Definition}[section]

\numberwithin{equation}{section}
\newcommand{\eqindent}{\displayindent0pt\displaywidth\textwidth}

\newcommand{\supp}[1]{\mathrm{supp}\left( {#1} \right)}
\newcommand{\jet}[2]{\mathscr{J}_{#1}{#2}}

\newcommand{\norm}[1]{\|{#1}\|}

\newcommand{\diam}[1]{\mathrm{diam}\left(#1\right)}
\newcommand{\set}[1]{\left\{#1\right\}}
\newcommand{\void}{\varnothing}
\newcommand{\abs}[1]{|#1|}
\newcommand{\dist}[2]{\mathrm{dist}\left( {#1}, {#2} \right)}
\newcommand{\E}{\mathcal{E}}
\newcommand{\Eb}{{\mathcal{E}}}
\newcommand{\Epmb}{{\mathcal{E}}}
\newcommand{\G}{\Gamma_+^{\sharp}}

\newcommand{\sk}{\sigma^{\sharp}}
\renewcommand{\ss}{\sk}
\renewcommand{\P}{\mathcal{P}}
\newcommand{\brac}[1]{\left(#1\right)}
\newcommand{\intr}{{interior}}
\newcommand{\cl}{{closure}}
\newcommand{\ksharp}{k^{\sharp}}
\newcommand{\csharp}{C}
\newcommand{\knice}{\Lambda_{\mathrm{nice}}^{(k)}}
\newcommand{\cnice}{C_{\mathrm{nice}}}
\newcommand{\Q}{\mathcal{Q}}
\newcommand{\R}{\mathbb{R}}
\renewcommand{\hat}[1]{\widehat{#1}}

\renewcommand{\d}{\partial}
\newcommand{\Rn}{\mathbb{R}^n}
\newcommand{\xqs}{x_Q^\sharp}
\renewcommand{\bar}[1]{\overline{#1}}

\begin{document}

	\title{Nonnegative $C^2(\mathbb{R}^2)$ Interpolation}

	\author{Fushuai Jiang\footnote{Department of Mathematics, University of California, Davis, One Shields Ave, Davis, CA 95616, USA \newblock fsjiang@math.ucdavis.edu} \and Garving~K. Luli\footnote{Department of Mathematics, University of California, Davis, One Shields Ave, Davis, CA 95616, USA
			\newblock kluli@math.ucdavis.edu}}

	\maketitle

	\begin{abstract}
		In this paper, we prove two improved versions of the Finiteness Principle for nonnegative $ C^2(\R^2) $ interpolation, previously proven by Fefferman, Israel, and Luli. The first version sharpens the finiteness constant to $ 64 $, and the second version carries better computational practicality. Along the way, we also provide a detailed construction of nonnegative $ C^2 $ interpolants in one-dimension, and prove the nonexistence of a bounded linear $ C^2 $-extension operator that preserves nonnegativity. 
	\end{abstract}

	

	\section{Introduction}
	
	For nonnegative integers $ m,n $, we write $ C^m(\Rn) $ to denote the Banach space of $ m $-times continuously differentiable real-valued functions such that the following norm is finite
	\begin{equation*}
	\norm{F}_{C^m(\Rn)} := \sup\limits_{x \in \Rn}\brac{  \sum\limits_{\abs{\alpha} \leq m} \abs{ \partial^\alpha F(x) }^2   }^{1/2}\,.
	\end{equation*}
	If $ S $ is a finite set, we write $ \#(S) $ to denote the number of elements in $ S $. We use $ C $ to denote constants that depend only on $ m $ and $ n $.

	\begin{problem}\label{prob.norm}
		Let $ E \subset \Rn $ be a finite set. Let $ f : E \to [0,\infty) $. Compute the order of magnitude of 
		\begin{equation}
		\norm{f}_{C^m_+(E)} := \inf \set{ \norm{F}_{C^m(\Rn)} : F|_E = f\text{ and } F \geq 0 }\,. \label{cm-negnorm}
		\end{equation}
	\end{problem}

	By ``order of magnitude'' we mean the following: Two quantities $ M $ and $ \tilde{M} $ determined by $ E,f,m,n $ are said to have the\underline{ same order of magnitude} provided that $ C^{-1}M \leq \tilde{M} \leq CM $, with $ C $ depending only on $ m $ and $ n $. To compute the order of magnitude of $ \tilde{M} $ is to compute a number $ M $ such that $ M $ and $ \tilde{M} $ have the same order of magnitude. 
	
	Problem \ref{prob.norm} without the nonnegative constraint has been extensively studied, see \cite{BS01,F05,FK09-DataI,FK09-DataII,F09-data-III,FI20}. 
	
	We also consider an open problem posed in \cite{FI20}.
	
	\begin{problem}\label{prob.interpolant}
		Let $ E \subset \Rn $ be a finite set. Let $ f : E \to [0,\infty) $. Compute a nonnegative function $ F \in C^m(\Rn) $ such that $ F|_E = f $ and $ \norm{F}_{C^m(\Rn)} \leq C\norm{f}_{C^m_+(E)} $. 
	\end{problem}

	\medspace
	
	We will present a brief history of Problem  \ref{prob.norm} and an overview of our results on Problems  \ref{prob.norm} and \ref{prob.interpolant}.
	
	We start with elementary background. Given a subset $E \subset \mathbb{R}^n$ and $f: E \rightarrow \mathbb{R}$, we define the \underline{trace norm} of $f$ as
	$$\norm{f}_{C^m(E)} := \inf \set{ \norm{F}_{C^m(\Rn)} : F|_E = f};$$
	we say that $F \in C^m(\mathbb{R}^n)$ is an \underline{almost optimal $C^m(\mathbb{R}^n)$ interpolant} if $F \in C^m(\mathbb{R}^n)$, $F|_E = f$, and $\|F \|_{C^m(\mathbb{R}^n)}\leq C(m,n) \norm{f}_{C^m (E)} $ for some constant $C(m,n)$ depending only on $m, n$. For nonnegative interpolants, one can define analogously the trace norm by requiring the interpolant to be nonnegative, see \eqref{cm-negnorm}.

	We recall the basic finiteness principle of \cite{F05}.

	\namedtheorem[Theorem 0.A]{\label{thm.fp-fef}For large enough $ \ksharp $ and $ \csharp $, both depending only on $ m $ and $ n $, the following holds:
		
		Let $f:E\rightarrow \mathbb{R}$ with $%
		E\subset \mathbb{R}^{n}$ finite. Suppose that for each $S\subset E$ with $%
		\#\left( S\right) \leq k^{\sharp}$ there exists $F^{S}\in C^{m}\left( \mathbb{R}%
		^{n}\right) $ with norm $\left\Vert F^{S}\right\Vert _{C^{m}\left( \mathbb{R}%
			^{n}\right) }\leq 1$, such that $F^{S}=f$ on $S$. Then there exists $F\in
		C^{m}\left( \mathbb{R}^{n}\right) $ with norm $\left\Vert F\right\Vert
		_{C^{m}\left( \mathbb{R}^{n}\right) }\leq C$, such that $F=f$ on $E$.}{Finiteness Principle}

	Theorem \hyperref[thm.fp-fef]{0.A} and several related results were first conjectured by Y.~Brudnyi and P.~Shvartsman\cite{BS94-a,BS94-b,S87}. The first nontrivial case $ C^2(\Rn) $ was proven by P.~Shvartsman\cite{S82,S87} with the sharp finiteness constant $ \ksharp = 3\cdot 2^{n-1} $. Theorem \hyperref[thm.fp-fef]{0.A} is further refined to a Sharp Finiteness Principle in \cite{F09-data-III}, which serves as the backbone for efficient algorithms for computing trace norms and almost optimal interpolants.

	For nonnegative smooth interpolation, in \cite{FIL16-2}, the authors proved the following theorem.
	
	\namedtheorem[Theorem 0.B]{\label{thm.fp-FIL}For large enough $ \ksharp $ and $ \csharp $, both depending only on $ m $ and $ n $, the following holds:

		Let $ f: E \to [0, \infty) $ with $ E \subset \Rn $ finite. Suppose that for each $ S \subset E $ with $ \#(S) \leq \ksharp $, there exists $ F^S \in C^m(\Rn) $ with norm $ \norm{F^S}_{C^m(\Rn)} \leq 1 $, such that $ F^S = f $ on $ S $ and $ F^S \geq 0 $ on $ \Rn $. Then there exists $ F \in C^m(\Rn) $ with norm $ \norm{F}_{C^m(\Rn)} \leq \csharp $, such that $ F = f $ on $ E $ and $ F \geq 0 $ on $ \Rn $.}{Finiteness Principle for Nonnegative Smooth Interpolation}

	The proof of Theorem \hyperref[thm.fp-FIL]{0.B} given in \cite{FIL16-2} depends on a refinement procedure for shape fields proven in \cite{FIL16}. As such, the construction of the interpolant is not very explicit, and the finiteness constant $\ksharp$ is larger than it is necessary. For example, for $m = 2, n = 2$, \cite{FIL16-2} gives $\ksharp \geq 100 + 5^{l_{*}+100}$ for some $l_{*} \geq 100$.

	In this paper, we begin by showing that for $ m= 2, n = 2 $, $ \ksharp = 64 $ is sufficient (see Theorem \ref{thm.2DFP}). Although not proven sharp here, it is a substantial improvement over the one given by \cite{FIL16-2}. 
	
	For a better finiteness constant than \cite{FIL16-2} and also ours, see \cite{Shv08} (which gives $ \ksharp = 8 $); however, the method in \cite{Shv08} assumes the validity of the Finiteness Principle and does not yield a construction for the interpolant. 
	
	With a more careful analysis of our proof for the Finiteness Principle, we are able to prove a Sharp Finiteness Principle analogous to the first one proven in \cite{F09-data-III} without the nonnegative constraint; the Sharp Finiteness Principle reads as follows: Given a finite set $ E \subset \R^2 $ with $ \#(E) = N $, we can produce a list of subsets $ S_1, \cdots, S_L $ such that $ E = \bigcup_{\ell = 1}^L S_\ell $, $ \#(S_\ell) \leq C $, and $ L \leq CN $ such that $ \norm{f}_{C^2_+(E)} $ and $ \max\limits_{\ell = 1, \cdots, L}\norm{f}_{C^2_+(S_\ell)} $ have the same order of magnitude. Thus, computing the order of magnitude of $ \norm{f}_{C^2_+(E)} $ amounts to computing each $ \norm{f}_{C^2_+(S_\ell)} $ for $ \ell = 1, \cdots, L $. In the forthcoming papers \cite{JL20-Ext,JL20-Alg}, we will use this result to provide efficient algorithms analogous to the Fefferman-Klartag algorithms \cite{FK09-DataI} for solving nonnegative interpolation problems. 
	
	Our two-dimensional results in this paper rely on their one-dimensional counterparts. We will provide a detailed analysis of the one-dimensional situation in Section \ref{section: 1D finiteness principle}. Along the way,  we also show the nonexistence of a bounded linear extension operator that preserves nonnegativity. This is the content of Theorem \ref{thm.nonlinear}. This is in sharp contrast to $ C^m(\Rn) $ extensions without the nonnegative constraint, for which there exists a bounded linear extension operator of bounded depth\cite{f-2005-b}.

	Our approach is inspired by \cite{f-2005-b,F09-data-III,F12,I13}. However, we will need new ingredients to apply the machinery adapted from the aforementioned references.

	Lastly, we remark that our approach can be adapted to treat nonnegative $ C^m(\R) $ ($ m > 2 $) extensions for finite sets $E$, and to prove the Finiteness Principle for nonnegative $ C^{1,\omega}(\R^2) $ extensions for arbitrary closed sets $ E $.
	
	\medspace
	
	Next, we sketch the main ideas for our approach, sacrificing accuracy for the ease of understanding.
	
	We begin with interpolation in one-dimension. For nonnegative $ C^2(\R) $ interpolation, we will show that, if one can interpolate three consecutive points, then one can interpolate any finite set of points by patching consecutive three-point interpolants together.\footnote{Here we mention that the finiteness constant $ \ksharp = 3 $ is sharp for nonnegative $ C^2(\R) $ interpolation. See \cite{BS01}.} To handle nonnegative $C^2(\mathbb{R}^2)$ interpolation, we will reduce local interpolation problems to the one-dimensional situation.

	To illustrate the idea, we assume that $E \subset Q_0:=[0,1) \times[0,1)$. For a square $ Q \subset \R^2 $, we write $ 2Q $ to denote the two times concentric dilation of $ Q $, and $ \delta_{Q} $ to denote the sidelength of $ Q $. We perform a Calder\'on-Zygmund decomposition to $ Q_0 $, bisecting $ Q_0 $ and its children, which we will call $ Q_\nu $, until the following conditions are satisfied: Any two nearby squares are comparable in size; $ E \cap 2Q_\nu $ lies on a curve with slope $ \leq C $ and curvature $ \leq C\delta_{Q_\nu}^{-1} $; and any two local solutions near $ Q_\nu $ are indistinguishable up to a Taylor error on the order of $ \delta_{Q_\nu} $. We then solve the local interpolation problem by straightening $ E \cap 2Q_\nu $ and treating it as a one-dimensional problem. To ensure two nearby local solutions are Whitney compatible when patched together by a partition of unity, we prescribe a collection of Whitney-compatible polynomials, denoted by $ P_\nu $, each based at a representative point $ x_\nu $ near the center of $ Q_\nu $, and force the local solution to take $ P_\nu $ as a jet at $ x_\nu $. 
	
	The two-dimensional Finiteness Principle is then a consequence of its one-dimensional counterpart and Helly's Theorem from combinatorial geometry. 
	
	In order to prove the Sharp Finiteness Principle, we need to localize the dependence of the $ P_\nu $'s on the given data $ (E,f) $. This involves a variant of Helly's Theorem, a careful analysis when $ f $ is locally small (on the order of $ \delta_Q^2 $), and the combinatorial properties of the Calder\'on-Zygmund squares.

	Here we have given an overly simplified account of our approach. In practice, we have to control derivatives on small scales and handle subtraction with great care in order to preserve nonnegativity. The technical matters will be handled in the sections below.

	Inspired by \cite{BS01}, we also pose the following question on the best finiteness constant for nonnegative $ C^2(\R^2) $ interpolation, and conjecture the answer to be in the positive.
	
	\begin{problem}
		For nonnegative $ C^2(\R^2) $ interpolation, can we take $ \ksharp = 6 $?
	\end{problem}
	
	It would be interesting to know more about the connection between the methods employed in this paper and the method of ``Lipschitz selection'' presented in \cite{BS01}.
	
	We end the introduction by announcing here our solutions to Problems \ref{prob.norm} and \ref{prob.interpolant}; the detail will be presented in the forthcoming papers \cite{JL20-Ext,JL20-Alg}. For a given $ E \subset \R^2 $ with $ \#(E) = N $, we can process $ E $ with at most $ CN\log N $ operations and $ CN $ storage. After that, we can compute the order of magnitude of $ \norm{f}_{C^2_+(E)} $ for any $ f : E \to [0,\infty) $ using at most $ CN $ operations. After preprocessing $ E $ using at most $ CN\log N $ operations and $ CN $ storage, we are able to receive further inputs, consisting of a function $ f : E \to [0,\infty) $ and a number $ M \geq 0 $. Then, given $ x \in \R^2 $, we are able to produce a list $ (f_\alpha(x) : \abs{\alpha} \leq 2) $ using at most $ C\log N $ operations. Suppose an \textbf{Oracle} tells us that $ \norm{f}_{C^2_+(E)} \leq M $. We can then guarantee the existence of a nonnegative function $ F \in C^2(\R^2) $ with $ \norm{F}_{C^2(\R^2)} \leq CM $ and $ F|_E = f $, such that $ \d^\alpha F(x) = f_\alpha(x) $ for $ \abs{\alpha} \leq 2 $.

	To the extend of our knowledge, there has been no previously known result on Problem \ref{prob.interpolant}.

	This paper is part of a literature on extension and interpolation, going back to the seminal works of H.~Whitney \cite{W34,W34-2,W34-3}. We refer the interested readers to  \cite{BS94-a,BS94-b,BS01,CFIK20,F05,f-2005-b,FK09-DataI,FK09-DataII,F09-data-III,FIL16,FIL16-2,FI20,F12,S82,S87,S20} and references therein for the history and related problems.
	
	\begin{acknowledgment}
		We express our gratitude to Charles Fefferman, Kevin O'Neill, and Pavel Shvartsman for their valuable comments. We also thank all the participants in the 11th Whitney workshop for fruitful discussions, and Trinity College Dublin for hosting the workshop. 
		
		The first author is supported by the UC Davis Summer Graduate Student Researcher Award and the Alice Leung Scholarship in Mathematics. The second author is supported by NSF Grant DMS-1554733 and the UC Davis Chancellor's Fellowship. 
	\end{acknowledgment}

	\section[Main results]{Statement of results}
	
	First we set up notations. Let $ n = 1,2 $. We write $ C^2_+(\Rn) $ to denote the collection of all functions $ F: \Rn \to [0, \infty) $ whose derivatives up to the second order are continuous and bounded. We write $ \d^m $ to denote the $ m $-th derivative of a single-variable function. 
	
	We begin with our results in one-dimension.

	\begin{subtheorem}{theorem}\label{thm. 1D}
		\begin{theorem}[\textbf{1-D Finiteness Principle}]\label{thm.FP-1d}
			There exists a constant $ \csharp > 0 $ such that the following holds. 
			
			Let $E  = \set{ x_1, \dots, x_N} \subset \R$ be a finite set with $ x_1 < \dots < x_N $ and $N \geq 3$. Let $f : E \to [0, \infty)$. Suppose
			\begin{enumerate}[(i)]
				\item For every consecutive three points $ E_j = \set{x_j, x_{j+1}, x_{j+2}} (j = 1, \dots, N-2 ) $ there exists a function $ F_j \in C^2_+(\R) $ such that $ F_j\big|_{E_j} = f $; and
				
				\item $ \norm{F_j}_{C^2(\R)} \leq M $.
			\end{enumerate}
			Then there exists $ F \in C^2_+(\R) $ with
			\begin{enumerate}[(A)]
				\item $ F|_E = f $, and
				
				\item  $ \norm{F}_{C^2(\R)} \leq \csharp M $. 
			\end{enumerate}
		\end{theorem}
		
		\begin{remark}
			In the present work, we do not pursue the minimal $ \csharp $. See, e.g. \cite{F12} for a discussion on best constants.
		\end{remark}

		We will also need the following variant of Theorem \ref{thm.FP-1d} in the proof of Lemma \ref{lem.graph}.

		\begin{theorem}\label{thm.FP-1dpm} 
			There exists a constant $ \csharp > 0 $ such that the following holds. 
			
			Let $E = \{x_1,\dots,x_N\}\subset \mathbb{R}$ be a finite set with $x_1<\dots<x_N$ and $N \geq 3$. Let $f: E \rightarrow \R $. Suppose 
			\begin{enumerate}[(i)]
				\item For every consecutive three points $E_j=\{x_j,x_{j+1},x_{j+2}\}$ $(j=1,\dots,N-2)$, there exists a function $F_j \in C^2(\mathbb{R})$ such that $F_j|_{E_j}=f$;
				\item $\abs{\d^mF_j} \leq A_m$ on $ \R $ for $ m = 0, 1, 2 $.

			\end{enumerate}
			Then there exists $F \in C^2(\R)$ such that 
			\begin{enumerate}[(A)]
				\item $F|_E=f$;
				\item $ \abs{\d^m F} \leq  C A_m $ on $ \R $ for $ m = 0,1,2 $. 
			\end{enumerate}
		\end{theorem}

	\end{subtheorem}

	\begin{remark}
		The proofs of Theorems \ref{thm.FP-1d} and \ref{thm.FP-1dpm} will be given in Section \ref{section: 1D finiteness principle}. 
	\end{remark}

	Let $ n = 1, 2 $. Given a finite set $ E \subset \R^n $, we write $ C^2(E) $ to denote all functions $ f: E \to \R $, equipped with the trace norm $ \norm{f}_{C^2(E)} := \inf\set{\norm{F}_{C^2(\Rn)} : F|_E = f} $. We write $ C^2_+(E) $ to denote all functions $ f : E \to [0,\infty) $, equipped with the ``trace norm'' $ \norm{f}_{C^2_+(E)} := \inf \set{\norm{F}_{C^2(\Rn)} : F|_E = f \text{ and } F \geq 0} $.

	The proofs of Theorems \ref{thm.FP-1d} and \ref{thm.FP-1dpm} along with an argument involving quadratic programming immediately give rise to the following results.

	\begin{subtheorem}{theorem}
		
		\begin{theorem}\label{thm.BD-1d}
			Let $ E \subset \R $ be a finite set. There exist universal constants $ C,D $ and an operator $ \Eb : C^2_+(E) \to C^2_+(\R) $ such that the following hold.
			\begin{enumerate}[(A)]
				\item $ \Eb (f) \big|_E = f $ for all $ f \in C^2_+(E) $.
				\item $ \norm{\Eb (f)}_{C^2(\R)} \leq C\norm{f}_{C^2_+(E)} $.
				\item 

				Moreover, for each $ x \in \R $, there exists $ S(x) \subset E $ with $ \#(S(x)) \leq D $, such that for all $ f,g \in C^2_+(E) $ with $ f|_{S(x)} = g|_{S(x)} $, we have
				\begin{equation*}
				\eqindent
				\d^m (\Eb(f))(x) = \d^m(\Eb (g))(x)
				\enskip
				\text{ for }
				m = 0,1,2\,.
				\end{equation*}

			\end{enumerate}
			
		\end{theorem}
		
		\begin{remark}
			
			In general, $ \Eb $ is not additive. See Theorem \ref{thm.nonlinear}. 
		\end{remark}

		 Theorem \ref{thm.BD-1d} holds in the absence of the nonnegative constraint. This is the content of the next theorem.

		\begin{theorem}\label{thm.BD-1dpm}
			Let $ E \subset \R $ be a finite set. There exist universal constants $ C,D $ and a linear operator $ \Epmb : C^2(E) \to C^2(\R) $ such that the following hold.
			\begin{enumerate}[(A)]
				\item $ \Epmb (f) \big|_E = f $ for all $ f \in C^2(E) $.
				\item $ \norm{\Epmb(f)}_{{C}^2(\R)} \leq C\norm{f}_{{C}^2(E)} $.

				\item Moreover, for each $ x \in \R $, there exists $ S(x) \subset E $ with $ \#(S(x)) \leq D $, such that for all $ f,g \in C^2_+(E) $ with $ f|_{S(x)} = g|_{S(x)} $, we have
				\begin{equation*}
				\d^m (\Epmb(f))(x) = \d^m(\Epmb (g))(x)
				\enskip
				\text{ for }
				m = 0,1,2\,.
				\end{equation*}

			\end{enumerate}
			
		\end{theorem}

	\end{subtheorem}

	\begin{remark}\label{rem.1D-depth}
		The number $ D $ in Theorems \ref{thm.BD-1d} and \ref{thm.BD-1dpm} is called the \underline{depth} of the operator $ \E $. The proofs of Theorems \ref{thm.BD-1d} and \ref{thm.BD-1dpm} will be given in Section \ref{section: 1D finiteness principle}. We also remark that the set $ S(x) $ takes a particularly simple form.
		\begin{itemize}
			\item Suppose $ \#(E) \leq 3 $. We take $ S(x) = E $.
			\item Suppose $ \#(E) \geq 4 $. Enumerate $ E = \set{x_1, \cdots, x_N} $ with $ x_1 < \cdots < x_N $.
			\begin{itemize}
				\item If $ x < x_1 $ or $ x> x_N $, we take $ S(x) $ to be the three points in $ E $ closest to $ x $.
				\item If $ x \in [x_1, x_2] $, we take $ S(x) = \set{x_1, x_2, x_3} $. 
				\item If $ x \in [x_{N-1}, x_N] $, we take $ S(x) = \set{x_{N-2}, x_{N-1}, x_N} $.
				\item Otherwise, we take $ S(x) = \set{x_1', x_2', x_3', x_4'} \subset E $ with $ x_1' < x_2' < x_3' < x_4' $ such that $ x \in [x_2', x_3'] $.

			\end{itemize}
		\end{itemize}

	\end{remark}

	It has been shown in \cite{W34-2} the existence of an extension operator satisfying (A,B) of Theorem \ref{thm.BD-1dpm}. We thank P. Shvartsman for bringing to our attention that an algorithm for constructing $ S(x) $ in a more general one-dimensional setting (without nonnegativity) was given in \cite{S20}, in which the interested readers will also find an informative account of the one-dimensional extension theory (without nonnegativity).

	\begin{theorem}\label{thm.nonlinear}
		Let $ E \subset \R $ be a finite set. There does not exist a map $ \E: C^2_+(E) \to C^2(\R) $ that satisfies both of the following.
		\begin{enumerate}[(A)]
			\item For all $ f \in C^2_+(E) $, we have $ \E(f)(x) = f(x) $ for all $ x \in E $, $ \E(f)\geq 0 $ on $ \R $, and $ \norm{\E(f)}_{C^2(\R)} \leq C\norm{f}_{C^2_+(E)} $ for some universal constant $ C $.

			\item $ \E(f+g) = \E(f) + \E(g) $ for all $ f,g \in C^2_+(E) $.
			
		\end{enumerate}
	\end{theorem}

	\begin{remark}
		By considering finite sets of the form $ E \subset \R^{n-1}\times\set{0} $, we can further conclude that, for $ C^2(\Rn) $ with $ n \geq 2 $, there does not exist a bounded additive extension operator that preserves nonnegativity. See Section \ref{section: 1D finiteness principle} for the proof.
	\end{remark}

	We now turn to our results in two-dimension.

	\begin{theorem}[\textbf{2-D Finiteness Principle}]
		\label{thm.2DFP}
		There exists a constant $ \csharp > 0 $ such that the following holds.
		
		Let $ f : E \to [0, +\infty) $ with $ E \subset \R^2 $ finite. Suppose for each $ S\subset E $ with $ \#(S) \leq 64 $, there exists $ F^S \in C^2_+(\R^2) $ such that
		\begin{enumerate}[(i)]
			\item $ \norm{F^S}_{C^2(\R^2)} \leq M $, and
			\item $ F^S\big|_S = f $.
		\end{enumerate}
		Then there exists $ F \in C^2_+(\R^2) $ such that
		\begin{enumerate}[(A)]
			\item $ F|_E = f $, and
			\item $ \norm{F}_{C^2(\R^2)} \leq \csharp M $.
		\end{enumerate}
	\end{theorem}

	\begin{remark}
		The proof of Theorem \ref{thm.2DFP} is given in Section \ref{section: 2DFPloc}. 
	\end{remark}
	
	We also have an improved version of Theorem \ref{thm.2DFP}.

	\begin{theorem}[\textbf{2-D Sharp Finiteness Principle}]\label{thm.SFP}
		Let $ E \subset \R^2 $ with $ \#(E) = N < \infty $. Then there exist universal constants $ C, C' ,C'' $ and a list of subsets $ S_1, S_2, \cdots, S_L \subset E $ satisfying the following.
		\begin{enumerate}[(A)]
			\item $ \#(S_\ell) \leq C  $ for each $ \ell = 1, \cdots, L $.
			\item $ L \leq C'N $.
			\item Given any $ f : E \to [0,\infty) $, we have
			\begin{equation*}
			\max_{\ell = 1, \cdots, L}\norm{f}_{C^2_+(S_\ell)} \leq \norm{f}_{C^2_+(E)} \leq C'' \max_{\ell = 1, \cdots, L}\norm{f}_{C^2_+(S_\ell)}\,.
			\end{equation*}
		\end{enumerate}
		
	\end{theorem}

	\begin{remark}
		The proof of Theorem \ref{thm.SFP} is given in Section \ref{sect.SFP}.

	\end{remark}

	\medspace

	In a forthcoming paper\cite{JL20-Ext}, we will prove the following result.

	\begin{theorem}\label{thm.BD}
		\newcommand{\ctp}{C^2_+}
		\newcommand{\ct}{C^2}
		\newcommand{\rt}{\mathbb{R}^2}
		\newcommand{\pos}{[0,\infty)}
		Let $ E \subset \rt $ be a finite set. There exist (universal) constants $C, D$, and a map $ \E :   \ctp(E) \times \pos \to \ctp(\rt)  $ such that the following hold. 
		
		\begin{enumerate}[(A)]
			\item Let $ M \geq 0 $. Then for all $ f \in \ctp(E) $ with $ \norm{f}_{\ctp(E)} \leq M $, we have $ \E(f,M) = f $ on $ E $ and 
			$ \norm{\E(f,M)}_{\ct(\rt)} \leq CM $.
			
			\item For each $ x \in \rt $, there exists a set $ S(x) \subset E $ with $ \# (S(x)) \leq D $ such that for all $ M \geq 0 $ and $ f, g \in \ctp(E) $ with $ \norm{f}_{\ctp(E)}, \norm{g}_{\ctp(E)} \leq M $ and $ f|_{S(x)} = g|_{S(x)} $, we have
			\begin{equation*}
			\d^\alpha \E(f,M)(x) = \d^\alpha \E(g,M)(x)
			\text{ for }
			\abs{\alpha} \leq 2\,.
			\end{equation*}
			
		\end{enumerate}
	\end{theorem}

	We will not use Theorem \ref{thm.BD} in this paper.

	As a consequence of Theorems \ref{thm.SFP} and \ref{thm.BD}, in \cite{JL20-Alg}, we will provide the following algorithms.

	\begin{algorithm}\label{alg.norm}
		\begin{center}
			Nonnegative $ C^2(\R^2) $ Interpolation Algorithm - Trace Norm
		\end{center}
		\begin{itemize}
			\item[] \textbf{DATA:} $ E \subset \R^2 $ finite with $ \#(E) = N $. $ f: E \to [0,\infty) $.
			\item[] \textbf{RESULT}: The order of magnitude of $ \norm{f}_{C^2_+(E)} $. More precisely, the algorithm outputs a number $ M \geq 0 $ such that both of the following hold.
			\begin{itemize}
				\item We guarantee the existence of a function $ F \in C^2_+(\R^2) $ such that $ F|_E = f $ and $ \norm{F}_{C^2(\R^2)} \leq CM $.
				\item We guarantee there exists no $ F \in C^2_+(\R^2) $ with norm at most $ M $ satisfying $ F|_E = f $.
			\end{itemize}
		\end{itemize}
			
	\end{algorithm}

	\begin{algorithm}\label{alg.interpolant}
		\begin{center}
			Nonnegative $ C^2(\R^2) $ Interpolation Algorithm - Interpolant
		\end{center}
		\begin{itemize}
			\item[] \textbf{DATA:} $ E \subset \R^2 $ finite with $ \#(E) = N $. $ f: E \to [0,\infty) $. $ M \geq 0 $. 
			\item[] \textbf{ORACLE:} $ \norm{f}_{C^2_+(E)} \leq M $. 
			\item[] \textbf{QUERY:} $ x \in \R^2 $. 
			\item[] \textbf{RESULT}: A list of numbers $ ( f_\alpha(x) : \abs{\alpha} \leq 2) $ that guarantees the following: There exists a function $ \tilde{F} \in C^2_+(\R^2) $ with $ \norm{\tilde{F}}_{C^2(\R^2)} \leq CM$ and $ \tilde{F}|_E = f $, such that $ \d^\alpha \tilde{F}(x) = f_\alpha(x) $ for $ \abs{\alpha} \leq 2 $. 
			
		\end{itemize}
	\end{algorithm}

	\medspace

	We will present the proofs for Theorems \ref{thm. 1D} - \ref{thm.SFP} in the sections below. We will start from scratch and introduce the relevant terminologies and notations in the next section.

	\section[Preliminaries]{Conventions and Preliminaries}
	\label{section: prelim}

	\subsection*{Constants} We use $ c_*, C_*, C, C' > 0 $, etc. to denote ``controlled'' universal constants. They may be different quantities in different instances. We will label them to avoid confusion when necessary.

	\subsection*{Coordinates and norms}
	
	We assume that we are given an ordered orthogonal coordinate system $ x = (s,t)_{\mathrm{standard}} $ on $ \R^2 $ a priori. We write $ B(x,r) $ to denote the open disc of radius $ r > 0 $ centered at $ x \in \R^2 $. 
	
	We use $ \alpha,\beta \in \mathbb{N}_0^2 $ etc. to denote multi-indices. We adopt the partial ordering $ \alpha \leq \beta $ if and only if $ \alpha_i \leq \beta_i $ for $ i = 1,2 $.

	\newcommand{\X}{\mathbb{X}}
	
	Let $ \Omega \subset \Rn $ be a set with nonempty interior $ \Omega^0 $. For positive integers $ m,n $, we write $ C^m(\Omega)$ to denote the vector space of $ m $-times continuously differentiable real-valued functions on $ \Omega^0 $ such that the following norm is finite:
	\begin{equation}\label{def: C2 norm}
	\norm{F}_{C^m(\Omega)} := \sup\limits_{x \in \Omega^0}\left( \sum_{\abs{\alpha}\leq m}  \abs{ \partial^\alpha F(x) }^2 \right)^{1/2}.
	\end{equation}
	We write $ C^m_+(\Omega)$ to denote the collection of functions $F \in C^m(\Omega)$ such that $ F \geq 0 $ on $ \Omega $. This is not a vector space.

	Let $E \subset \Rn $ be finite. We define
	\begin{equation*}
	C^m(E) := \set{F|_E : F \in C^m(\Rn)}\,.
	\end{equation*}
	$ C^m(E) $ is a vector space that can be equipped with a seminorm, which we will called the \underline{trace norm} of $ f \in C^m(E) $:
	\begin{equation*}
	\norm{f}_{C^m(E)} := \inf \set{\norm{F}_{C^m(\Rn)} : F \in C^m(\Rn)\enskip \text{ and } F|_E = f}\,.
	\end{equation*}
	Similarly, we define
	\begin{equation*}
	C^m_+(E) := \set{F|_E : F\in C^m_+(\Rn)}\,.
	\end{equation*}
	We will abuse terminology and refer to the following as the \underline{(nonnegative) trace norm} of $ f \in C^m_+(E) $:
	\begin{equation*}
	\norm{f}_{C^m_+(E)} := \inf \set{\norm{F}_{C^m(\Rn)} : F \in C^m_+(\Rn) \enskip \text{ and } F|_E = f}\,.
	\end{equation*}

	\subsection*{Jets}
	
	\newcommand{\Rj}{\mathcal{R}}
	
	We write $ \P $ to denote the space of degree one polynomials on $\R^2$. It is a three-dimensional vector space.
	
	For $ x_{0} = (s_0, t_0) \in \R^2 $ and a continuously differentiable function $F$ on $\R^2$, the $ 1 $-jet of $F$ at $ x_{0} \in \R^2 $ is given by
	\begin{equation*}
	\jet{x_{0}}{F}(x) := F(x_0) + \nabla F(x_0) \cdot (x - x_0)\, . 
	\end{equation*}
	We use $ \Rj_{x_{0}} $ to denote the vector space of $ 1 $-jets at $ x_0 \in \R^2 $. $ \Rj_{x_{0}} $ inherits a norm from $ \R^3 $ via the identification
	\begin{equation}
	I_{x_0}\, : \, a(s - s_0) + b(t - t_0) + c \, \mapsto \, (a,b,c). 
	\label{eq. jet and euclidean id}
	\end{equation}

	\subsection*{Calder\'on-Zygmund squares}

	A square $ Q \subset \R^2 $ is of the form $ Q = [s_0, s_0 + \delta) \times [t_0, t_0 + \delta) $, where $ \delta > 0 $ and $ s_0, t_0 \in \R $. 
	
	For a square $ Q \subset \R^2 $, $ \lambda Q $ denotes the concentric dilation of $ Q $ by a factor of $ \lambda > 0 $. Let $ Q^* = 2Q $. $\delta_Q$ denotes the side length of $Q$. 
	
	For a square $ Q_0 \in \R^2 $, by a dyadic bisection of $ Q_0 $, we mean dividing $ Q_0 $ into four mutually disjoint congruent squares $ Q_1, Q_2, Q_3, Q_4 $ such that $ Q_0 = \bigcup_{i = 1}^4Q_i $. $ Q_0 $ is called the \underline{dyadic parent} of $ Q_1, \dots, Q_4 $.  In this case, we write $ Q_i^+ = Q_0 $ for $ i = 1, \dots, 4 $. A dyadic parent for a dyadic square is unique if it exists. 
	
	Two squares $ Q $ and $ Q' $ are neighbors if one of the following holds.
	\begin{itemize}
		\item $ Q = Q' $; or
		
		\item $ \cl(Q)\cap \cl(Q') \neq \void $, but $ \intr(Q) \cap \intr(Q') = \void $.
	\end{itemize}
	If $ Q $ and $ Q' $ are neighbors, we write $ Q \leftrightarrow Q' $. 
	
	A collection of mutually disjoint squares $ \Lambda = \set{Q} $ is a \underline{Calder\'on-Zygmund (CZ) covering} of $ \R^2 $ if $ \R^2 = \bigcup_{Q \in \Lambda}Q $, and
	\begin{equation}
	\text{if $ Q \leftrightarrow Q'$, then $ \frac{1}{4}\delta_Q \leq \delta_{Q'} \leq 4 \delta_Q $}\,.
	\label{eq.CZ-def}
	\end{equation}
	It is easy to see that \eqref{eq.CZ-def} implies that a CZ covering satisfies the \underline{bounded intersection property}: If $ Q \in \Lambda $, then
	\begin{equation}\label{def: BIP}
	\#\left( \set{  Q' \in \Lambda  : \frac{9}{8}Q' \cap \frac{9}{8}Q \neq \void   }    \right) \leq 21.
	\end{equation}

	We will only consider \emph{nonnegative} (smooth) cutoff functions and partition of unity. A $C^2$-partition of unity $ \set{\theta_Q} $ subordinate to a CZ covering $ \Lambda = \set{Q} $ of $ \R^2 $ is \underline{CZ-compatible with $ \Lambda $} if
	\begin{equation}\label{def: partition of 1 }
	\theta_Q \geq 0
	, \ 
	\supp{\theta_Q} \subset \frac{9}{8}Q
	, \ 
	\abs{\partial^\alpha\theta_Q} \leq C\delta_Q^{-\abs{\alpha}}
	\ \forall \abs{\alpha} \leq 2
	\text{, and }
	\sum_{Q \in \Lambda} \theta_Q \equiv 1.
	\end{equation}
	Here $ C $ is some universal constant. Such partition of unity exists, see e.g.\cite{W34}.

	\section[Convex sets and polynomials]{Basic convex sets and Whitney fields}
	\label{section: convex}
	
	\begin{definition}\label{def.Gamma}
		Let $E \subset \R^2$ be a finite set. Let $f : E \to [0,+\infty)$. For a point $ x \in \R^2 $, a subset $ S \subset E $, and a real number $ M \geq 0 $, we introduce the following objects:
		\begin{equation}
		\Gamma_+(x, S, M): =  \set{ P \in \P :\, \begin{matrix*}[l]
			\text{There exists }F^S \in C^2_+(\R^2) \text{ such that }\\ \norm{F^S}_{C^2(\R^2)} \leq M, F^S\big|_S = f, \text{ and } \jet{x}{F^S} = P.
			\end{matrix*}  } \,,  
		\label{eq.Gamma+def}
		\end{equation}
		and
		\begin{equation}
		\sigma(x,S): = \set{ P \in\P :\, \begin{matrix*}[l]
			\text{There exist }F^S \in C^2(\R^2)
			\text{ such that } \\
			F^S\big|_S = 0,
			\norm{F^S}_{C^2(\R^2)} \leq 1,
			\text{ and } \jet{x}{F^S} = P.
			\end{matrix*}  }\, .
		\label{eq.sigma-def}
		\end{equation}
		Given an integer $ k \geq 0 $ and a number $ M \geq 0 $, we define
		\begin{equation}
		\G(x,k,M) := \bigcap_{S \subset E,\, \#(S) \leq k}\Gamma_+(x, S, M)\,,
		\label{eq.Gk-def}
		\end{equation}
		and
		\begin{equation}
		\sk(x,k):= \bigcap_{S \subset E,\, \#(S) \leq k} \sigma(x,S)\,.
		\label{eq.sk-def}
		\end{equation}

	\end{definition}

	Since $ \#(E) < \infty $, for sufficiently large $ M \geq 0 $ depending on $ E $ and $ f $, $ \Gamma_+(x,S,M) \neq \void $ for any $ S \subset E $. As a consequence, for a specific $ k $, $ \G(x,k,M) \neq \void $ if $ M $ is sufficiently large. 
	
	It is easy to see that $ \Gamma_+ $, $ \G $, $ \sigma $, and $ \sk $ are convex and bounded (as subsets of $ \R^3 $ via the identification \eqref{eq. jet and euclidean id}). We can easily see from \eqref{eq.sigma-def} and \eqref{eq.sk-def} that $ \sigma $ and $ \sk $ are \emph{symmetric about the origin}. Since $ E $ is finite, for each fixed $ x \in \R^2 $ and $ M > 0 $, there are only finitely many distinct $ \sigma(x,S) $ and $ \Gamma(x,S,M) $ (for a fixed $ M $). Therefore, we may apply the finite version of Helly's Theorem (see Section \ref{subsect:convex} for the statement). Both $ \sk $ and $ \G $ are monotone decreasing (with respect to set inclusion $ \subset $) in $ k $. Furthermore, $ \G $ is monotone increasing in $ M $. 
	
	Since $ \sigma $ and $ \sk $ contain the zero polynomial, they are never empty. 
	
	
	Understanding the shapes of $ \G $ and $ \sk $ is the key to proving Theorems \ref{thm.2DFP}, \ref{thm.SFP}, and \ref{thm.BD}. 
	
	\newcommand{\B}{\mathcal{B}}
	
	We will also be working with the following object.
	
	\begin{definition}\label{def.Taylor ball}
		Given $ x \in \R^2 $ and $ \delta > 0 $, we introduce the following object
		\begin{equation}
		\B(x,\delta) := \set{P \in \P : \abs{\d^\alpha P(x)} \leq \delta^{2-\abs{\alpha}}}.
		\end{equation}
	\end{definition}
	
	To understand the significance of $ \B(x,\delta) $, we point out that Taylor's theorem can be reformulated in the following way: Given $ F \in C^2(\R^2) $ with $ \norm{F}_{C^2(\R^2)} \leq M $, then $ \jet{x}{F} - \jet{y}{F} \in CM \cdot \B(x,\abs{x-y}) $ for any $ x,y \in \R^2 $.

	\subsection{Lemmas on convex sets}\label{subsect:convex}
	
	\begin{lemma}\label{lemma: G - G}
		$ \G(x,k, M) - \G(x,k,M) \subset 2M\cdot\sk(x,k)$. The minus sign denotes vector subtraction. 
	\end{lemma}
	
	\begin{proof}
		
		Let $ P_1,P_2 \in \G(x,k,M) $. For each $ S \subset E $ with $ \#(S) \leq k $, there exist $ F_1^S, F_2^S \in C^2_+(\R^2) $ such that for $ i = 1, 2 $, $ F_i^S\big|_S= f $, $ \norm{F_i^S}_{C^2(\R^2)} \leq M $, and $ \jet{x}{F_i^S} = P_i $.
		Then $ (F_1^S - F_2^S)\big|_{S} = 0$, $ \norm{F_1^S - F_2^S}_{C^2(\R^2)} \leq 2M $, and $ \jet{x}{F_1^S - F_2^S} = P_1 - P_2 $. Since $ S $ is arbitrary, $ P_1 - P_2 \in \sk(x,k, 2M) = 2M\cdot\sk(x,k) $. 
	\end{proof}

	We recall a classical result by Helly, the proof of which can be found in \cite{rF97}.
	
	\namedtheorem[Helly's Theorem]{\label{thm.Helly}Let $ \mathcal{F} $ be a finite collection of convex sets in $ \R^D $. Suppose every subcollection of $ \mathcal{F} $ of cardinality at most $ (D+1) $ has nonempty intersection. Then the whole collection has nonempty intersection.}{}

	The following lemma states that we can control polynomials in $ \G $ based at some point by polynomials that are based at a different point but are ``less universal'' (in the sense that it is the jet for an interpolant for fewer points).
	
	\begin{lemma}\label{lemma: Gamma trade off}
		There exists a universal constant $ C $ such that the following holds. Let $ x, x' \in \R^2 $. Let $ k_1 \geq 4k_2 $. Let $ M \geq 0 $.
		Given $ P \in \G(x,k_1,M) $, there exists $ P' \in \G(x', k_2, M) $ such that
		\begin{equation*}
		\abs{\partial^\alpha(P - P')(x)}\enskip,\enskip \abs{\partial^\alpha(P - P')(x')} \leq CM\abs{x - x'}^{2 - \abs{\alpha}}
		\text{ for } \abs{\alpha} \leq 1.
		\end{equation*}
	\end{lemma}
	
	\begin{proof}
		\newcommand{\gtemp}{\Gamma_+^{\mathrm{temp}}}
		\newcommand{\gtems}{\Gamma_+^{\mathrm{temp}, \#}}
		
		Fix $ P $ and $ M $ as in the hypothesis of the lemma. For each $ S \subset E $, we define
		\begin{equation*}
		\gtemp(S) := \set{ P' \in \P : \,\, \begin{matrix*}[l]
			\text{There exists } F^S \in C^2_+(\R^2) \text{ such that } \norm{F^S}_{C^2(\R^2)} \leq M,\\
			F^S\big|_S = f, \jet{x}{F^S} = P, \text{ and } \jet{x'}{F^S} = P'.
			\end{matrix*}  } \, .
		\end{equation*}
		Then $ \gtemp $ is a convex and bounded subset of $ \P $. Notice that
		\begin{equation}
		S \subset \widetilde{S} \text{ implies } \gtemp(\widetilde{S}) \subset \gtemp(S).
		\label{eq.gt-nest}
		\end{equation}
		It also follows from the definition of $ \G(x,k_1,M) $ that
		\begin{equation}
		\text{ if }  \#(S) \leq k_1, \text{ then } \gtemp(S) \neq \void\,.
		\label{eq.gt-notvoid}
		\end{equation}
		Let $ S_1, ..., S_4 \subset E $ be given with $ \#(S_i) \leq k_2 $ for each $ i $. Let $ S = \bigcup_{i = 1}^4S_i $. Then $ \#(S) \leq 4k_2 \leq k_1 $. Thanks to \eqref{eq.gt-notvoid}, $ \gtemp(S) \neq \void $. 
		Since $ S_i \subset S $, \eqref{eq.gt-nest} implies that $ \gtemp(S) \subset \gtemp(S_i) $. Therefore,
		\begin{equation*}
		\bigcap_{i = 1}^4 \gtemp(S_i) \supset \gtemp(S)  \neq \void.
		\end{equation*}
		Since $ \set{S_i}_{i = 1}^4 $ are arbitrary, applying \hyperref[thm.Helly]{Helly's Theorem} to the convex sets $ \gtemp(S_i)\subset \P $ (with $ \dim \P = 3 $), we have
		\begin{equation*}
		\bigcap_{S \subset E, \#(S) \leq k_2}\gtemp(S) \neq \void.
		\end{equation*}
		Let $ P' \in \bigcap\limits_{S \subset E, \#(S) \leq k_2}\gtemp(S) $. By definition, $ P' \in \G(x', k_2, M) $. Setting $ S = \void $, we see that there exists $ F \in C^2_+(\R^2) $ with
		\begin{itemize}
			\item $ \norm{F}_{C^2(\R^2)} \leq M $; and
			\item $ \jet{x}{F} = P $ and $ \jet{x'}{F} = P' $. 
		\end{itemize}
		By Taylor's theorem, we have
		\begin{equation*}
		\abs{\partial^\alpha(P - P')(x)}
		= \abs{\partial^\alpha(\jet{x}{F} - \jet{x'}{F})(x)}
		\leq \abs{\partial^\alpha (F - \jet{x'}{F})(x)}
		\leq CM\abs{x - x'}^{2 - \abs{\alpha}}\,.
		\end{equation*}
		The estimate for $ 	\abs{\partial^\alpha(P - P')(x')} $ is similar. 
	\end{proof}

	\begin{lemma}\label{lemma: 17}
		Under the hypothesis of Theorem \ref{thm.2DFP}, $ \G(x, 16, M) \neq \void $ for all $ x \in \R^2 $. 
	\end{lemma}
	
	\begin{proof}
		Recall that $ \Gamma_+(\cdot, \cdot, \cdot) $ is a convex set in a three-dimensional vector space $ \P $. By \hyperref[thm.Helly]{Helly's Theorem}, it suffices to show that the intersection of any four-element subfamily is nonempty. To this end, fix $ x \in \R^2 $, let $ S_1, \cdots, S_4 \subset E $ with $ \#(S_i) \leq 16 $, and let $ S = \bigcup_{i = 1}^4S_i $. We have 
		\begin{equation}
		\Gamma_+(x, S, M) \subset \bigcap_{i = 1}^4 \Gamma_+(x, S_i, M).
		\label{eq. 17}
		\end{equation}
		Since $ \#(S) \leq 64 $, the hypothesis of Theorem \ref{thm.2DFP} implies that $ \Gamma_+(x, S, M) \neq \void $, and hence, the intersection on the right hand side of \eqref{eq. 17} is nonempty. This concludes the proof.
	\end{proof}

	The following variant of \hyperref[thm.Helly]{Helly's Theorem} can be found in Section 3 of \cite{f-2005-b}.

	\begin{lemma}\label{lem.cHel}
		Let $ \mathcal{F} $ be a finite collection of compact, convex, and symmetric subsets of $ \R^D $. Suppose $ 0 $ is an interior point for each $ K \in \mathcal{F} $. Then there exist $ K_1, \cdots, K_{D(D+1)} \in \mathcal{F} $ such that
		\begin{equation*}
		K_1 \cap \cdots \cap K_{D(D+1)} \subset C_D \cdot \brac{\bigcap\limits_{K \in \mathcal{F}} K } \,.
		\end{equation*}
		Here, $ C_D $ is a constant that depends only on $ D $.
	\end{lemma}

	\begin{lemma}\label{lem.sigma-localized}
		There exists a universal constant $ C $ such that the following holds. Let $ x \in \R^2 $. Then given $ k \geq 0 $, there exist $ S_1, \cdots, S_{12} \subset E $, with $ \#(S_i) \leq k $ for each $ i $, such that 
		\begin{equation*}
		\bigcap_{i = 1}^{12}\sigma(x,S_i) \subset C \cdot \brac{ \bigcap_{S \subset E, \#(S) \leq k} \sigma(x,S) } = C \cdot \sk(x,k)\,.
		\end{equation*}
		
	\end{lemma}

	\begin{proof}
		Let $ x \in \R^2 $. Note that $ \sk(x,k) $ has nonempty interior (in the relative topology of the maximal affine space that it spans).
		We apply Lemma \ref{lem.cHel} (with $ D \leq \dim \P = 3 $) to $ \cl(\sigma(x,S)) $. Thus, there exist $ S_1, \cdots, S_{12} \subset E $ with $ \#(S_i) \leq k $ for each $ i = 1, \cdots, 12 $, such that
		\begin{equation*}
		\bigcap_{i = 1}^{12}\sigma(x,S_i) \subset C_D \cdot \brac{ \bigcap_{S \subset E, \#(S) \leq k} \cl\brac{\sigma(x,S) }}\,.
		\end{equation*}
		Therefore, 
		\begin{equation}
		\bigcap_{i = 1}^{12}\sigma(x,S_i) \subset 2C_D \cdot \brac{ \bigcap_{S \subset E, \#(S) \leq k} \sigma(x,S) } = 2C_D \cdot \sk(x,k)\,.
		\label{eq.polarrep}
		\end{equation}
		This proves the lemma. 
	\end{proof}

		\subsection{Whitney fields}
		\label{section.WF}
	
	\newcommand{\cone}{\mathcal{C}_+}
	
	In this subsection, we assume $ n = 1 $ or $ 2 $. We use $ \P $ to denote the space of polynomials on $ \Rn $ with degree no greater than one.
	
	We now recall the notion of a Whitney field.

	Let $ S \subset \Rn $ be a finite set. We use $ W^2(S) $ to denote the (finite dimensional) vector space of sections of $ S \times \P $. An element $ \vec{P} \in W^2(S) $ is called a \underline{Whitney field}, and has the form $ \vec{P} = (P^x)_{x \in S} $. $ W^2(S) $ can be endowed with a norm
	\begin{equation}
	\norm{\vec{P}}_{W^2(S)} := \max_{x \in S}\brac{ \sum_{\abs{\alpha} \leq 1}\abs{\d^\alpha P^x(x)}^2 }^{1/2} + \max_{\substack{x,y\in S\\x \neq y}}\,\brac{ \sum_{\abs{\alpha} \leq 1} \brac{ \frac{\abs{\d^\alpha (P^x - P^y)(x)}}{\abs{x-y}^{2-\abs{\alpha}}} }^2  }^{1/2}\,.  
	\label{eq.W2-norm}
	\end{equation}
	
	We are interested in jets that can be extended to nonnegative $ C^2 $ functions. For $ x \in \Rn $ and $ M \geq 0 $, we define
	\begin{equation}
	\cone(x,M) : = \set{ P \in \P : \brac{\sum_{\abs{\alpha}\leq 1}\abs{\d^\alpha P(x)}^2}^{1/2} \leq M,\, P(x) \geq 0,\, 
		\text{ and } 
		\abs{\nabla P} \leq \sqrt{4M\cdot P(x)} }\,.
	\label{eq.cone-def}
	\end{equation}

	The next lemma tells us how to approximate $ \Gamma_+ $.

	\begin{lemma}\label{lem.cone-equiv}
		There exists a universal constant $ C $ such that the following holds. Given $ M \geq 0 $, we have 
		\begin{equation}
		\Gamma_+(x,\void,C^{-1}M) \subset \cone(x,M) \subset \Gamma_+(x,\void,CM)\,.	\footnote{Here, when $ n = 1 $, $ \Gamma_+(x,\void,M) $ is defined to be $ \set{\bar{\mathscr{J}}_x F : F \in C^2(\R), F \geq 0, \text{ and } \norm{F}_{C^2(\R)} \leq M} $, where $ \bar{\mathscr{J}}_x $ is the first degree Taylor expansion about the point $ x $ for a single-variable function. }
		\label{eq.cone-equiv}
		\end{equation}
	
	\end{lemma}

	\begin{proof}
		The statement is clear for $ M = 0 $. 
		
		Suppose $ M > 0 $.
		
		The first inclusion follows immediately from Taylor's theorem. We prove the second inclusion.
		
		Without loss of generality, we may assume $ x = 0 $. 
		
		Pick $ P \in \cone(0,M) $. We have
		\begin{equation*}
		\sum_{\abs{\alpha}\leq 1}\abs{\d^\alpha P(0)}^2 \leq M^2
		\end{equation*}
		and
		\begin{equation}
		\abs{\nabla P}^2 \leq 4M\cdot P(0).
		\label{eq.cone-equiv-tc1}
		\end{equation}
		Restricting $ P $ to each one-dimensional subspace of $ \Rn $ and using \eqref{eq.cone-equiv-tc1}, we see that
		\begin{equation*}
		\tilde{P} := M\abs{x}^2 + P(x) = M\abs{x}^2 + \nabla P \cdot x + P(0) \geq 0
		\text{ for } x \in \Rn. 
		\end{equation*}
		Let $ B $ be the unit disc in $ \Rn $. Let $ \theta \in C^2_+(\Rn) $ be a cutoff function satisfying
		\begin{equation*}
		\supp{\theta} \subset B\,,\,
		\theta \equiv 1\text{ near } 0\,,
		\,
		\abs{\d^\alpha \theta} \leq C
		\text{ for }\abs{\alpha} \leq 2.
		\label{eq.chiest}
		\end{equation*}
		We define
		\begin{equation*}
		F := \theta \cdot \tilde{P} = M\theta \abs{x}^2 + \theta \nabla P \cdot x + \theta P(0)\,.
		\end{equation*}
		Immediately, we have $ \jet{0}{F} = \jet{0}{\tilde{P}} = \jet{0}{P} = P $ and $ F \geq 0 $ on $ \Rn $. Moreover,
		\begin{equation}
		\abs{\d^\alpha F(x)} \leq CM
		\text{ for } x \in B
		\text{ and }\abs{\alpha} \leq 2\,.
		\label{eq.rLrL}
		\end{equation}
		Since $ \theta $ is supported in $ B $, we can conclude that, $ \norm{F}_{C^2(\Rn)} \leq CM $ and $ \jet{0}{F} \in \Gamma_+(0,\void,A) $ for $ A = CM $. This concludes the proof.	
		
	\end{proof}

	\newcommand{\rmin}{\mathcal{M}}

	\begin{definition}
		Recall the definition of $ \cone $ in \eqref{eq.cone-def}. Given a finite set $ S \subset \Rn $, we define 
		\begin{equation}
		W^2_+(S) := \set{\vec{P} = (P^x)_{x \in S} \in W^2(S) : \begin{matrix}
			\text{ There exists } M \geq 0
			\text{ such that  }\\P^x \in \cone(x,M) \text{ for each } x \in S.
			\end{matrix}
		}\,.
		\label{eq.def-W2+S}
		\end{equation}
		For $ \vec{P} \in W^2_+(S) $, we define
		\begin{equation}
		\norm{\vec{P}}_{W^2_+(S)} := \norm{\vec{P}}_{W^2(S)} + \rmin(\vec{P}),\,
		\end{equation}
		where $ \norm{\vec{P}}_{W^2(S)} $ is defined in \eqref{eq.W2-norm} and
		\begin{equation}
		\rmin(\vec{P}) := \inf\set{ M > 0  :  	\abs{\nabla P^x} \leq \sqrt{4M P^x(x)} \text{ for each } x \in S}\,.
		\label{eq.rhomin}
		\end{equation}
	\end{definition}

	\begin{remark}
		The definition of $ \rmin $ is motivated by the estimate \eqref{eq.rLrL}.
	\end{remark}

	The following is immediate from Taylor's theorem and Lemma \ref{lem.cone-equiv}.
	
	\begin{lemma}\label{lem.Taylor-WF}
		Let $ F \in C^2_+(\Rn) $. Let $ S \subset \Rn $ be a finite set. For each $ x \in S $, let $ P^x := \jet{x}{F} $. Let $ \vec{P} := (P^x)_{x \in S} $. Then $ \vec{P} \in W^2_+(S) $ with $ \norm{\vec{P}}_{W^2_+(S)} \leq C\norm{F}_{C^2(\Rn)} $ for some constant $ C $ depending only on $ n $. 
	\end{lemma}

	The next lemma follows immediately from Lemma \ref{lem.Taylor-WF}.
	
	\begin{lemma}\label{lem.TC-W}
		Let $ S \subset \Rn $ be a finite set. Given any $ f \in C^2_+(E) $, there exists $ \vec{P} \in W^2_+(S) $ such that $ \norm{\vec{P}}_{W^2_+(S)} \leq C \norm{f}_{C^2_+(S)} $ and $ P^x(x) = f(x) $ for each $ x \in S $. The constant $ C $ depends only on $ n $.
	\end{lemma}

	\newcommand{\ew}{\mathcal{W}}
	\begin{lemma}[Whitney extension theorem for finite sets]\label{lem.WE-map}
		Let $ S \subset \Rn $ be a finite set. There exist a constant $ C $ depending only on $ n $, and a map $ \ew^S : W^2_+(S) \to C^2_+(\Rn) $ such that the following hold.
		\begin{enumerate}[(A)]
			\item $ \norm{\ew^S(\vec{P})}_{C^2(\Rn)} \leq C \norm{\vec{P}}_{W^2_+(S)} $.
			
			\item $ \jet{x}{\ew^S(\vec{P})} = P^x $ for each $ x \in S $.
		\end{enumerate} 
	\end{lemma}

	\begin{proof}[Sketch of proof]
		We begin by assuming $ S = \set{y} $. We write $ * $ instead of $ \set{*} $ in certain places to avoid cumbersome notation.
		
		Let $ \vec{P} = P \in W^2_+(y) $. 
		
		Suppose $ P (y) = 0 $. Since $ P \in W^2(y) $, we must have $ \nabla P \equiv 0 $. Therefore, we simply set
		\begin{equation*}
		\ew^y (P) \equiv 0\,.
		\end{equation*}
		Conclusions (A) and (B) are satisfied. 
		
		Suppose $ P(y) > 0 $. 
		By definition,
		\begin{equation*}
		P \in \cone(y,M),
		\text{ where } M := \max\set{  \brac{\sum_{\abs{\alpha}\leq 1}\abs{\d^\alpha P(y)}^2}^{1/2}\enskip,\enskip \frac{\abs{\nabla P}^2}{4{P(y)}} }
		\,.
		\end{equation*}
		Thus, $ \tilde{P}(x) := P(x) + M\abs{x-y}^2 \geq 0 $ for all $ x \in \Rn $.

		Let $ \chi $ be a cutoff function that satisfies $ \chi \equiv 1 $ near $ y $, $ \supp{\chi} \subset B(y,1) $, and $ \abs{\d^\alpha \chi} \leq C $ for $ \abs{\alpha} \leq 2 $. Define
		\begin{equation}
		\ew^{{y}}(P) := \chi \cdot \tilde{P}.
		\label{eq.WyP}
		\end{equation}
		It is clear that $ \ew^{{y}}(P) \geq 0 $ and $ \jet{y}{\ew^{y}({P})} = \jet{y}{\tilde{P}} = P $. 
		Moreover, for $ x \in B(y,1) $ and $ \abs{\alpha} \leq 2 $
		\begin{equation*}
		\abs{\d^\alpha\ew^{{y}}({P})} \leq CM. 
		\end{equation*}
		Therefore,
		\begin{equation}
		\norm{\ew^{{y}}({P})}_{C^2(\Rn)} \leq CM \leq C'\norm{{P}}_{W^2_+(S)}.
		\label{eq.wmap-est-2}
		\end{equation}

		Next, we sketch the proof of the lemma for general $ S $,
		
		Let $ WC $ be a Whitney cover of $ \Rn $ associated with the set $ S $, and let $ \set{\theta_Q} $ be a partition of unity compatible with $ WC $. See \cite{FIL16-2}.
		
		In particular, $ WC $ and $ \set{\theta_Q} $ satisfy the following properties.
		\begin{itemize}
			\item $ \Rn = \bigcup_{Q \in WC}Q $;
			\item $ Q \in WC $ if and only if $ Q $ satisfies one of the following:
			\begin{itemize}
				\item $ \delta_Q = 1 $ and $ S \cap Q^* \leq 1 $ (recall that $ Q^* = 2Q $);
				\item $ \delta_Q < 1 $, $ S \cap Q^* \leq 1 $, and $ S \cap (Q^+)^* > 1 $ (recall that $ Q^+ $ is the dyadic parent of $ Q $). 
			\end{itemize}
			\item If $ Q, Q' \in WC $ and $ Q \leftrightarrow Q' $ (i.e. the closures of $ Q $ and $ Q' $ have nonempty intersection), then $ C^{-1}\delta_Q \leq \delta_{Q'} \leq C\delta_Q $. 
			\item $ \sum_{Q \in WC}\theta_Q \equiv 1 $,
			\item $ \supp{\theta_Q} \in Q^* $ for each $ Q \in WC $, and
			\item $ \abs{\d^\alpha \theta_Q} \leq C\delta_Q^{-\abs{\alpha}} $ for $ \abs{\alpha} \leq 2 $ and $ Q \in WC $.
		\end{itemize}
		For each $ Q \in WC $, we consider three different cases.
		\begin{enumerate}[\text{Case} 1]
			\item When $ S \cap Q^* \neq \void $, we set $ \ew^Q := \ew^{{y}} $ where $ y \in S \cap Q^* $ and $ \ew^y $ is defined in \eqref{eq.WyP}. We set $ P^Q := P^y $.
			\item When $ S \cap Q^* = \void $ and $ \delta_Q < 1 $, we may pick $ y \in S \cap (Q^+)^* $. We set $ \ew^Q := \ew^{{y}} $ and set $ P^Q := P^y $.
			\item When $ S \cap Q^* = \void $ and $ \delta_Q = 1 $, we set $ \ew^Q \equiv 0 $ and $ P^Q :\equiv 0 $. 
		\end{enumerate}
		Finally, we set
		\begin{equation*}
		\ew^S(\vec{P}) := \sum_{Q \in WC}\theta_Q \cdot \ew^Q(P^Q).
		\end{equation*}
		One then verifies that $ \ew^S(\vec{P}) \geq 0 $ and $ \norm{\ew^S(\vec{P})}_{C^2(\Rn)} \leq C \norm{\vec{P}}_{W^2_+(S)} $ via Lemma \ref{lem.cone-equiv} and a routine argument from the classical Whitney extension theorem. See \cite{Ste-70} for details.
	\end{proof}

	\section[CZ decomposition]{Calder\'on-Zygmund squares} 
	\label{sect.CZ}

	\subsection{Calder\'on-Zygmund decomposition of $ \mathbb{R}^2 $}

	\begin{definition}
		Let $ \cnice >  0 $ and $ k \geq 1 $. Recall the notation $ Q^* = 2Q $. We say a dyadic square $ Q $ is \underline{$ k $-nice} if for all $ x \in E \cap Q^* $,
		\begin{equation}\label{def: k nice}
		\diam{\sk(x,k)} \geq \cnice\delta_Q.
		\end{equation}
	\end{definition}
	
	We now describe our decomposition procedure. 
	
	\begin{CZalg}
		Let $ Q $ be a square.
		\begin{itemize}
			\item If $ Q $ is $ k $-nice, then return $ \Lambda_Q^{(k)} = \set{Q} $;
			
			\item otherwise, return
			\begin{equation*}
			\Lambda_Q^{(k)} := \bigcup \set{ \Lambda_{Q'}^{(k)} : Q' \text{ dyadic and } (Q')^+ = Q  }.
			\end{equation*}
		\end{itemize}
	\end{CZalg}

	\begin{remark}
		The algorithm terminates after finitely many steps for each unit square. To see this, notice that $ E $ is finite, and for fixed $ k $ and $ \cnice $, (\ref{def: k nice}) clearly holds for sufficiently small squares containing no more than one point. Moreover, since $ \sk $ does not depend on $ f $, the complexity of our algorithm depends solely on the set $ E $. 
	\end{remark}

	\begin{definition}
		For a particular choice of $ \cnice > 0 $ and $ k \geq 1 $, we use $ \knice =\set{Q_i} $ to denote the collection of $ k $-nice squares obtained from applying the algorithm above to each of the unit squares with their vertices on the integer lattice.
	\end{definition}

	\begin{lemma}\label{lemma: CZ}
		$\knice$ is a CZ covering of $ \R^2 $.
	\end{lemma}

	\begin{proof}
		Since we obtain $ \knice $ by applying the algorithm to each square of the unit grid, $ \knice $ is indeed a covering of $ \R^2 $.
		
		Suppose \eqref{eq.CZ-def} fails, i.e., there exist some $ Q, Q' \in \knice $ with $ Q \leftrightarrow Q' $ but 
		\begin{equation*}
		\delta_Q \leq \frac{1}{8} \delta_{Q'}.
		\end{equation*}
		Then $ (Q^+)^* \subset (Q')^* $. Since $ Q^+ $ is not $ k $-nice, there exists $ \hat{x} \in E\cap (Q^+)^* \setminus Q^* $ such that 
		\begin{equation*}
		\diam{\sk(\hat{x},k)} < 2\cnice\delta_Q.
		\end{equation*}
		On the other hand, 
		\begin{equation*}
		\cnice\delta_{Q'} \leq  \diam{\sk(\hat{x},k)}.
		\end{equation*}
		A contradiction is reached once we combine all the inequalities above, because $ Q' $ is $ k $-nice.
	\end{proof}

	Our main goal is to construct a local interpolant for each $k$-nice square and then to patch these local solutions together. We need several lemmas that guarantee the consistency of our operation. 
	
	The following lemma states that polynomials in $ \G $ with the same base point $ x $ control each other in the Whitney sense after our decomposition.

	\begin{lemma}\label{lemma: estimate for difference in Gamma}
		Let $ \cnice, k \geq 1 $, $ Q \in \knice $, $ x \in E \cap Q^* $, and $ 0 \leq \abs{\alpha} \leq 1 $. If $ P, P' \in \G(x,k, M) $, then
		\begin{equation}\label{eq. G-G est}
		\abs{\partial^\alpha (P - P')(x)} \leq 14\cnice M \delta_Q^{2-\abs{\alpha}}.
		\end{equation}
		
	\end{lemma}

	\begin{proof}
		
		Note that \eqref{eq. G-G est} is immediate if $ \delta_Q = 1 $ or $ \alpha = (0,0) $. Therefore, we only need to consider the case when $ \delta_Q < 1 $ and $ \abs{\alpha} = 1 $. The assumption $ \delta_Q < 1 $ implies that there exists $ y \in E \cap (Q^+)^* $ such that $ \diam{\sk(y, k)} < 2\cnice\delta_Q $. Fix such $ y $. 
		
		Suppose toward a contradiction, that we can find a point $ x \in E \cap Q^* $ and $ P, P' \in \G(x,k, M) $ such that (\ref{eq. G-G est}) is false for some $ \abs{\alpha} = 1 $. Fix such $ \alpha $. 
		
		By Lemma \ref{lemma: G - G}, $ P - P' \in 2M\cdot\sk(x, k) $. By definition, for any $ S \subset E $ with $ \#(S) \leq k $, there exists $ F^{S} \in C^2(\R^2) $ such that
		\begin{itemize}
			\item $ F^S\big|_S = 0 $,
			\item $ \norm{F^{S}}_{C^2(\R^2)} \leq 2M $, and
			\item $ \partial^\alpha(\jet{x}{F^{S}}) = \partial^\alpha(P - P') $.
		\end{itemize}
		By assumption, $ \abs{\partial^\alpha F^{S}(x)} > 14\cnice M\delta_Q $. Since, $ x, y \in (Q^+)^* $, we have $ \abs{x - y} < 6\delta_Q $. Therefore,
		\begin{equation*}
		\begin{split}
		\abs{\partial^\alpha\jet{y}{F^S}(y)}
		=\abs{\partial^\alpha F^{S}(y)}
		\geq \abs{\partial^\alpha F^{S}(x)} - \norm{F^{S}}_{C^2(\R^2)}\abs{x - y}
		\geq 2\cnice M\delta_Q.
		\end{split}
		\end{equation*}
		Since $ S $ is arbitrary, we have $ \diam{\sk(y,k)} \geq 2\cnice\delta_Q$. A contradiction.
	\end{proof}

	\begin{lemma}\label{lemma: neighboring polynomials}
		Let $ \cnice, k \geq 1 $. There exists a universal constant $ C $ such that the following holds. Let $ Q, Q' \in \knice $. Let $ x_Q \in Q $ and $ x_{Q'} \in Q' $. Let $ M \geq 0 $. Let $ P_Q \in \G(x_Q, 4k,M) $ and $ P_{Q'} \in \G(x_{Q'},4k,M) $. Then for $ \abs{\alpha} \leq 1 $ and $ x \in 100Q \cup 100Q' $, 
		\begin{equation}\label{eq: neighboring polynomials}
		\abs{\partial^\alpha(P_Q - P_{Q'})(x)} \leq CM\cdot \max\set{\abs{x_Q - x_{Q'}}, \delta_Q, \delta_{Q'}}^{2 - \abs{\alpha}}.
		\end{equation}
	\end{lemma}
	
	\newcommand{\ptemp}{P_{\mathrm{temp}}}
	\newcommand{\seg}{\mathrm{seg}}
	
	\begin{proof}
		\newcommand{\dm}{\delta_{\infty}}
		Set
		\begin{equation*}
		\dm := \max\set{\abs{x_Q - x_{Q'}}, \delta_Q, \delta_{Q'}}.
		\end{equation*}
		By \eqref{eq.CZ-def}, we have
		\begin{equation*}
		\abs{x_Q - x}, \abs{x_{Q'}- x}, \abs{x_Q - x_{Q'}} \leq C\dm
		\text{ for }
		x \in 100Q \cup 100Q'.
		\end{equation*}
		
		By Lemma \ref{lemma: Gamma trade off}, there exists a $ \ptemp \in \G(x_{Q'}, k, M) $ with
		\begin{equation}\label{eq: pi - ptemp}
		\abs{\partial^\alpha(P_Q - \ptemp)(x_{Q'})} \leq CM\abs{x_Q - x_{Q'}}^{2 - \abs{\alpha}} \leq CM\dm^{2 - \abs{\alpha}}. 
		\end{equation} 
		
		Since $ P_{Q'} \in \G(x_{Q'}, 4k, M) \subset \G(x_{Q'}, k, M) $, Lemma \ref{lemma: estimate for difference in Gamma} applied to $ P_{Q'} $ and $ \ptemp $ gives
		\begin{equation}\label{eq: pj - ptemp}
		\abs{\partial^\alpha(P_{Q'} - \ptemp)(x_{Q'})} \leq CM\delta_{Q'}^{2 - \abs{\alpha}} \leq CM\dm^{2 - \abs{\alpha}}
		\text{ for }\abs{\alpha} \leq 1.
		\end{equation}
		
		Combining (\ref{eq: pi - ptemp}) and (\ref{eq: pj - ptemp}), we have
		\begin{equation}\label{eq: Pi - Pj}
		\abs{\partial^\alpha(P_Q- P_{Q'})(x_{Q'})} \leq CM\dm^{2 - \abs{\alpha}}
		\text{ for }\abs{\alpha} \leq 1.
		\end{equation}
		
		Since $ P_Q$ and $ P_{Q'} $ are affine polynomials, (\ref{eq: neighboring polynomials}) follows from (\ref{eq: Pi - Pj}) in the case $ \abs{\alpha} = 1 $. By the fundamental theorem of calculus, we have
		\begin{equation}\label{eq: FTC}
		(P_Q- P_{Q'})(x) = (P_Q - P_{Q'})(x_{Q'}) + \int_{ \seg(x_{Q'} \to x)}\nabla(P_Q - P_{Q'})\,,
		\end{equation}
		where $ \seg(x_{Q'} \to x) $ is the straight line segment from $ x_{Q'} $ to $ x $. Note that $ \nabla(P_Q - P_{Q'}) $ is a constant vector since both $ P_Q $ and $ P_{Q'} $ are affine. Taking the absolute value of (\ref{eq: FTC}) and applying \eqref{eq: Pi - Pj} with $ \abs{\alpha} = 1 $, we conclude that (\ref{eq: neighboring polynomials}) holds for $ \abs{\alpha} = 0 $.
	\end{proof}

	\subsection{Local geometry}
	
	\newcommand{\gr}{\mathrm{Graph}}
	
	The goal of this section is to show that according to our decomposition, we have partitioned the data points into clusters whose geometry is essentially one-dimensional. To proceed, we introduce some notations.
	
	Note that the $ C^2 $ norm we are using in \eqref{def: C2 norm} is rotationally invariant. Let $ \omega \in [-\pi/2, \pi/2] $. We associate with $ \omega $ a coordinate system obtained by rotating the plane counterclockwise about the origin by an angle of $ \omega $. Thus, for $ x \in \R^2 $,
	\begin{equation*}
	x = (s,t)_{\mathrm{standard}} = (x^{(1)}_\omega, x^{(2)}_\omega)_\omega,
	\end{equation*}
	where $ x^{(1)}_\omega = s\cos\omega + t \sin \omega $ and $ x^{(2)}_\omega = -s\sin\omega + t\cos\omega $. When the choice of $ \omega $ is clear, we write $ \partial_1, \partial_2 $ to denote the partial derivatives with respect to the first, second variable, respectively. They coincide with the directional derivatives along $ \omega $ and $ \omega^\perp $, if we also treat $ \omega $ as a unit vector. 
	
	If $ \phi: I \to \R $ is a function defined on $ I \subset \R $, we denote by $ \gr(\phi; I, \omega) $ the graph of $\phi$ over $I$ (with respect to the standard coordinate system) rotated by the angle $\omega$.

	\begin{lemma}\label{lem.graph}
		Let $ k \geq 4 $ and let $ \cnice $ be sufficiently large. Suppose $ Q \in \knice $. Then there exist $ \omega \in [-\pi/2, \pi/2] $ and a twice continuously differentiable function $ \phi: \R \to \R $ such that
		\begin{itemize}
			
			\item $ E \cap Q^* \subset \gr(\phi; \R, \omega) $;
			
			\item $ \abs{\phi'} \leq 1 $, and
			
			\item $ \abs{\phi''} \leq \delta_Q^{-1} $.
			
		\end{itemize}
		The constant $ C $ depends only on $ \cnice $. 
	\end{lemma}

	\begin{proof}
		If $ E \cap Q^* = \void $, there is nothing to prove. From now on, we assume $ E \cap Q^* \neq \void $.
		
		Fix $x_0 \in E \cap Q^*$. Let $\delta = \delta_Q$. Since $ Q \in \knice $, we have $ \diam{\sk (x_0, k)} \geq \cnice\delta$.
		Since $ \sk $ is symmetric about the origin, there exist  $P^{x_0} \in \sk(x_0, k)$ and $ \omega = \frac{I_{x_0} (P^{x_0}) }{\norm{ I_{x_0} (P^{x_0}) }} $ (where $ I_{x_0} $ is the identification map in \eqref{eq. jet and euclidean id} and $ \norm{\,\cdot\,} $ is the Euclidean norm) such that 
		\begin{equation}
		\abs{\partial_2 P^{x_0}(x_0)} \geq \cnice\delta/2\,
		\label{eq. jet est 1}
		\end{equation}
		and
		\begin{align}
		\partial_1 P^{x_0}(x_0) &= 0\, . 
		\label{eq. jet est 2}
		\end{align}
		Here, $ \partial_i = \partial_{x_\omega^{(i)}} $ for $ i = 1,2 $. 
		
		\begin{claim}\label{claim: S}
			Given any $ \epsilon_0 > 0 $, we may pick $ \cnice > 0 $ large enough such that the following holds.
			
			For any $ S \subset E \cap Q^* $ containing $ x_0 $ with $ \#(S) \leq k $, there exists $ \phi^S \in C^2(\R) $ such that
			\begin{enumerate}[(i)]
				\item $ S \subset \gr \brac{ \phi^S; I^S, \omega } $,
				\item $ \abs{\brac{\phi^S}'} \leq \epsilon_0 $ on $ I^S $, and
				\item $ \abs{\brac{\phi^S}''} \leq \epsilon_0{\delta_Q^{-1}} $ on $ I^S $. 
			\end{enumerate}
			
		\end{claim}
		
		\begin{proof}[Proof of Claim \ref{claim: S}]
			Let $ S \subset E \cap Q^* $ be such that $ x_0 \in S $ and $ \#(S) \leq k $. 
			
			Since $P^{x_0} \in \sk(x_0, k)$, there exists $F^S \in C^2(\R^2)$ such that
			\begin{enumerate}[(i)]
				\item $ F^S\big|_S = 0$,
				\item $\norm{F^S}_{C^2(\mathbb{R}^2)} \leq 1$, and
				\item $\jet{x_0}{F^S}= P^{x_0}$.
			\end{enumerate}
			By \eqref{eq. jet est 1}, we have 
			\begin{equation}
			\abs{ \partial_2 F^S(x_0)  } = \abs{  \partial_2 P^{x_0} } \geq \cnice\delta_Q/2. 
			\label{eq. del 2 est.}
			\end{equation}
			Now, for all $x \in Q^*$, we have $\abs{x_0 - x} \leq 3\delta_Q$. Hence, for all $x \in Q^*$, by \eqref{eq. jet est 2} and property (ii) of $F^{S}$, we have  
			\begin{equation}
			\abs{\partial_1 F^S(x)} \leq \norm{F^S}_{C^2(\mathbb{R}^2)} \abs{x_0 - x} \leq 3\delta_Q.
			\label{eq. del 1 est}
			\end{equation} 
			From \eqref{eq. del 2 est.}, we also have, for all $x \in Q^*$, 
			\begin{equation*}
			\abs{\partial_2 F^S(x)} \geq \abs{\partial_2 F^S(x_0)} - \norm{F^S}_{C^2(\mathbb{R}^2)}\abs{x_0 -x} \geq (\cnice/2 - 3) \delta_Q.
			\end{equation*}
			Therefore, if $\cnice$ is sufficiently large, the implicit function theorem yields a function $\phi^S \in C^2(I^S)$ for some open interval $ I^S $ such that $ S \subset \gr\brac{\phi^S; I^S, \omega} $. 
			
			First we compute the derivatives of $\phi^S$: 
			\begin{align}
			\brac{\phi^S} '(x^{(1)}_\omega) &= -\frac{\partial_1 F^S(x)}{\partial_2 F^S(x)} \label{phi1} \\
			\brac{\phi^S} ''(x^{(1)}_\omega) &= \frac{-\brac{\partial_2 F^S(x)}^2\partial_1^2 F^S(x)+2\partial_1 F^S(x) \partial_2 F^S(x)\partial_{12}^2 F^S(x)-\brac{\partial_1 F^S(x)}^2 \partial_2^2 F^S(x)}{\brac{\partial_2 F^S(x)}^3}. \label{phi2}
			\end{align}
			From \eqref{eq. del 1 est} - \eqref{phi2}, we conclude that, for sufficiently large $ \cnice $, 
			\begin{equation*}
			\big|{({\phi^S})'}\big| \leq  \epsilon_0
			\text{ on } I^S
			\enskip
			\text{ and }
			\enskip
			\big|{({\phi^S})''}\big| \leq \epsilon_0{\delta_Q^{-1}}
			\text{ on } I^S\,.
			\end{equation*}

			This concludes the proof of the claim.
			
		\end{proof}

		Next, we define the projections $\pi_i: \mathbb{R}^2 \rightarrow \mathbb{R}$ by $\pi_i((x^{(1)}_\omega,x^{(2)}_\omega))=x^{(i)}_\omega$, for $i=1,2$. By Claim \ref{claim: S}, we know that $\pi_1|_{E \cap Q^*}$ is a one-to-one map. Therefore, $E \cap Q^*$ lies on a graph with respect to the $x^{(1)}_{\omega}$-axis. 
		
		It remains to see that the graph can be taken to have controlled derivatives. 
		
		For simplicity of notation, we suppress $\omega$ in the subscript.
		
		Let $x_{0} = (x_{0}^{(1)}, x_{0}^{(2)})$. We may assume without loss of generality that $ \pi_{1}(E \cap
		Q^{*}) = \set{ x_{0}^{(1)}, x_{1}^{(1)}, \dots, x_{L-1}^{(1)} } $ such that $ x_{0}^{(1)} <  x_{1}^{(1)}  < \dots <  x_{L-1}^{(1)} $, where $ L = \#(E \cap
		Q^{*}) $. Let $ \pi_{2}(E \cap Q^{*}) = \set{  x_{0}^{(2)}, x_{1}^{(2)}, \dots, x_{L-1}^{(2)} } $, where $ x_{i}^{(2)} = \pi_{2}\circ \pi_{1}^{-1} (x_{i}^{(1)}) $ for $i = 1, \dots, L-1$.
		
		Let $E_{j} = \set{ x_{j}^{(1)}, x_{j+1}^{(1)}, x_{j+2}^{(1)}  }$ for $j = 1, \dots, L - 3$. Let $S_j=\pi_1^{-1}(E_j)\cup \{x_{0}\}$. By Claim \ref{claim: S}, we know that there exist $\phi^{S_j} \in C^2(I_j)$ and a constant $ C $, depending only on $ \cnice $, such that 
		\begin{itemize}
			\item $\phi^{S_j}\big|_{E_j} = \pi_2\circ \pi_1^{-1}$,
			\item  $\abs{  ({\phi^{S_j}}) '(x^{(1)})} \leq \epsilon_0$ for all $x^{(1)} \in [x^{(1)}_{j}, x^{(1)}_{j+2}]$, and 
			\item $\abs{ ({\phi^{S_j}} )''(x^{(1)})} \leq \epsilon_0{\delta_Q^{-1}}$ for all $x^{(1)} \in [x^{(1)}_{j},x^{(1)}_{j+2}]$.
		\end{itemize}
		Therefore, by Theorem \ref{thm.FP-1dpm} and the fact that $ \delta_Q \leq 1 $, we may choose $ \epsilon_0 $ sufficiently small such that there exists $\phi \in C^2(\mathbb{R})$ such that 
		\begin{itemize}
			\item $\phi|_{E \cap Q^*}=\pi_2\circ \pi_1^{-1}$, 
			\item  $ \norm{\phi'}_{C^{0}(\R)}  \leq 1$, and 
			\item $ \norm{\phi''}_{C^{0}(\R)} \leq \delta_Q^{-1}$.
		\end{itemize}
		This completes the proof of the lemma.
	\end{proof}

	For future reference, we make the following definition.
	
	\begin{definition}\label{def.good geo}
		A pair $ (k, \cnice) $ \textbf{guarantees good geometry} if the following hold:
		\begin{itemize}
			\item $ k \geq 4 $; and
			\item $ \cnice $ is sufficiently large such that Lemma \ref{lem.graph} holds.
		\end{itemize}
	\end{definition}

	\begin{lemma}\label{lem.diffeo}
		Let $ (k, \cnice) $ guarantee good geometry. Let $Q \in \knice $. There exist a universal constant $ C $ and a diffeomorphism $\Phi = \Phi_Q \in C^2(\R^2, \R^2)$, such that the following hold.
		\begin{enumerate}[(A)]
			\item $ \Phi(E \cap Q^*) \subset \R \times \set{0} $;
			\item $ \big\| \nabla \Phi \big\|, \norm{\nabla \Phi^{-1}} \leq 2 $; and
			\item $ \norm{\nabla^2\Phi}, \norm{\nabla^2\Phi^{-1}} \leq C \delta_Q^{-1} $. 
		\end{enumerate}
		Here, $ \norm{\,\cdot\,} $ denotes the Euclidean norm.
	\end{lemma}
	
	\begin{proof}
		We may compose on the right by a rotation $ \omega $ if necessary, and assume $ \omega = 0 $. Such rotation will not affect the Euclidean norm. Let $ \phi $ be as in Lemma \ref{lem.graph}. Put
		\begin{equation}\label{eq. Phi def}
		\Phi (s,t) := (s,t-\phi(s))
		\ \text{and} \ 
		\Phi^{-1}(\hat{s},\hat{t}) := (\hat{s}, \hat{t}+\phi(\hat{s})).
		\end{equation}
		They are clearly inverses of each other and are twice continuously differentiable. 
		
		Property (A) follows from how we construct $ \phi $ (see Lemma \ref{lem.graph}).

		To see (B), we note that
		\begin{equation}
		\nabla \Phi (s,t) = 
		\left(\begin{matrix}
		1&0\\-\phi'(s)&1
		\end{matrix}\right)
		\text{ and }
		\nabla \Phi^{-1} (\hat{s},\hat{t}) = 
		\left(\begin{matrix}
		1&0\\\phi'(\hat{s})&1
		\end{matrix}\right)\,.
		\label{eq. Phi der est}
		\end{equation}
		Property (B) then follows from \eqref{eq. Phi der est} and the first derivative estimate of $ \phi $ in Lemma \ref{lem.graph}. 
		
		Further differentiating each matrix in \eqref{eq. Phi der est}, we see that the only nonzero terms occur when $ \partial_s $ is applied to the bottom left entries and yields $ \mp \phi'' $. Conclusion (C) then follows from the second derivative estimate of $ \phi $ in Lemma \ref{lem.graph}.
	\end{proof}

	\newcommand{\crep}{c_{\mathrm{rep}}}
	
	\begin{lemma}\label{lem.rep}
		Let $ (k, \cnice) $ guarantee good geometry. There exists a universal constant $ \crep $ such that the following holds. Let $Q \in \knice $. Then there exists $ \xqs \in Q $ with $ \dist{\xqs}{E} \geq \crep\delta_Q $. 
	\end{lemma}
	
	\begin{proof}
		If $ E \cap \frac{1}{2}Q = \void $, we may pick $ \xqs $ to be the center of $ Q $ and let $ \crep = 1/4 $. 
		
		Suppose $ E \cap \frac{1}{2}Q \neq \void $. Fix $ \hat{x} \in E \cap \frac{1}{2}Q $. There exists a universal constant $ c_1 > 0 $ such that $ B(\hat{x}, c_1\delta_Q) \subset Q $, where $ B(\hat{x}, c_1\delta_Q) $ is the ball of radius $ c_1\delta_Q $ centered at $ \hat{x} $. Let $ \Phi $ be as in Lemma \ref{lem.diffeo}. (Again, we may assume $ \omega = 0 $.) By (B) of Lemma \ref{lem.diffeo}, there exists a constant $ c_2 > 0 $, depending only on $ \cnice $, such that $ B(\Phi(\hat{x}), c_2\delta_Q) \subset \Phi(B(\hat{x}, c_1\delta_Q)) $. Recall that $ \Phi(E \cap Q^*) \subset \R \times \set{0} $. Let $ \bar{x}_Q^{\sharp} := \Phi(\hat{x}) + (0, c_2\delta_Q/2) $. Then $ \dist{\bar{x}_Q^{\sharp}}{\Phi(E \cap Q^*)} \geq c_2\delta_Q/2 $. Let $ \xqs = \Phi^{-1}(\bar{x}_Q^{\sharp}) $. By (B) of Lemma \ref{lem.diffeo} again, $ \dist{\xqs}{E \cap Q^*} \geq c_3\delta_Q $ for some $ c_3 > 0 $ depending only on $ \cnice $. Finally, since $ \xqs \in Q $, $ \dist{\xqs}{E \setminus Q^*} \geq \delta_Q/2 $. This concludes the proof of the lemma.

	\end{proof}

	Recall Definition \ref{def.Taylor ball}.

	\begin{lemma}\label{lem.sk-ball}
		Let $ (k,\cnice) $ guarantee good geometry. Let $ Q \in \knice $. Let $ \xqs $ be as in Lemma \ref{lem.rep}. Then 
		\begin{equation*}
		\sk(\xqs, 4k) \subset C\cdot \B(\xqs,\delta_Q)
		\end{equation*}
		for some universal constant $ C $. 
	\end{lemma}

	\begin{proof}
		If $ \delta_Q = 1 $, then the lemma follows from the definitions of $ \sk $ and $ \B $. 
		
		Suppose $ \delta_Q < 1 $. Then $ Q^+ $ exists and is not $ k $-nice, meaning that there exists $ \hat{x} \in  E \cap (Q^+)^* $ such that 
		\begin{equation}
		\diam{\sk(\hat{x},k)} < 2\cnice\delta_Q\,.
		\label{eq.small-diam}
		\end{equation}
		Fix such $ \hat{x} $. By our choice of $ \xqs $ (see Lemma \ref{lem.rep}), we have that
		\begin{equation*}
		\abs{\hat{x} - \xqs} \leq C\delta_Q\,.
		\end{equation*}

		Let $ P \in \sk(\xqs,4k) $. The argument in the proof of Lemma \ref{lemma: Gamma trade off} applied to $ \sk $ yields $ P' \in \sk(\hat{x},k) $ such that
		\begin{equation}
		P - P' \in C\cdot \brac{\B(\xqs,\delta_Q) \cap \B(\hat{x},\delta_Q)}\,.
		\label{eq.set-inc1}
		\end{equation}
		Moreover, since $ P' \in \sk(\hat{x},k) $, by the definition of $ \sk $, we have $ P'(\hat{x}) = 0 $. Thanks to \eqref{eq.small-diam}, we also have $ \abs{\nabla P'} \leq C\delta_Q $. Therefore, we can conclude that
		\begin{equation}
		P' \in C \cdot \B(\hat{x},\delta_Q).
		\label{eq.set-inc2}
		\end{equation}
		Taylor's theorem, together with \eqref{eq.set-inc1} and \eqref{eq.set-inc2}, implies that $ P \in C\cdot \B(\xqs,\delta_Q) $. Since $ P $ is an arbitrary element in $ \sk(\xqs,4k) $, the lemma follows.
	\end{proof}

	\section{1-D Results} 
	\label{section: 1D finiteness principle}

	In this section, we provide the proofs for our one-dimensional results. First, we will prove Theorem \ref{thm.FP-1dpm} and indicate how the proof of Theorem \ref{thm.FP-1d} follows. Then, we will sketch a proof for Theorem \ref{thm.BD-1d}. The proof for Theorem \ref{thm.BD-1dpm} uses the same idea but with easier intermediate steps.

	We will use $ x,y $ to denote points on $ \R $, and $ \partial^m $ to denote the $ m $-th derivative of a single-variable function. When $ m = 1 $, we simply write $ \partial $ instead of $ \partial^1 $. We use $ \P $ to denote the vector space of one-variable polynomials with degree no greater than one.

	\subsection{Finiteness Principles for $ C^2(\R) $ and $ C^2_+(\R) $}

	\begin{proof}[Proof of Theorem \ref{thm.FP-1dpm}]
		For $ N \geq 3 $, let
		$ I_1 = (-\infty, x_3] $, $ I_2 = [x_2, x_4] $, $ \dots $, $ I_{N-3} = [x_{N-3}, x_{N-1}] $, and $ I_{N-2} = [x_{N-2}, +\infty) $. By assumption, for each $ j $, there exists $ F_j \in C^2_+(\R) $ with $ F_j\big|_{E_j}= f $ and 
		\begin{equation}\label{eq. 3pt est}
		\abs{F_j}\leq A_0
		\,,
		\abs{\partial F_j} \leq A_1
		\,,
		\abs{\partial^2F_j} \leq A_2
		\,.
		\end{equation}
		
		We introduce a partition of unity $\set{\theta_j}$ that satisfies
		\begin{enumerate}[(i)]
			\item $ \sum_{j = 1}^{n-2}\theta_j \equiv 1 $ on $ \R $;
			
			\item $ \supp{\theta_j} \subset I_j $ for each $ j = 1,\dots, N-2 $; and
			
			\item \footnote[2]{
				For the existence of such partition function, see e.g. \cite{W34}.
			}
			for each $ 1 \leq k \leq 2 $ and $ 1 \leq j \leq N-2 $, 
			\begin{equation}\label{eq. theta_j est}
			\abs{\partial^k \theta_j(x)} \leq 
			\begin{cases}
			C\abs{x_{j+1} - x_j}^{-k} &\text{ if } x \in [x_{j}, x_{j+1}]\\
			C\abs{x_{j+2} - x_{j+1}}^{-k}&\text{ if } x \in [x_{j+1}, x_{j+2}]
			\end{cases}\,.
			\end{equation}

		\end{enumerate}
		Notice that the interior of $ I_i \cap I_j $ supports at most two partition functions ($ \theta_i $ and $ \theta_j $). 
		
		Define
		\begin{equation}\label{eq: 1D ext map def}
		F(x) = \sum_{j = 1}^{N-2}\theta_j(x)F_j(x).
		\end{equation}
		Clearly, $ F|_E = f $, $ F $ is twice continuously differentiable, and 
		\begin{equation}\label{eq: sup norm of f}
		\abs{F} \leq 2A_0.
		\end{equation}
		
		Observe that (\ref{eq. 3pt est}) and condition (ii) of $ \set{\theta_j} $ imply
		\begin{equation}\label{eq. F_k est-far}
		\abs{\partial^mF} \leq A_m
		\text{ on }
		(-\infty, x_2] \cup [x_{N-1}, +\infty)\,.
		\end{equation}
		
		Suppose $ x \in (x_2, x_{N-1}) $. Let $ j $ be the least integer such that $ x \in I_j $. The only partition functions possibly nonzero at $ x $ are $ \theta_j $ and $ \theta_{j+1} $. Since $ \theta_j(x) + \theta_{j+1}(x) \equiv 1 $, we have $ \partial^k\theta_j(x) = -\partial^k\theta_{j+1}(x) $ for $ k = 1,2 $. Thus,
		\begin{equation}\label{eq. Fk-expand}
		\partial^k F(x)
		= \partial^k F_j(x)\theta_j(x) + \partial^k F_{j+1}(x)\theta_{j+1}(x)  + \sum_{l = 0}^{k-1} \binom{k}{l}\partial^{l}(F_j -F_{j+1})(x)\partial^{k - l}\theta_j(x).
		\end{equation}
		
		\begin{claim}\label{claim 5}
			Let $ x \in I_{j} \cap I_{j+1} $. Then 
			\begin{equation}
			\abs{(F_j - F_{j+1})(x)} \leq 2A_1 \abs{x_{j+1} - x_{j}}\,. \label{eq. alt est 1}
			\end{equation}
			
			For $ l = 0, 1 $, we also have 
			\begin{equation}
			\abs{\partial^l(F_j - F_{j+1})(x)} 
			\leq 2A_2 \abs{x_{j+1} - x_{j}}^{2-l} \,.
			\label{eq. neighbor est}
			\end{equation}
			
		\end{claim}
		
		\begin{proof}[Proof of Claim \ref{claim 5}]
			Note that by construction, $ I_j \cap I_{j+1} = [x_{j+1}, x_{j+2}]$. 
			
			Observe that \eqref{eq. alt est 1} is an immediate consequence of the mean value theorem.
			
			It remains to show \eqref{eq. neighbor est}.
			
			Observe that $ (F_j - F_{j+1})(x_{j+1}) = (F_j - F_{j+1})(x_{j+2}) = 0 $. By Rolle's theorem, there exists $ \hat{x}_j \in (x_j, x_{j+1}) $ such that $ \partial(F_j - F_{j+1})(\hat{x}_j) = 0 $. By the fundamental theorem of calculus and triangle inequality, we have
			\begin{equation*}
			\begin{split}
			\abs{\partial(F_j - F_{j+1})(x)} 
			\leq \int_{\hat{x}_j}^{x} \abs{\partial^2(F_j - F_{j+1})(y)}dy 
			\leq 2A_2\abs{x_{j+2} - x_{j+1}}  	\text{ for all } x \in I_j \cap I_{j+1}.
			\end{split}
			\end{equation*}
			Similar calculations yield the case $ l = 0 $. (\ref{eq. neighbor est}) is proven.
		\end{proof}
		
		Now, \eqref{eq. theta_j est}, \eqref{eq. F_k est-far}, \eqref{eq. Fk-expand}, and \eqref{eq. alt est 1} imply that
		\begin{equation}\label{eq. conc1}
		\abs{\partial F} \leq CA_1.
		\end{equation}
		
		Likewise, \eqref{eq. theta_j est}, \eqref{eq. F_k est-far}, \eqref{eq. Fk-expand}, and \eqref{eq. neighbor est} imply that
		\begin{equation}\label{eq. conc2}
		\abs{\partial^2 F} \leq CA_2.
		\end{equation}
		
		In view of \eqref{eq: sup norm of f}, \eqref{eq. conc1}, and \eqref{eq. conc2}, we conclude the proof of the theorem.

	\end{proof}

	\begin{proof}[Proof of Theorem \ref{thm.FP-1d}]
		We simply take $ A_m = 1 $ for $ m = 0,1,2 $ in the above proof of Theorem \ref{thm.FP-1dpm}, and note that $ F(x) $ defined by \eqref{eq: 1D ext map def} is nonnegative if all of the $ F_j $'s are nonnegative.
	\end{proof}

	\subsection{$ C^2(\R) $ and $ C^2_+(\R) $ extension operators of bounded depth}

	Now we explain the proof of Theorem \ref{thm.BD-1d}. 
	
	Let $ E \subset \R $ be a finite set. We enumerate $ E = \set{x_1, \cdots, x_N} $ with $ x_1 < \cdots < x_N $. Let $ E_i := \set{x_i, x_{i+1}, x_{i+2}} $ for $ i = 1, \cdots, N-2 $. Suppose for each $ i $, we are given an extension operator $ \Eb_i: C^2_+(E_i) \to C^2_+(\R) $ with $ \norm{\Eb_i (f)}_{C^2(\R)} \leq C\norm{f}_{C^2_+(E_i)} $ and $ \brac{\Eb_i (f)}\big|_{E_i} = f $. Let $ \set{I_i} $ and $ \set{\theta_i} $ be as in the proof of Theorem \ref{thm.FP-1d}. We define
	\begin{equation}
	\Eb(f)(x) := \sum_{i = 1}^{N-2}\theta_i(x) \cdot \Eb_i(f)(x)\,. 
	\label{eq.1Dext}
	\end{equation}
	Conclusions (A) and (B) of Theorem \ref{thm.BD-1d} follow from the same argument as in the proof of Theorem \ref{thm.FP-1d}. Moreover, by assumption, $ \Eb_i(f) $ depends only on $ \set{f(x_i), f(x_{i+1}), f(x_{i+2})} $ for each $ i $, and the $ \theta_{i} $'s have bounded overlap. Therefore, conclusion (C) and Remark \ref{rem.1D-depth} follow.
	
	Hence, in order to construct a bounded extension operator with bounded depth in dimension one, it suffices to construct a bounded extension operator for every consecutive three points. This is a routine linear algebra problem and is readily solvable via the nonnegative Whitney extension theorem (see Lemma \ref{lem.WE-map}). We leave the details to the interested readers.
	
	For Theorem \ref{thm.BD-1dpm}, we simply replace each summand on the right-hand side in \eqref{eq.1Dext} with $ \theta_i \cdot \Epmb_i $, where $ \Epmb_i $ is an extension operator associated with $ E_i $ without the nonnegative constraints.

	\subsection{Non-additivity}

	In this section, we use the following notations
	\begin{equation*}
	\begin{split}
	\norm{f}_{\dot{C}^m(E)} &:= \inf\set{\norm{F}_{\dot{C}^m(\R)} : F \in C^m(\R)\enskip\text{ and }\enskip F|_E = f}
	\text{ and }
	\\
	\norm{f}_{\dot{C}^m_+(E)} &:= \inf\set{\norm{F}_{\dot{C}^m(\R)} : F \in C^m_+(\R)\enskip\text{ and }\enskip F|_E = f}\,.
	\end{split}
	\end{equation*}

	\begin{proof}[Proof of Theorem \ref{thm.nonlinear}]
		
		Let $ \epsilon > 0 $ be a sufficiently small number. We use $ C, C', C_* $ etc.~to denote universal constants.
		
		Consider $ E = \set{x_1, x_2, x_3} \subset \R $, where $ x_j = (j-1)\epsilon $ for $ j = 1, 2, 3 $. Suppose toward a contradiction, that $ \Eb : C^2_+(E) \to C^2_+(\R) $ is a bounded extension map that is additive. That is, $ \Eb(f+g) = \Eb(f) + \Eb(g) $ for all $ f, g \in C^2_+(E) $, and
		\begin{equation*}
		C^{-1}\norm{\Eb (f)}_{C^2(\R)} \leq \norm{f}_{C^2_+(E)} \leq C \norm{\Eb (f)}_{C^2(\R)}\,. 
		\end{equation*}
		For $ j = 1, 2, 3 $, we define
		\begin{equation*}
		f(x_j) := (j-1)\epsilon
		\,\,\,\,\text{ and }\,\,\,\,
		g(x_j) := 1 - f(x_j)\,.
		\end{equation*}
		Then $ f, g \in C^2_+(E) $, and $ f + g \equiv 1 $. It is easy to see that
		\begin{equation*}
		\norm{f + g}_{C^2_+(E)} = 1\,.
		\end{equation*}
		In fact,
		\begin{equation*}
		\norm{f+g}_{\dot{C}^m(E)} = 0
		\text{ for }
		m = 1,2\,.
		\end{equation*}
		Since $ \Eb $ is bounded, we have
		\begin{equation}
		1 \leq \norm{\E(f+g)}_{C^2(\R)} \leq C\,.
		\label{eq.Tfg}
		\end{equation}

		We analyze the trace norms of $ f $ and $ g $.

		We begin with $ f $. By calculating the divided difference, we see that
		\begin{equation*}
		\norm{f}_{\dot{C}^1_+(E)} \geq 1\,.
		\end{equation*}
		By Taylor's theorem (see Lemma \ref{lem.cone-equiv}), we have
		\begin{equation}
		\norm{f}_{\dot{C}^2_+(E)} \geq C\cdot  \frac{\norm{f}_{\dot{C}^1(E)}^2}{f(x_j)} \geq  C'\epsilon^{-1}\,.
		\label{eq.contra}
		\end{equation}
		Here, we use nonnegativity of the extension.
		
		In particular, \eqref{eq.contra} implies
		\begin{equation}
		\abs{\d^2\E f(x_0)} \geq C_0\epsilon^{-1}
		\label{eq.fx0}
		\end{equation}
		for some $ x_0 $. Fix such $ x_0 $.

		Now we turn to $ g $. 
		
		Let $ \psi $ be a cutoff function such that $ \psi \equiv 1 $ in a neighborhood of $ [0,2\epsilon] $, $ \supp{\psi} \subset [-1,1]  $, and $ \abs{\d^m\chi} \leq C $ for $ m = 0,1,2 $. Consider the function $ \tilde{g} $ defined by
		\begin{equation*}
		\tilde{g}(x) := \psi(x)\cdot(1-x)
		\end{equation*}
		It is clear that $ \tilde{g} \in C^2_+(\R) $ with $ \tilde{g}|_E = g $. Moreover, 
		\begin{equation*}
		\norm{\tilde{g}}_{C^2(\R)} \leq C\,.
		\end{equation*}
		Therefore, 
		\begin{equation*}
		\norm{g}_{C^2_+(E)} \leq C\,.
		\end{equation*}
		Since $ \Eb $ is bounded, we know that, for $ x_0 $ as in \eqref{eq.fx0},
		\begin{equation}
		\abs{\d^2 \Eb (g)(x_0)} \leq C_1\,.
		\label{eq.gx0}
		\end{equation}
		Therefore, we have, with $ C_0 $ and $ C_1 $ as in \eqref{eq.fx0} and in \eqref{eq.gx0},
		\begin{equation*}
		\abs{\d^2(\E f + \E g)(x_0)} \geq C_0\epsilon^{-1} - C_1 \geq C\epsilon^{-1}\,.
		\end{equation*}
		For sufficiently small $ \epsilon $, this would contradict \eqref{eq.Tfg}. 
	\end{proof}

	\section[2-D Finiteness Principles]{2-D Finiteness Principle}
	\label{section: 2DFPloc}

	\subsection{Statement of the main local lemma} The goal of this section is to prove a local version of the finiteness principle,  which produces a nonnegative local interpolant taking a jet in some prescribed $ \G $ (see \eqref{eq.Gk-def}) at a point sufficiently far away from the data. We will use these jets as transitions in our estimates. 
	
	Recall Definition \ref{def.good geo}. Also recall that Lemma \ref{lem.rep} produces a point $ \xqs \in Q $ such that 
	\begin{equation}
	\dist{\xqs}{E} \geq \crep\delta_Q
	\label{eq. xsharp}
	\end{equation}
	for each $ Q \in \knice $ given that $ (k, \cnice) $ guarantees good geometry. We fix the number $ \crep $. 
	
	\newcommand{\kloc}{k_{\mathrm{loc}}}
	\begin{lemma}\label{lem.FP-loc}
		Let $ E \subset \R^2 $ be finite, and let $ f : E \to [0,\infty) $. Let $ (k, \cnice) $ guarantee good geometry and $ Q \in \knice $. Let $ \xqs \in Q $ satisfy \eqref{eq. xsharp}. Let $ \kloc \geq 3 $. Suppose $ \G(\xqs, \kloc, M) \neq \void $. Then there exist a universal constant $ C $ and a function $ F^{\sharp}_Q \in C^2_+(100Q) $ such that the following hold.
		\begin{enumerate}[(A)]
			\item $ F^{}_Q\big|_{E \cap Q^*} = f $,
			
			\item $ \norm{F_Q}_{C^2(100Q)} \leq CM $, and
			
			\item $ \jet{\xqs}{F_Q} \in \G(\xqs, \kloc, CM) $.
		\end{enumerate}
		
	\end{lemma}

	Note that if $ \#(E \cap Q^*) \leq \kloc $, the conclusion follows immediately. 
	
	Hereafter, we assume $ \#(E \cap Q^*) > \kloc \geq 3 $.
	
	The main idea of the proof is to treat the local interpolation problem differently depending on whether the local data is big or small. For big local data, we solve the problem as if there were no nonnegative constraints. For small local data, we simply prescribe a zero jet.

	Below we give a more detailed overview of our strategy, still without dwelling into the technicalities. 
	
	Our approach relies on three crucial lemmas. The first one (Lemma \ref{lem.TC}) describes the relationships among the value, gradient, and zero set of a jet generated by a nonnegative function. The second one (Lemma \ref{lem.perturb}) is a perturbation lemma, which specifies the conditions under which we are allowed to modify an element in $ \G(\xqs, \,\cdot\,\,,\,\cdot\,\,) $. We emphasize the importance of the choice of base point $ x^{\sharp}_Q $, which is far away enough from all the data points (on the order of $ \delta_Q $) so that we have room to modify the interpolants' behavior near $ \xqs $. The third one (Lemma \ref{lem.BS}) tells us that the local data is either uniformly big or uniformly small (on the order of $ \delta_Q^2 $).
	
	We begin the proof of Lemma \ref{lem.FP-loc} by first tackling a one-dimensional interpolation problem. Recall that, thanks to Lemma \ref{lem.graph}, the data points locally lie on a curve. The interpolation problem along this curve is essentially one-dimensional and readily solved, thanks to Theorems \ref{thm.FP-1d}, \ref{thm.FP-1dpm}, and Lemma \ref{lem.diffeo}. 
	
	We then solve the local problem when the local data is uniformly large, namely, $ \min_{x \in E \cap Q^*} f(x) \geq B\delta_Q^2 $ for some universal $ B > 0 $ to be determined. We replace the local data $ f|_{E \cap Q^*} $ by $ g(x) = f(x) - P^\sharp(x) $ for $ x \in E \cap Q^* $, where $ P^\sharp $ is a suitable element in $ \G(x^{\sharp}_Q, \kloc, C) $ such that $ g $ achieves two zeros and that $ P^\sharp \geq B'\delta_Q^2 $ on $ 100Q $ for some $ B' > 0 $ depending only on $ B $. Thanks to Rolle's theorem, the resulting one-dimensional $ g $-interpolant, although not necessarily nonnegative, will be uniformly small on the order of $ \delta_Q^2 $, and in particular, bounded from below by $ -c\delta_Q^2 $. 
	Now, we are in the suitable order of magnitude to force a zero jet at $ x^{\sharp} $. 
	To do this, we simply extend the one-dimensional interpolant in the normal direction by constant, and use a bump function to damp out the function at $ \xqs $. 
	If we choose $ B $ such that $ B' $ is bigger than $ c $, we may add $ P^\sharp $ back to the zero-jet interpolant while preserving nonnegativity of the sum on $ 100Q $, and solve the local problem.
	
	Next, we solve the local problem when the data is not uniformly big. Thanks to Lemma \ref{lem.BS}, the local data has to be uniformly small, i.e., $ \max_{x \in E \cap Q^*}f(x) \leq B''\delta_Q^2 $ for some $ B'' > 0 $ depending only on $ B $. Therefore, we are in the correct order of magnitude to force a zero jet as in the previous step. Thanks to the perturbation lemma (Lemma \ref{lem.perturb}), the zero jet in this case is indeed a $ \kloc $-point jet, and the problem is solved.
	
	Sections \ref{sect:lemma-loc} and \ref{sect:lip} will be devoted to the proof of Lemma \ref{lem.FP-loc}.

	\subsection{Key lemmas}\label{sect:lemma-loc}
	
	In this section, we use Cartesian coordinates $ x = (s,t) $ on $ \R^2 $. We also write $ \xqs = x^{\sharp}_Q = (s_Q^\sharp, t_Q^\sharp) $.

	\begin{lemma}\label{lem.TC}
		There exist universal constants $ C, C', C'' $ such that the following hold. Suppose $ P \in \Gamma_+(x, \void, M) $. Then 
		\begin{align}
		P(y) + CM\abs{y-x}^2 &\geq 0
		\, \text{ for all } y \in \R^2,
		\label{eq. TC}\\
		\abs{\nabla P} &\leq C'\sqrt{MP(x)}\,\text{ and }
		\label{eq. TC1}
		\\
		\dist{x}{\set{P = 0}} &\geq C''M^{-1/2}\sqrt{P(x)} \,.\label{eq. TC2}
		\end{align}
	
	\end{lemma}

	\begin{proof}
		\eqref{eq. TC} is a direct consequence of Taylor's Theorem.
		
		To see \eqref{eq. TC1}, we simply compute the discriminants of the left hand side of \eqref{eq. TC} restricted to the $ s $ and $ t $-directions.
		
		Now we prove \eqref{eq. TC2}. If $ P(x) = 0 $ or $ P $ is a constant polynomial, the inequality is obvious.  Assume that $ P(x) > 0 $ and $ P $ is nonconstant.
		
		Since $ P $ is an affine function and the gradient points toward the direction of maximal increase, we have
		\begin{equation}
		\abs{\nabla P}  = \frac{P(x)}{\dist{x}{\set{P = 0}}}
		\,.
		\label{eq. TC3}
		\end{equation}
		From \eqref{eq. TC1} and \eqref{eq. TC3}, we have the desired estimate.
	\end{proof}

	\begin{lemma}\label{lem.perturb}

		Let $ M > 0 $. Let $ (k,\cnice) $ guarantee good geometry (see Definition \ref{def.good geo}), let $ k' \geq 1 $, and let $ Q \in \knice $. Let $ \xqs $ be as in Lemma \ref{lem.rep}. Suppose $ E \cap Q^* \neq \void $. Suppose $ \G(\xqs,k',M) \neq \void $.
		
		\begin{enumerate}[(A)]

			\item There exists a number $ B > 0 $ exceeding a large universal constant such that the following holds. Suppose $ f(x) \geq BM\delta_Q^2 $ for each $ x \in E \cap Q^* $. Then 
			\begin{equation*}
			\G(\xqs,k',M) + M\cdot \B(\xqs,\delta_Q) \subset \G(\xqs,k',CM). 
			\end{equation*}

			\item Let $ A > 0 $. Suppose $ f(x) \leq AM\delta_Q^2 $ for some $ x \in E \cap Q^* $. Then
			\begin{equation*}
			0 \in \G(\xqs,k',A'M)\,.
			\end{equation*}
			The number $ A' $ depends only on $ A $. 
		\end{enumerate}

	\end{lemma}

	\begin{proof}
		We prove (A) first.
		
		Let $ B > 0 $ be sufficiently large.

		\begin{claim}\label{claim.222}
			Under the hypothesis of (A). Given any $ P \in \G(\xqs,k',M) $, we have
			\begin{equation*}
			P(\xqs) \geq B_0M\delta_Q^2\,,
			\end{equation*} where we can take $ B_0 = C(\sqrt{B}-1/2)^2 $. 
		\end{claim}
		
		\begin{proof}[Proof of Claim \ref{claim.222}]
			We repeat proof of Claim \ref{claim: b} with more control on the parameters.
			
			Let $ x \in E\cap Q^* $. By definition, there exists $ F \in C^2_+(\R^2) $ with $ \norm{F}_{C^2(\R^2)} \leq M $, $ \jet{\xqs}{F} = P $, and
			\begin{equation}
			F(x) = f(x) \geq BM\delta_Q^2\,.
			\label{eq.Fbig}
			\end{equation}
			Suppose toward a contradiction, that $ F(\xqs) < B_0M\delta_Q^2 $. We see from \eqref{eq. TC1} that $ \abs{\nabla F(\xqs)} \leq C\sqrt{B_0}M\delta_Q $. By the fundamental theorem of calculus, we have
			\begin{equation*}
			\abs{\nabla F(x)} \leq \abs{\nabla F(\xqs)} + C\norm{F}_{C^2(\R^2)}\delta_Q
			\leq C'\brac{\sqrt{B_0} + \frac{1}{4}}M\delta_Q
			\text{ on }
			Q^*\,.
			\end{equation*}
			By the fundamental theorem of calculus again, we have
			\begin{equation*}
			F(x) \leq F(\xqs) + C\delta_Q\cdot\sup_{x \in Q^*} \abs{\nabla F(x)}
			\leq C'\brac{B_0 + \sqrt{B_0} + \frac{1}{4}}M\delta_Q^2 = C'\brac{\sqrt{B_0} + \frac{1}{2}}^2M\delta_Q^2
			\text{ on }Q^*\,.
			\end{equation*}
			If we pick $ B_0 $ to be so small that $ \sqrt{B_0} < \frac{\sqrt{B} - \frac{1}{2}}{C'} $ with $ C' $ as above, we will contradict \eqref{eq.Fbig}. This proves the claim.
			
		\end{proof}
		
		Pick $ P \in \G(\xqs,k',M) $. By the claim, we know that $ P(\xqs) \geq B_0M\delta_Q^2 $. 
		
		Let $ \tilde{P} \in M\cdot\B(\xqs,\delta_Q) $. By definition, we have
		\begin{equation*}
		\abs{\d^\alpha \tilde{P}(\xqs)} \leq M\delta_Q^{2-\abs{\alpha}}
		\text{ for }
		\abs{\alpha} \leq 2\,.
		\end{equation*}

		Let $ S \subset E $ with $ \#(S) \leq k' $. We want to show that there exists $ F \in C^2_+(\R^2) $ with $ \norm{F}_{C^2(\R^2)} \leq CM $, $ F(x) = f(x) $ for each $ x \in S$, and $ \jet{\xqs}{F} = P + \tilde{P} $.
		
		We enumerate $ S = \set{x_1, \cdots, x_{k'}} $. We let $ \tilde{S} = \set{x_0, x_1, \cdots, x_{k'}} $ with $ x_0 := \xqs $.

		By the definition of $ \G $, there exists 
		\begin{equation}
		F^S \in C^2_+(\R^2)
		\text{ with }
		\norm{F^S}_{C^2(\R^2)} \leq M,\,
		F^S \big|_S = f,\,
		\text{ and }
		\jet{\xqs}{F^S} = P\,.
		\label{eq.goodFS}
		\end{equation}
		For $ i = 1, \cdots, k' $, we set
		\begin{equation*}
		P^{x_i} := \jet{x_i}{F^S}\,.
		\end{equation*}
		We also set
		\begin{equation}
		P^{x_0} := P + \tilde{P}.
		\label{eq.7.3.1}
		\end{equation}
		We put
		\begin{equation*}
		\vec{P} := \brac{P^{x_i}}_{i = 0}^{k'} \in W^2(\tilde{S}).
		\end{equation*}
		Thanks to Lemma \ref{lem.Taylor-WF} and Lemma \ref{lem.WE-map}, it suffices to show that $ \vec{P} \in W^2_+(\tilde{S}) $ and $ \norm{\vec{P}}_{W^2_+(S)} \leq CM $. 
		
		Thanks to \eqref{eq.goodFS}, we have
		\begin{equation}
		P^{x_i} \in \cone(x_i,CM)
		\text{ for all }
		i = 1, \cdots, k'\,,
		\label{eq.dum}
		\end{equation}
		and
		\begin{equation}
		P^{x_i} - P^{x_j} \in CM\cdot \B(x_i,\abs{x_i - x_j})
		\text{ for all } i, j = 1, \cdots, k'\,.
		\label{eq.Pxi}
		\end{equation}
		On the other hand, thanks to Claim \ref{claim.222}, we have
		\begin{equation*}
		P^{x_0}(x_0) = P(x_0) + \tilde{P}(x_0) \geq (B_0-1)M\delta_Q^2 \geq 0\,,
		\end{equation*}
		and
		\begin{equation*}
		\abs{\nabla P^{x_0}} \leq \abs{\nabla P} + \abs{\nabla \tilde{P}} \leq C\sqrt{MP(\xqs)} + M\delta_Q \leq C'\sqrt{M(P + \tilde{P})(\xqs) }\,.
		\end{equation*}
		This, combined with \eqref{eq.goodFS}, shows that
		\begin{equation}
		P^{x_0} \in \cone(x_0,CM)\,.
		\label{eq.dumdum}
		\end{equation}
		
		It remains to estimate $ \norm{\vec{P}}_{W^2(\tilde{S})} $.

		By Taylor's theorem, we have
		\begin{equation}
		P^{x_i} - P \in CM\cdot \brac{\B(x_i,\abs{x_i - x_0})\cap CM \cdot \B(x_0,\abs{x_i - x_0})}\,.
		\label{eq.7.3.2}
		\end{equation} 
		
		By Lemma \ref{lem.rep}, we have
		\begin{equation*}
		\dist{\xqs}{E} \geq C\delta_Q\,.
		\end{equation*}
		This, together with Taylor's theorem and the fact that $ \tilde{P} \in M\B(x_0,\delta_Q) $, implies
		\begin{equation}
		\tilde{P} \in CM\cdot \B(x_i,\abs{x_i - x_0})
		\text{ for all } i = 1, \cdots, k'\,.
		\label{eq.7.3.3}
		\end{equation}
		Therefore,
		\begin{equation}
		\begin{matrix*}[l]
		&P^{x_i} - P^{x_0} &= P^{x_i} - P - \tilde{P} \enskip\enskip&\text{(by \eqref{eq.7.3.1})} \\
		&{}&\in (-\tilde{P}) + CM\cdot\brac{ \B(x_i,\delta_Q)\cap\B(x_0,\delta_Q) }\enskip\enskip&\text{(by \eqref{eq.7.3.2})} \\
		&{}&\subset C' M\cdot \brac{ \B(x_i,\delta_Q)\cap\B(x_0,\delta_Q) }\enskip\enskip&\text{(by \eqref{eq.7.3.3})}\,.
		\end{matrix*}
		\label{eq.Px0}
		\end{equation}
		From \eqref{eq.dum}-\eqref{eq.Px0}, we can conclude that $ \norm{\vec{P}}_{W^2_+(S)} \leq CM $. This concludes the proof of (A).

		Now we turn to the proof of (B).
		
		\begin{claim}\label{Claim.Psmall}
			Assume the hypothesis of (B). Let $ P \in \G(\xqs,k',M) $. Then $ P(\xqs) \leq C(\sqrt{A}+1)^2M\delta_Q^2 $. 
		\end{claim}

		\begin{proof}[Proof of Claim \ref{Claim.Psmall}]
			
			Fix $ \hat{x} \in E \cap Q^* $ such that $ f(\hat{x}) \leq BM\delta_Q^2 $. Since $ k' \geq 1 $, by the definition of $ \G $, there exists a function $ F \in C^2_+(\R^2) $ with $ F(\hat{x}) = f(\hat{x}) \leq AM\delta_Q^2 $, $ \norm{F}_{C^2(\R^2)} \leq M$, and $ \jet{\xqs}{F} = P $. By Lemma \ref{lem.TC}, we have
			\begin{equation*}
			\abs{\nabla F(\hat{x})} = \abs{\nabla \jet{\hat{x}}{F}} \leq \sqrt{A}M\delta_Q\,.
			\end{equation*}
			By Taylor's theorem, we see that
			\begin{equation*}
			\abs{\nabla F(x)} \leq C(\sqrt{A} + 1)M\delta_Q
			\enskip\text{ for }\enskip
			x \in Q^*\,.
			\end{equation*}
			By the fundamental theorem of calculus, we see that
			\begin{equation*}
			P(\xqs) = F(\xqs) \leq C(A + \sqrt{A} + 1)M\delta_Q^2 \leq C'(\sqrt{A}+1)^2M\delta_Q^2\,.
			\end{equation*}
			This finishes the proof of Claim \ref{Claim.Psmall}.
		\end{proof}

		It remains to show that $ 0 \in \Gamma_+(\xqs, S, A'M) $ for each $ S \subset E $ with $ \#(S) \leq k' $. 
		
		We use $ A_0, A_1, $ etc. to denotes quantities that depend only on $ A $. 
		
		Fix $ P \in \G(\xqs,k',M) $. Let $ S \subset E $ satisfy $ \#(S) \leq k' $. By definition, there exists $ F^S \in C^2_+(\R^2) $, such that $ F^S \big|_S = f $, $ \norm{F^S}_{C^2(\R^2)} \leq 1 $, and $ \jet{\xqs}{F^S} = P $. By Claim \ref{Claim.Psmall}, we see that $ P(\xqs) \leq A_0M\delta_Q^2 $ and by \eqref{eq. TC1} that $ \abs{\nabla P} \leq A_1M\delta_Q $. In other words, 
		\begin{equation*}
		F^S(\xqs) \leq A_0M\delta_Q^2
		\, \text{ and }
		\abs{\nabla F^S(\xqs)} \leq A_1M\delta_Q.
		\end{equation*}
		The fundamental theorem of calculus then implies
		\begin{equation}
		\abs{\nabla F^S(x)} \leq A_2\delta_Q
		\,\text{ and } \,
		F^S(x) \leq A_2\delta_Q^2
		\,
		\text{ for all } x \in B(\xqs, \frac{\crep\delta_Q}{100}).
		\label{eq. FS small}
		\end{equation}
		Let $ \psi \in C^2_+(\R^2) $ be a cutoff function such that 
		\begin{equation}
		0 \leq \psi \leq 1
		\,,\,
		\psi \equiv 1 \text{ near } \xqs
		\,,\,
		\supp{\psi} \subset B(\xqs, \frac{\crep\delta_Q}{100})
		\,,\,
		\abs{\partial^\alpha \psi} \leq C\delta_Q^{-\abs{\alpha}}
		\text{ for }
		\abs{\alpha} \leq 2.
		\label{eq. psi}
		\end{equation}
		Let
		\begin{equation*}
		\tilde{F}^S := (1 - \psi) F^S. 
		\end{equation*}
		We have the following.
		\begin{itemize}
			\item By \eqref{eq. xsharp} and the fact that $ \supp{\psi} \subset B(\xqs, \frac{\crep \delta_Q}{100}) $, we have $ \tilde{F}^S\big|_S = f $.
			\item By \eqref{eq. psi} and the assumption that $ F^S \geq 0 $, we have $ \tilde{F}^S \geq 0 $ on $ \R^2 $. 
			\item Thanks to \eqref{eq. FS small} and \eqref{eq. psi}, $ \norm{\tilde{F}^S}_{C^2(\R^2)} \leq A_2M $.
			\item Since $ \psi \equiv 1 $ near $ \xqs $, we have  $ \jet{\xqs}{\tilde{F}^S} \equiv 0 $.
		\end{itemize}
		Since $ S $ is arbitrary, we have $ 0 \in \G(\xqs, k', A_2M) $. This completes the proof of (B) and the proof of the lemma.

	\end{proof}

		\begin{lemma}
		There exists a universal constant $ B > 0 $ such that the following holds.
		Let $ M > 0 $. Let $ (k,\cnice) $ guarantee good geometry (see Definition \ref{def.good geo}). Let $ Q \in \knice $. Let $ \xqs $ be as in Lemma \ref{lem.rep}. Suppose $ E \cap Q^* \neq \void $ and $ f(x) \geq BM\delta_Q^2 $ for all $ x \in E \cap Q^* $. Let $ k' \geq 0 $. Then
		\begin{equation*}
		\G(\xqs,k',M) + M\cdot\sk(\xqs,4k) \subset \G(\xqs,k',CM).
		\end{equation*}
	\end{lemma}

	\begin{proof}
		This is a direct consequence of Lemma \ref{lem.sk-ball} and Lemma \ref{lem.perturb}\,.
	\end{proof}

	\newcommand{\bmin}{B_{\min}}
	\newcommand{\bmax}{B_{\max}}

	\begin{lemma}\label{lem.BS}
		For each $ \bmin > 0 $, we can find $ \bmax$, depending only on $ \bmin $, such that the following holds. 
		
		Let $ E \subset \R^2 $ be a finite set. Let $ f : E \to [0,\infty) $. Let $ k' \geq 2 $. Suppose $ \G(x, k', M) \neq \void $ for all $ x \in \R^2 $. Let $ (k,\cnice) $ guarantee good geometry. Let $ Q \in \knice $. 
		Then at least one of the following holds.
		\begin{enumerate}[(A)]
			\item $ f(x) \leq \bmax M\delta_Q^2 $ for all $ x \in E \cap Q^* $.
			\item $ f(x) \geq \bmin M\delta_Q^2 $ for all $ x \in E \cap Q^* $. 
		\end{enumerate}
	\end{lemma}
	
	\begin{proof}
		Fix $ \bmin > 0 $. We use $ B, B' $, etc. to denote quantities that depend only on $ \bmin $. 
		
		Without loss of generality, we may assume $ M = 1 $. 
		
		If $ \min\limits_{x \in E \cap Q^*} f(x) \geq \bmin\delta_Q^2 $, there is nothing to prove. 
		
		Suppose there exists $ x_0 \in E \cap Q^* $ such that $ f(x_0) < \bmin\delta_Q^2 $. Fix such $ x_0 $. 
		
		Let $ S \subset E\cap Q^* $ satisfy $ x_0 \in S $ and $ \#(S) \leq k' $. Since $ \G(x, k', 1) \neq \void $ for each $ x \in \R^2 $, there exists $ F^S \in C^2_+(\R^2) $ such that $ F^S \big|_S = f $ and $ \norm{F^S}_{C^2(\R^2)} \leq 1 $. 
		
		Since $ F^S(x_0) < \bmin\delta_Q^2 $, \eqref{eq. TC1} implies there exists $ B > 0 $ such that 
		\begin{equation*}
		\abs{ \nabla F^S (x_0) } = \abs{\nabla \jet{x_0}{F^S} } \leq B\delta_Q.
		\end{equation*}
		Therefore, since $ \norm{F^S}_{C^2(\R^2)} \leq 1 $, we have $ \abs{\nabla F^S(x)} \leq B'\delta_Q $ for all $ x \in Q^* $. By the fundamental theorem of calculus, since $ F^S(x_0) < \bmin \delta_Q^2 $, we must have $ \abs{F^S(x)} \leq B''\delta_Q^2 $ for all $ x \in Q^* $. In particular,
		\begin{equation*}
		\abs{F^S(x)} \leq B''\delta_Q^2
		\,\text{ for all }
		x \in S.
		\end{equation*}
		Let $ \bmax := B'' $. Since $ S $ is arbitrary and is allowed to contain more than one point, we may conclude the proof of the lemma once we let $ S $ range over all $ k' $-point subsets of $ E \cap Q^* $ containing $ x_0 $. 
	\end{proof}

	\newcommand{\Gt}{\widetilde{\Gamma}_+^\sharp}
	\newcommand{\etot}{E}
	\newcommand{\eloc}{E_{\mathrm{loc}}}
	\newcommand{\nmax}{\nu_{\max}}

	\subsection{Solving the local problem}\label{sect:lip}
	
	In this subsection, we prove Lemma \ref{lem.FP-loc}. We fix the local data structure for the rest of the section.

	\begin{center}
		\underline{Local Data Structure (LDS)}
		\begin{itemize}
			\item A lengthscale $ \delta \leq 1$.
			
			\item A square $ Q \subset \R^2 $ with $ \delta_Q = \delta $.
			
			\item A representative point $ x^{\sharp} \in Q $ such that $ \dist{x^{\sharp}}{E} \geq \crep \delta $.
			
			\item A function $ \phi \in C^2(\R) $ that satisfies $ \abs{\phi^{(k)}} \leq \delta^{1-k} $ for $ k = 1, 2 $.
			
			\item A diffeomorphism $ \Phi : \R^2 \to \R^2 $ given by $ \Phi(s,t) = (s, t - \phi(s)) $. 
			
			\item $ \eloc = E \cap Q^* $ such that $ \eloc \subset \set{(s,\phi(s)) : s \in \R} $.

		\end{itemize}
	\end{center}

	Any $ Q \in \knice $ with $ (k,\cnice) $ guaranteeing good geometry admits the local data structure, thanks to Lemmas \ref{lem.graph}, \ref{lem.diffeo}, and \ref{lem.rep}.

	We have shown in Lemma \ref{lem.BS} that each local interpolation problem belongs to at least one of the two categories: The function's local values are uniformly big ($ \min_{x \in \eloc}f(x) \geq \bmin \delta^2 $), or are uniformly small ($ \max_{x \in \eloc}f(x) \leq \bmax\delta^2 $). The next lemma solves the former case.

	\begin{lemma}\label{lem.FP-loc big}
		There exists a sufficiently large $ \bmin > 0 $ such that the following holds.
		
		Let LDS be given. Let $ \kloc \geq 3 $. Suppose $ \G(x^{\sharp}, \kloc, M) \neq \void $, and $ f \geq \bmin\delta^2 $ on $\eloc$. Then there exists $ F \in C^2_+(\R^2) $ with $ F\big|_{\eloc} = f $, $ \norm{F}_{C^2(100Q)} \leq CM $, and $ \jet{x^{\sharp}}{F} \in \G(x^{\sharp}, \kloc, CM) $. 
	\end{lemma}

	\begin{proof}
		Without loss of generality, we may assume $ M = 1 $. 
		
		We will use $ b, B, B', $ etc. to denote quantities that depend only on $ \bmin $, and $ c, C, C', $ etc. to denote universal constants.
		
		Let $ P \in \G(x^{\sharp}, \kloc, 1) $. Pick distinct $ x_1, x_2 \in\eloc$. Let $ P^\sharp $ be the unique affine polynomial that passes through $ (x_1, f(x_1)), (x_2, f(x_2)) $, and $ (x^{\sharp}, P(x^{\sharp})) $. We first prove two claims about $ P^\sharp $.

		\begin{claim}\label{clam: Psharp}
			We have
			\begin{equation}
			\abs{ \nabla \left( P - P^\sharp \right)} \leq C\diam{ \mathrm{Triangle} (x_1, x_2, x^{\sharp})} \leq C\delta.
			\label{eq. Psh}
			\end{equation}
			As a consequence, $ P^\sharp \in \G(x^{\sharp}, \kloc, C) $.
		\end{claim}

		\begin{proof}[Proof of Claim \ref{clam: Psharp}]
			
			For convenience of notation, we temporarily label $ x_0 := x^{\sharp} $.
			
			Let $ S = \set{x_1, x_2} $. Since $ P \in \G(x_0, \kloc, 1) $ with $ \kloc \geq 3 $, there exists $ F^S \in C^2_+(\R^2) $ with $ F^S\big|_S = f $, $ \norm{F^S}_{C^2(\R^2)} \leq 1 $, and $ \jet{x_0}{F^S} = P $. In particular, $ F^S $ agrees with $ P^\sharp $ at $ x_i $ for $ i = 0, 1, 2 $.
			
			Let $ L_{ij} $ be the (open) segment connecting $ x_i $ and $ x_j $. The $ L_{ij} $'s are the sides of $ \mathrm{Triangle}(x_0,x_1,x_2) $. Let $ u_{ij} = \frac{x_j - x_i}{\abs{x_j - x_i}} $. Rolle's theorem implies that there exist $ \xi_{ij} \in L_{ij} $ such that 
			\begin{equation*}
			\nabla(F^S - P^\sharp)(\xi_{i j}) \cdot u_{ij} = 0
			\,\text{ for all } \,
			0 \leq i, j \leq 2.
			\end{equation*}
			Since $ \norm{F^S}_{C^2(\R^2)} \leq 1 $ and $ P^\sharp $ is an affine polynomial, we have
			\begin{equation}
			\abs{ \nabla (F^S - P^\sharp)(x_0) \cdot u_{ij} } \leq C\diam{\mathrm{Triangle}(x_0,x_1,x_2) }.
			\label{eq. FSPsh}
			\end{equation}
			Since $ \dist{x_0}{\etot} \geq \crep\delta $ and $ \eloc $ lies on the graph of $ \phi $ with $ \abs{\phi'} \leq 1 $, we have
			\begin{equation}
			\mathrm{Angle}\,(u_{01}, u_{12})
			\geq \gamma
			\label{eq. angle}
			\end{equation} 
			for some $ \gamma > 0 $ depending only on $ \crep $. 
			
			Let $ \omega $ be any unit vector. \eqref{eq. angle} implies that we can write
			\begin{equation}
			\omega = R_{\omega, 1}u_{01} + R_{\omega, 2}u_{12}
			,\,\,
			\abs{R_{\omega,i}} \leq C
			\text{ for } i = 1,2.
			\label{eq. omega}
			\end{equation}
			Here, $ C $ is a constant depending only on $ \gamma $. \eqref{eq. omega} implies that
			\begin{equation*}
			\abs{ \nabla (F^S - P^\sharp)(x_0) \cdot \omega } \leq C\diam{\mathrm{Triangle}(x_0,x_1,x_2)} \leq C\delta\,.
			\end{equation*}
			We conclude \eqref{eq. Psh} by letting $ \omega $ range over all unit vectors. Thanks to Lemma \ref{lem.perturb}, $ P^\sharp \in \G(x^{\sharp},\kloc,C) $. This proves the claim.
			
		\end{proof}

		\begin{claim}\label{claim: b}
			Suppose $ \bmin $ is sufficiently large. Then
			\begin{equation}
			P^\sharp (x^{\sharp}) \geq C\bmin\delta^2.
			\label{eq. Pbig}
			\end{equation}
		\end{claim}
		
		\begin{proof}[Proof of Claim \ref{claim: b}]
			
			The proof is identical to the proof of Claim \ref{claim.222}.
			
		\end{proof}

		Recall from LDS that $ \eloc $ lies on the graph of a $ C^2 $ function $ \phi $. Therefore, we may write $ \eloc = \set{z_i = (s_i, \phi(s_i)) : 1 \leq i \leq N }$ with $ s_i < s_{i+1} $ for all $ i = 1, \cdots,N-1 $. 
		
		For $ i = 1, \cdots, N-2 $, let $ S_i = \set{z_i, z_{i+1}, z_{i+2}} $. By Claim \ref{clam: Psharp}, $ P^\sharp \in \G(x^{\sharp}, \kloc, C) $. By definition, there exists $ F^{S_i} \in C^2_+(\R^2) $ such that $ F^{S_i}\big|_{S_i} = f $, $ \norm{F^{S_i}}_{C^2(\R^2)} \leq C $, and $ \jet{x^{\sharp}}{F^{S_i}} = P^\sharp $.
		
		Define $ g : \eloc \to \R $ by 
		\begin{equation*}
		g := f - \left(P^\sharp \big|_{\eloc}\right).
		\end{equation*}
		Note that $ g $ is not necessarily nonnegative. 
		
		Define $ G^{S_i} : \R^2 \to \R $ by
		\begin{equation*}
		G^{S_i} := F^{S_i} - P^\sharp.
		\end{equation*}
		Then immediately, we have
		\begin{align}
		&G^{S_i}\big|_{S_i} = g\, \text{ and }
		\label{eq. Gsi-int}
		\\
		&\jet{x^{\sharp}}{G^{S_i}} \equiv 0\,.
		\label{eq. Gsi-jet}
		\end{align}
		Since $ P^\sharp \in \G(x^{\sharp}, \kloc, C) $, we have $ \norm{P^\sharp}_{C^2(100Q)} \leq C $. 
		From this, together with the condition $ \norm{F^{S_i}}_{C^2(\R^2)} \leq C $, we learn that
		\begin{equation}
		\norm{G^{S_i}}_{C^2(100Q)} \leq C\,.
		\label{eq. Gsi-est}
		\end{equation}
		
		Thanks to \eqref{eq. Gsi-jet}, \eqref{eq. Gsi-est}, and the fundamental theorem of calculus, we have
		\begin{equation}
		\abs{\nabla G^{S_i}(s,t)} \leq C\delta
		\text{ for all } (s,t) \in 100Q.
		\label{eq. dtsmall}
		\end{equation}
		
		Let $ I_1 = (-\infty, s_3],I_2 = [s_2, s_4], \cdots, I_{N-3} = [s_{N-3}, s_{N-1}] $, and $ I_{N-2} = [s_{N-2}, +\infty) $. Let $ \set{\bar{\theta}_i : \R \to \R} $ be a partition of unity subordinate to $ \set{I_i}$ such that 
		\begin{equation}
		\supp{\bar{\theta}_i} \subset I_i
		\text{ and }
		\abs{\partial_s^k\bar{\theta}_i(s)} \leq \begin{cases}
		C\abs{s_{i+1} - s_i}^{-k}&\text{if } s \in [s_i, s_{i+1}]\\
		C\abs{s_{i+2}-s_{i+1}}^{-k}&\text{if } s \in [s_{i+1}, s_{i+2}]
		\end{cases} \text{ for } 1\leq k \leq 2
		\,.
		\label{eq. te}
		\end{equation}
		Note that the interior (in the topology of $ \R $) of $ I_i \cap I_j $ supports at most two partition functions. 
		
		Let
		\begin{equation*}
		\theta_i(s,t) := \bar{\theta}_i(s)
		\text{ for } i = 1, \cdots, N-2.
		\end{equation*}
		It follows immediately that the interior of $ (I_i\times\R) \cap (I_j \times \R) $ supports at most two partition functions. It is also clear that
		\begin{equation}
		\partial_t\theta_i \equiv 0
		\text{ for } i = 1, \cdots, N-2.
		\label{eq. dttheta0}
		\end{equation}
		
		Recall $ \Phi $ as in LDS. Define
		\begin{equation*}
		G :=  \brac{ \sum_{i = 1}^{N-2} \left(G^{S_i}\circ\Phi^{-1}\right) \bigg|_{t = 0} \cdot \theta_i } \circ \Phi = \sum_{i = 1}^{N-2} \left[G^{S_i}\circ\gamma\right] \cdot \left[\theta_i\circ\Phi\right],
		\end{equation*}
		where $ \gamma(s) := (s, \phi(s)) $ is a parametrization of the graph of $ \phi $. 
		
		\begin{claim}\label{claim: G}
			The function $ G $ satisfies $ G\big|_{\eloc} = g $ and $ \norm{G}_{C^2(100Q)} \leq C $.
		\end{claim}

		\begin{proof}[Proof of Claim \ref{claim: G}]
			It is clear from \eqref{eq. Gsi-int} that $ G\big|_{\eloc} = g $. 
			
			Now we estimate $ \norm{G}_{C^2(100Q)} $. 
			
			Thanks to \eqref{eq. dttheta0}, $ \supp{\theta_i \circ \Phi} \subset I_i \times \R $. Hence, the support of the $ \theta_i\circ\Phi $'s have bounded overlap. Since $ 0 \leq \theta_i \circ \Phi \leq 1 $, \eqref{eq. Gsi-est} implies that
			\begin{equation}
			\abs{G} \leq C
			\text{ on }
			100Q.
			\label{eq. Gsize}
			\end{equation}
			
			Now we compute the derivatives of $ G $. Thanks to \eqref{eq. dttheta0}, we have
			\begin{equation}
			\partial_t^k G \equiv 0
			\text{ for } k =1,2;
			\text{ and }
			\partial_{st}^2G \equiv 0
			\,.
			\label{eq. t0}
			\end{equation}
			Therefore, it remains to estimate the pure $ s $-derivatives of $ G $. 
			
			First of all, thanks to \eqref{eq. dttheta0}, we have
			\begin{equation*}
			\partial_s^k(\theta_i\circ\Phi)
			= \partial_s^k\theta_i\circ\Phi = \partial_s^k\bar{\theta}_i
			\text{ for } k = 1,2
			\,.
			\end{equation*}
			It follows from \eqref{eq. te} that
			\begin{equation}
			\abs{\partial_s^k(\theta_i \circ \Phi)(s,t)} \leq
			\begin{cases}
			C\abs{s_{i+1} - s_i}^{-k}
			&\text{if } (s,t)  \in [s_i, s_{i+1}]\times\R
			\\
			C\abs{s_{i+2}-s_{i+1}}^{-k}
			&\text{if } (s,t) \in [s_{i+1}, s_{i+2}]\times\R
			\end{cases} 
			\text{ for } 1\leq k \leq 2
			\,.
			\label{eq. tp}
			\end{equation}
			
			Now we compute the $ s $-derivatives of $ G^{S_i} \circ \gamma(s) $.
			\begin{equation*}
			\begin{split}
			\partial_s (G^{S_i} \circ \gamma) 
			&= 
			\partial_s G^{S_i} \circ \gamma + \phi'\partial_tG^{S_i}\circ\gamma\,,\\
			\partial_s^2(G^{S_i} \circ \gamma) 
			&= 
			\partial_s^2G^{S_i}\circ\gamma + (\phi')^2\partial_t^2G^{S_i}\circ\gamma + 2\phi'\partial_{st}^2G^{S_i}\circ\gamma + \phi''\partial_tG^{S_i}\circ\gamma\,.
			\end{split}
			\end{equation*}
			Recall that $ \abs{\phi^{(k)}} \leq \delta^{1-k} $ for $ k = 1,2 $. Applying \eqref{eq. dtsmall} to the last term of the second identity and \eqref{eq. Gsi-est} to the rest of the terms, we conclude that 
			\begin{equation}
			\norm{G^{S_i} \circ \gamma}_{C^2[-50,50]} \leq C.
			\label{eq. Giest}
			\end{equation}
			
			Since $ G(s,t) = G^{S_1}\circ\gamma(s) $ or $ G(s,t) = G^{S_{N-2}}\circ\gamma(s) $ outside of the strip $ [s_2, s_{N-1}]\times\R $, \eqref{eq. Gsi-est} implies
			\begin{equation}
			\abs{\partial_s^k G(s,t)} \leq C
			\text{ for }
			(s,t) \notin [s_2, s_{N-1}]\times\R, \, k = 1,2\,.
			\label{eq. oe}
			\end{equation}
			\renewcommand{\abs}[1]{\left|#1\right|}
			Suppose $ (s,t) \in [s_2, s_{N-1}]\times\R $. Let $ j $ be the least integer such that $ s \in I_j $. Then
			\begin{equation}
			\begin{split}
			\partial_s^k G = \left(\partial_s^k(G^{S_j}\circ\gamma)\right)\cdot\brac{\theta_j\circ\Phi} + \left(\partial_s^k(G^{S_{j+1}}\circ\gamma)\right)\cdot\brac{\theta_{j+1}\circ\Phi}  
			\\
			+\sum_{l = 0}^{k-l}\binom{k}{l}\brac{\partial_s^l(G^{S_j}\circ\gamma - G^{S_{j+1}}\circ\gamma)} \cdot \brac{\partial_s^{k-l}\theta_j\circ\Phi}.
			\end{split}
			\label{eq. G}
			\end{equation}
			By an argument similar to the proof of Claim \ref{claim 5}, combined with estimate \eqref{eq. Giest}, we have, for $ s \in I_j \cap I_{j+1} $,
			\begin{equation}
			\begin{split}
			\abs{\brac{G^{S_j}\circ\gamma - G^{S_{j+1}}\circ\gamma}(s)} &\leq C\abs{s_{j+1}-s_j}
			\text{ and }
			\\
			\abs{\partial_s^l \left( G^{S_j}\circ\gamma - G^{S_{j+1}}\circ\gamma \right)(s)} 
			&\leq C\abs{s_{j+1}-s_j}^{2-l}
			\text{ for } l = 0,1\,.
			\end{split}
			\label{eq. Gest}
			\end{equation}
			
			To estimate \eqref{eq. G}, we apply \eqref{eq. Giest} to the first two terms (note that $ 0 \leq \theta_j \circ \Phi \leq 1 $), and apply \eqref{eq. tp} and \eqref{eq. Gest} to the last term. Hence,
			\begin{equation}
			\abs{\partial^k_sG(s,t)} \leq C
			\text{ for }
			(s,t) \in [s_2, s_{N-1}]\times\R, \, k = 1,2\,.
			\label{eq. ie}
			\end{equation}
			The claim follows from \eqref{eq. Gsize}, \eqref{eq. t0}, \eqref{eq. oe}, and \eqref{eq. ie}.
			
		\end{proof}

		\newcommand{\clb}{C_{\mathrm{lb}}}
		
		Recall that, by construction, $ g(x_1) = g(x_2) = 0 $. Since $ G $ is constant in the $ t $-direction, Rolle's theorem implies that 
		\begin{equation}
		\abs{\d^\alpha G} \leq \clb\delta^{2-\abs{\alpha}}
		\,\text{ on }\,
		100Q \text{ for } \abs{\alpha} \leq 2\,.
		\label{eq. Gsmall}
		\end{equation}
		In particular, 
		\begin{equation}
		G \geq -\clb\delta^2\,\text{ on } 100Q.
		\label{eq. Gb}
		\end{equation}	
		
		Let $ \psi \in C^2_+(\R^2) $ be a cutoff function that satisfies the following. 
		\begin{itemize}
			\item $ 0 \leq \psi \leq 1 $ on $ \R^2 $, $ \psi \equiv 1 $ near $ x^{\sharp} $, $ \supp{\psi} \subset B(x^{\sharp}, \frac{\crep\delta}{100}) $; and
			\item $ \abs{\d^\alpha \psi} \leq C\delta^{2-\abs{\alpha}} $ for $ \abs{\alpha} \leq 2 $.
		\end{itemize}
		Define $ \tilde{G} := \R^2 \to \R $ by 
		\begin{equation*}
		\tilde{G} := (1 - \psi) G
		\,.
		\end{equation*}
		Then we have the following.
		\begin{itemize}
			\item Thanks to Claim \ref{claim: G}, \eqref{eq. Gsmall}, and the second property of $ \psi $, we have
			
		\end{itemize}
		\begin{equation}
		\norm{\tilde{G}}_{C^2(100Q)} \leq C\,.
		\label{eq. tG}
		\end{equation}
		\begin{itemize}

			\item $ \tilde{G} |_{\eloc} = g $, since $ \dist{x^{\sharp}}{\etot} \geq \crep\delta $ and $ \supp{\psi} \subset B(x^{\sharp}, \frac{\crep\delta}{100}) $.
			
			\item $ \jet{x^{\sharp}}{\tilde{G}} \equiv 0 $, since $ \psi \equiv 1 $ near $ x^{\sharp} $.

		\end{itemize}
		\begin{itemize}
			\item Moreover, since $ 0 \leq \psi \leq 1 $, \eqref{eq. Gb} implies
		\end{itemize}
		\begin{equation}
		\tilde{G} \geq -\clb\delta^2 \, \text{ on } 100Q. 
		\label{eq. R+2}
		\end{equation}
		
		Finally, define $ F : \R^2 \to \R $ by
		\begin{equation*}
		F := \tilde{G} + P^\sharp
		\,.
		\end{equation*}
		Then the following are immediate.
		\begin{itemize}
			\item $ F \big|_{\eloc} = g + \left(P^\sharp \big|_{\eloc}\right) = f $, and
			\item $ \jet{x^{\sharp}}{F} = \jet{x^{\sharp}}{\tilde{G}} + \jet{x^{\sharp}}{P^\sharp} = P^\sharp \in \G(x^{\sharp}, \kloc, C) $.
		\end{itemize}
		It remains to show that $ \norm{F}_{C^2(100Q)} $ is universally bounded and that $ F $ is nonnegative on $ 100Q $.
		
		Recall that $ P^\sharp \in \G(x^{\sharp}, \kloc, C) $, so we have $ \norm{P^\sharp}_{C^2(100Q)} \leq C $. It follows from \eqref{eq. tG} that
		\begin{equation*}
		\norm{F}_{C^2(100Q)} \leq C.
		\end{equation*}
		
		It remains to show that $ F $ is nonnegative on $ 100Q $.

		To this end, observe that \eqref{eq. TC1} and \eqref{eq. Pbig} imply
		\begin{equation}
		\dist{ x^{\sharp} }{\set{P^\sharp = 0}} \geq C\sqrt{P^\sharp (x^{\sharp})} = C\sqrt{\bmin}\cdot\delta.
		\label{eq. Pfar}
		\end{equation}
		Therefore, for sufficiently large $ \bmin $, \eqref{eq. Pbig} and \eqref{eq. Pfar} yield
		\begin{equation}
		P^\sharp \geq \clb\delta^2
		\,\text{ on } 100Q.
		\label{eq. R+1}
		\end{equation}
		Therefore, \eqref{eq. R+2} and \eqref{eq. R+1} imply that 
		\begin{equation*}
		F \geq 0
		\text{ on }
		100Q.
		\end{equation*}

		This concludes the proof of Lemma \ref{lem.FP-loc big}.

	\end{proof}

	Fix $ \bmin $ as in Lemma \ref{lem.FP-loc big}. The following lemma complements Lemma \ref{lem.FP-loc big}.
	
	\begin{lemma}\label{lem.FP-loc small}
		Let LDS be given. Let $ \kloc \geq 3 $.
		Suppose $ \G(x^{\sharp}, \kloc, M) \neq \void $, and that there exists $ x \in\eloc$ such that $ f(x) < \bmin M\delta^2 $. Then there exists $ F \in C^2_+(100Q) $ such that $ F \big|_{\eloc} = f $, $ \norm{F}_{C^2(100Q)} \leq BM $, and $ \jet{x^{\sharp}}{F} \in \G(x^{\sharp}, \kloc, BM) $. The number $ B $ depends only on $ \bmin $.
	\end{lemma}

	\begin{proof}
		Without loss of generality, we may assume $ M = 1 $. 
		
		We write $ B_1, B_2, $ etc. to denote quantities depending only on $ \bmin $.
		
		By Lemma \ref{lem.BS}, there exists $ \bmax > 0 $, depending only on $ \bmin $ such that 
		\begin{equation*}
		f(x) \leq \bmax\delta^2
		\,\text{ for all }
		x \in \eloc. 
		\end{equation*}
		By Lemma \ref{lem.perturb}, we have
		\begin{equation}
		0 \in \G(x^{\sharp}, \kloc, B_1).
		\label{eq. zj}
		\end{equation}
		
		Recall that $ \eloc $ lies on the graph of a $ C^2 $ function $ \phi $. Write $ \eloc = \set{z_i = (s_i, \phi(s_i)) : 1 \leq i \leq N }$ with $ s_i < s_{i+1} $ for all $ i = 1, \cdots,N-1 $. 
		
		For $ i = 1, \cdots, N-2 $, let $ S_i = \set{z_i, z_{i+1}, z_{i+2}} $. By \eqref{eq. zj}, there exists $ F^{S_i} \in C^2_+(\R^2) $ such that $ F^{S_i}\big|_{S_i} = f $, $ \norm{F^{S_i}}_{C^2(\R^2)} \leq B_4 $, and $ \jet{x^{\sharp}}{F^{S_i}} \equiv 0 $. 
		
		Let $ I_1 = (-\infty, s_3],I_2 = [s_2, s_4], \cdots, I_{N-3} = [s_{N-3}, s_{N-1}] $, and $ I_{N-2} = [s_{N-2}, +\infty) $. Let $ \set{\bar{\theta}_i} $ be a partition of unity subordinate to $ \set{I_i}$ such that 
		\begin{equation*}
		\supp{\bar{\theta}_i} \subset I_i
		\text{ and }
		\abs{\partial_s^k\bar{\theta}_i(s)} \leq \begin{cases}
		C\abs{s_{i+1} - s_i}^{-k}&\text{if } s \in [s_i, s_{i+1}]\\
		C\abs{s_{i+2}-s_{i+1}}^{-k}&\text{if } s \in [s_{i+1}, s_{i+2}]
		\end{cases} \text{ for } 1\leq k \leq 2
		\,.
		\end{equation*}
		Put
		\begin{equation*}
		\theta_i(s,t) := \bar{\theta}_i(s)
		\text{ for } i = 1, \cdots, N-2\,.
		\end{equation*}
		
		Recall the diffeomorphism $ \Phi $ in LDS. Define $ \bar{F}: \R^2 \to \R $ by
		\begin{equation*}
		\bar{F} := \brac{ \sum_{i = 1}^{N-2} \left(F^{S_i}\circ\Phi^{-1}\right) \bigg|_{t = 0} \cdot \theta_i } \circ \Phi
		\,.
		\end{equation*}
		It is clear that $ \bar{F} \geq 0 $. By the same argument as in the proof of Claim \ref{claim: G}, we have $ \bar{F}|_{\eloc} = f $ and $ \norm{\bar{F}}_{C^2(100Q)} \leq
		B_2 $.
		
		Since $ f \leq \bmax\delta^2 $ on $ \eloc $ and $ \bar{F} $ is constant in the $ t $-direction, we also have
		\begin{equation*}
		\abs{\bar{F}} \leq C\bmax\delta^2
		\text{ on }
		[-50,50).
		\end{equation*}
		
		Let $ \psi \in C^2_+(\R^2) $ be a cutoff function such that 
		\begin{itemize}
			\item $ 0 \leq \psi \leq 1 $ on $ \R^2 $, $ \psi \equiv 1 $ near $ x^{\sharp} $, $ \supp{\psi} \subset B(x^{\sharp}, \frac{\crep\delta}{100}) $; and
			\item $ \abs{\d^\alpha\psi} \leq C\delta^{2-\abs{\alpha}} $ for $ \abs{\alpha} \leq 2 $. 
		\end{itemize}
		
		Define $ F : \R^2 \to \R $ by
		\begin{equation*}
		F := (1 - \psi) \bar{F}\,. 
		\end{equation*}
		
		The following hold.
		\begin{itemize}
			\item $ F \geq 0 $, since $ \bar{F} \geq 0 $ and $ 0 \leq \psi \leq 1 $.
			
			\item $ F\big|_{\eloc} = f $, thanks to \eqref{eq. xsharp} and the fact that $ \supp{\psi} \subset B(x^{\sharp}, \frac{\crep\delta}{100}) $.
			
			\item $ 
			\jet{x^{\sharp}}{F} \equiv 0 \in \G(x^{\sharp}, \kloc, B_4) $, since $ \psi \equiv 1 $ near $ x^{\sharp} $.
			
			\item $ \norm{F}_{C^2(100Q)} \leq B_3 $. To see this, we note that since $ F \geq 0 $ on $ \R^2 $, $ \norm{F}_{C^2(100Q)} \leq B_2 $, and $ \abs{F} \leq C\bmax\delta^2 $ on $ 100Q $, \eqref{eq. TC1} implies that $ \abs{\nabla F} \leq B_4\delta $ on $ 100Q $. Thanks to the second condition on $ \psi $, the conclusion follows.
		\end{itemize}
		
		This proves the lemma.
	\end{proof}

	\begin{proof}[Proof of Lemma \ref{lem.FP-loc}]
		Fix $ \bmin $ as in Lemma \ref{lem.FP-loc big}. The lemma follows from Lemma \ref{lem.FP-loc big} and Lemma \ref{lem.FP-loc small}. 
	\end{proof}

	\subsection{Proof of Theorem \ref{thm.2DFP}}
	\label{subsection.proof-FP}

	\newcommand{\kknice}{\Lambda^{(4)}_{\mathrm{nice}}}
	
	Before proceeding to the proof of Theorem \ref{thm.2DFP}, we make a brief comment on the finiteness constant $ 64 $. Lemma \ref{lemma: Gamma trade off} and Lemma \ref{lemma: neighboring polynomials} state that jets of $ 4k $-point interpolants based in neighboring squares from $ \knice $ are compatible in the Whitney sense (see (\ref{eq: neighboring polynomials})); Lemma \ref{lem.graph} states that the geometry of data points in each square of $ \knice $ is sufficiently nice when $ k \geq 4 $; Lemma \ref{lem.FP-loc} states that in such case, a local version of the extension problem is readily solved. Hence, if we pick $ k = 4 $, we may use the jets of $ 4 \cdot 4 = 16 $-point interpolants (if they exist) to guarantee compatibility of nearby local extensions. By Lemma \ref{lemma: 17}, such jets exist.

	Now, we examine compatibility of the local interpolants constructed in Lemma \ref{lem.FP-loc}.

	\begin{proof}[Proof of Theorem \ref{thm.2DFP}]
		Without loss of generality, we may assume $ M = 1 $. 
		
		Set $ k = 4 $. Pick $ \cnice $ so that $ (4, \cnice) $ guarantees good geometry. 
		
		By Lemma \ref{lemma: CZ} $ \kknice $ is a CZ covering of $ \R^2 $. 
		
		By Lemma \ref{lemma: 17}, $ \G(x, 16, 1) \neq \void $ for any $ x \in \R^2 $. 
		
		We distinguish three types of squares in $ \knice $. 
		
		\begin{enumerate}[\text{Type }1.]
			\item Suppose $E \cap Q^* \neq \void$. Let $ F^\sharp_Q := F_{Q} $, where $ F_{Q} $ is as in Lemma \ref{lem.FP-loc} with $ \kloc = 16 $. Let $ P^\sharp_Q := \jet{x^{\sharp}_Q}{F^\sharp_Q} $. We have $ P^\sharp_Q \in \G(x^{\sharp}_Q, 16, C) $.

			\item Suppose $ E \cap Q^* = \void $ but $\delta_Q < 1 $. Pick $ P_Q^\sharp \in \G(x^\sharp_Q, 16, 1) $, and set $ F_Q^\sharp := \ew^{\{x_Q^\sharp\}} (P_Q^\sharp) $, where $ \ew^{\{x_Q^\sharp\}} $ is as in Lemma \ref{lem.WE-map} with $ S = \{x_Q^\sharp\} $.

			\item Suppose $ E \cap Q^* = \void $ and $\delta_Q = 1$. Set $F^\sharp_Q \equiv  0$.
		\end{enumerate}

		By Lemma \ref{lem.FP-loc} and Lemma \ref{lem.WE-map}, $ F^\sharp_Q \in C^2_+(100Q) $, $ F^\sharp_Q\big|_{E \cap Q^*} = f $, and
		\begin{equation}\label{eq: F_i est}
		\norm{F^\sharp_Q}_{C^2(100Q)} \leq C. 
		\end{equation}
		
		\begin{claim}\label{claim: neighboring est}
			If $Q \leftrightarrow Q'$, then for each $x \in \frac{9}{8}Q \cup \frac{9}{8}Q'$ and $ 0 \leq \abs{\alpha} \leq 1 $,
			\begin{equation}\label{eq. neighboring est}
			\abs{ \partial^\alpha (F^\sharp_Q - F^\sharp_{Q'})(x)  } \leq C\delta_Q^{2 - \abs{\alpha}}.
			\end{equation}
			The constant $ C $ is universal. 
		\end{claim}
		
		\begin{proof}[Proof of Claim \ref{claim: neighboring est}]
			Temporarily fix $ x \in \frac{9}{8}Q \cup \frac{9}{8}Q' $ for $ Q \leftrightarrow Q' $. 
			
			Assume that either $Q$ or $Q'$ is of Type 3, then (\ref{eq. neighboring est}) follows from \eqref{eq.CZ-def} and (\ref{eq: F_i est}). 
			
			Suppose neither $ Q $ nor $ Q'$ is of Type 3. Thanks to \eqref{eq.CZ-def} and our choice of $ x^{\sharp}_Q $ and $ x^\sharp_{Q'} $ in Lemma \ref{lem.rep}, we have $ \abs{x^{\sharp}_Q - x}, \abs{x^{\sharp}_{Q'} - x}, \abs{x^{\sharp}_Q - x^{\sharp}_{Q'}} \leq C\delta_Q $. 
			
			Recall from Lemma \ref{lem.FP-loc} and Lemma \ref{lem.WE-map} that 
			\begin{equation*}
			 \jet{x^{\sharp}_Q}{F^\sharp_Q} = P^\sharp_Q \in \G(x^{\sharp}_Q, 16, C) 
			 \text{ and }
			 \jet{x^{\sharp}_{Q'}}{F^\sharp_{Q'}} = P^\sharp_{Q'} \in \G(x^{\sharp}_{Q'}, 16, C).
			\end{equation*}

			By Taylor's theorem, 
			\begin{equation}\label{eq: taylor 1}
			\begin{split}
			\abs{ \partial^\alpha  \left( F^\sharp_Q - P^\sharp_Q \right) (x) } \leq C\delta_Q^{2 - \abs{\alpha}} \text{ and }
			\abs{ \partial^\alpha  \left( F^\sharp_{Q'} - P^\sharp_{Q'} \right) (x) } \leq C\delta_{Q'}^{2 - \abs{\alpha}} \leq C\delta_Q^{2-\abs{\alpha}}\,.
			\end{split}
			\end{equation}
			
			By Lemma \ref{lemma: neighboring polynomials}, 
			\begin{equation}\label{eq: nbp}
			\abs{\d^\alpha(P^\sharp_Q - P^\sharp_{Q'})(x)} \leq C\delta_Q^{2 - \abs{\alpha}}. 
			\end{equation}
			
			Now, (\ref{eq. neighboring est}) follows from (\ref{eq: taylor 1}) and (\ref{eq: nbp}).
			
		\end{proof}

		Let $ \set{\theta_Q} $ be a partition of unity that is CZ-compatible with $\kknice$. Define
		\begin{equation*}
		F(x) := \sum_{Q \in \kknice} \theta_Q(x)\cdot F^\sharp_Q(x)
		\end{equation*}
		It is clear that $ F \geq 0 $, $ F|_E = f $, and $ F $ is twice continuously differentiable. 
		For $\abs{\alpha} \leq 2$ and $x \in Q$,
		\begin{equation}\label{eq: exp 1}
		\partial^\alpha F(x) = \sum_{Q \in \kknice}\partial^\alpha F^\sharp_Q(x)\cdot \theta_Q(x) + \sum_{Q' \leftrightarrow Q} \sum_{0 < {\beta} \leq {\alpha}}\binom{\alpha}{\beta}\partial^{\alpha - \beta}(F^\sharp_{Q} - F^\sharp_{Q'})(x)\cdot \partial^{\beta}\theta_{Q'}(x).
		\end{equation}
		Applying (\ref{def: BIP}), (\ref{def: partition of 1 }), (\ref{eq: F_i est}), and (\ref{eq. neighboring est}) to (\ref{eq: exp 1}), we can conclude that
		\begin{equation*}
		\norm{F}_{C^2(\R^2)} \leq \csharp.
		\end{equation*}
		
	\end{proof}

	\section{Sharp Finiteness Principle}
	\label{sect.SFP}

	In this section, we give the proof of Theorem \ref{thm.SFP}. Here we remind the readers the statement of the theorem. 
	
	\namedtheorem[Theorem 5]{\label{thm5}Let $ E \subset \R^2 $ with $ \#(E) = N < \infty $. Then there exist universal constants $ C, C' ,C'' $ and a list of subsets $ S_1, S_2, \cdots, S_L \subset E $ satisfying the following.
		\begin{enumerate}[(A)]
			\item $ \#(S_\ell) \leq C  $ for each $ \ell = 1, \cdots, L $.
			\item $ L \leq C'N $.
			\item Given any $ f : E \to [0,\infty) $, we have
			\begin{equation*}
			\max_{\ell = 1, \cdots, L}\norm{f}_{C^2_+(S_\ell)} \leq \norm{f}_{C^2_+(E)} \leq C'' \max_{\ell = 1, \cdots, L}\norm{f}_{C^2_+(S_\ell)}\,.
			\end{equation*}
	\end{enumerate} }{2-D Sharp Finiteness Principle}
	
	Before we proceed to the proof, we briefly explain the clusters $ S_\ell $'s in the statement. 
	
	For each square $ Q \in \knice $, we associate to it a basic cluster $ S(\xqs) $ (see Definition \ref{def.Sx}) that guarantees internal Whitney compatibility.
	
	The clusters in Theorem \hyperref[thm5]{5} can be classified into three types.
	
	\begin{itemize}
		\item The first type is the union of a ``consecutive'' three-point cluster (since $ E $ locally lies on a curve with controlled geometry), nearby basic clusters, and nearby ``keystone'' clusters (see next bullet point). This is the ``largest'' type of clusters, since it plays the key role of relaying information about $ E $ to various lengthscales.

		\item The second type is the basic cluster for each ``keystone square'' (see Definition \ref{def.CZ2}). Keystone squares are locally the smallest squares and they play an important role in relaying information to nearby small squares containing no data point.
		
		\item The third type is the union of keystone square clusters (see the second bullet point above) that are associated with each ``special square'' (see Lemma \ref{lem.Lsp}). This type of clusters is used to eliminate the ambiguity in how these special squares receive information from $ E $.
		
	\end{itemize}

	We now give the full account.\\
	\medspace

	\newcommand{\Ls}{\Lambda^\sharp}
	\newcommand{\Lz}{\Lambda_0}
	\renewcommand{\knice}{\Lambda_0}
	\newcommand{\dq}{\delta_Q}
	\newcommand{\Lks}{\Lambda_{\rm KS}}
	
	\subsection{CZ squares and clusters}
	
	Let $ (k,\cnice) $ guarantee good geometry (Definition \ref{def.good geo}). We fix such $ (k,\cnice) $ for the rest of the section. We may assume, for instance, 
	\begin{equation*}
	k = 4
	\text{ and }
	\cnice = 1000\,.
	\end{equation*}

	\begin{definition}\label{def.CZ2}
		We define the following objects.
		\begin{itemize}
			\item We set
			\begin{equation}
			\eqindent
			\Lz := \Lambda^{(k)}_{\rm nice}.
			\label{eq.Lz-def}
			\end{equation}
			\item We also set
			\begin{equation}
			\eqindent
			\Ls:= \set{Q \in \knice : E \cap Q^* \neq \void}\,.
			\label{eq.Ls-def}
			\end{equation}
			Note that $ \Ls $ coincides with Type 1 squares in the proof of Theorem \ref{thm.2DFP} (Section \ref{subsection.proof-FP}). 
			
			\item We say $ Q \in \Lz $ is a \underline{keystone square} if $ \dq < 1 $ and for any $ Q' \in \Lz $ with $ Q' \cap 100Q \neq \void $, we have $ \delta_{Q'} \geq \dq $. The collection of keystone squares is denoted by $ \Lks $. 
		\end{itemize}
	\end{definition}

	Keystone squares first appear in the work of Sobolev extension\cite{I13}. See also \cite{FIL14} for a more thorough discussion.

	\begin{lemma}\label{lem.keystone}
		Let $ \Lks $ be as in Definition \ref{def.CZ2}. Then
		\begin{equation*}
		\#(\Lks) \leq C\cdot \#(E).
		\end{equation*}
%
%
%
		
	\end{lemma}

	The proof of Lemma \ref{lem.keystone} can be found in Section 4 of \cite{I13} and Section 7 of \cite{FIL14}.

	Next, we define the basic cluster associated with each square in $ \Lz $. 
	
	\begin{definition}\label{def.Sx}
		Let $ Q \in \Lz $ and let $ \xqs $ be as in Lemma \ref{lem.rep}. (Note that $ \xqs $ is a representative point ``far'' from the data on the lengthscale $ \delta_Q $.) Let $ S_1, \cdots, S_{12} \subset E $ be as in Lemma \ref{lem.sigma-localized} (with $ x = \xqs $ and $ 4k $ in place of $ k $). We define
		\begin{equation}
		S(\xqs) :=  \bigcup_{i = 1}^{12} S_i\,.
		\label{eq.Sx-def}
		\end{equation}
	\end{definition}

	Since $ \#(S_i) \leq C $ for $ i = 1, \cdots, 12 $ (Lemma \ref{lem.sigma-localized}), we have
	\begin{equation}
	\#(S(\xqs)) \leq C\,.
	\label{eq.Sxqs-bd}
	\end{equation}
	By Lemma \ref{lem.sigma-localized}, we have
	\begin{equation}
	\sigma(\xqs,S(\xqs)) \subset C\cdot \ss(\xqs,4k)\,.
	\label{eq.lalala}
	\end{equation}

	Next, we state a key lemma that allows us to relay information from keystone squares to small squares in $ \Lz $ whose neighborhood contains no points from $ E $. The latter requires separate attention for the following reason: Suppose $ Q \in \Lz $ with $ \dq < 1 $ and $ E \cap Q^* = \void $. Then $ (Q^+)^* $ may intersect an uncontrolled number of squares in $ \Ls $. Keystone squares are designed partially to deal with such situations. See \cite{I13,FIL14} for further discussion.

	\newcommand{\Lsp}{\Lambda_{special}}
	
	\begin{lemma}\label{lem.Lsp}
		Let $ \Lz, \Lks $ be as in Definition \ref{def.CZ2}. We can find a subset $ \Lsp \subset \Lz $ and a map $ \mu : \Lz \to \Lks $ such that the following holds for some universal constant $ C $.
		\begin{enumerate}[(A)]
			\item $ \#(\Lsp) \leq C\cdot \#(E) $.
			
			\item $ \mu(Q) \in \Lks $, where $ \Lks $ is as in Definition \ref{def.CZ2}. Moreover, $ \dist{Q}{\mu(Q)} \leq C\dq $.
			
			\item Suppose $ Q, Q' \in \Lz \setminus \Lsp $ and $ Q \leftrightarrow Q' $, then $ \mu(Q) = \mu(Q') $. 
		\end{enumerate}
	\end{lemma}
	
	The proof of Lemma \ref{lem.Lsp} can be found in Section 6 of \cite{FIL14}.

	\newcommand{\touch}{\leftrightarrow}
	
	\newcommand{\Sks}{S_{\rm KS}}
	\newcommand{\Ssp}{S_{\rm special}}
	
	\begin{definition}\label{def.KS-cluster}
		Recall $ \Lz, \Ls, \Lks $ as in Definition \ref{def.CZ2}. Recall the representative point $ \xqs $ as in Lemma \ref{lem.rep}. Let $ Q \in \Lks $. We define
		\begin{equation}
		\Sks(Q) := S(x_{Q}^\sharp),
		\label{eq.Sks}
		\end{equation} 
		where $ S(\xqs) $ is as in \eqref{eq.Sx-def}. Recall $ \Lsp, \mu $ as in Lemma \ref{lem.Lsp}. Let $ Q \in \Lsp $. We define
		\begin{equation}
		\eqindent
		\Ssp(Q) := \bigcup_{Q' \leftrightarrow Q, Q' \in \Lz}S(x_{\mu(Q')}^\sharp),
		\label{eq.Ssp}
		\end{equation}
		where $ x_{\mu(Q')}^\sharp $ is as in Lemma \ref{lem.rep} and $ S(x_{\mu(Q')}^\sharp) $ is as in \eqref{eq.Sx-def}.

	\end{definition}

	Recall from Lemma \ref{lemma: CZ} that $ \Lz $ is a CZ covering of $ \R^2 $. In particular, $ \Lz $ enjoys the bounded intersection property \eqref{def: BIP}. Together with \eqref{eq.Sxqs-bd} and the definitions of $ \Sks, \Ssp $ in \eqref{eq.Sks},\eqref{eq.Ssp}, we see that
	\begin{align}
	\#(\Sks(Q)) &\leq C
	\text{ for each } Q \in \Lks, \text{ and }\label{eq.Sks-bound}\\
	\#(\Ssp(Q)) &\leq C
	\text{ for each } Q \in \Lsp.
	\label{eq.Ssp-bound}
	\end{align}

	Now we turn our attention to clusters associated with $ \Ls $.

	\newcommand{\Sb}{\overline{S}}
	For convenience, we set, for each $ Q \in \Ls $, 
	\begin{equation*}
	N(Q) := \#(E \cap Q^*)\,.
	\end{equation*}
	Thanks to Lemma \ref{lem.graph}, we know that for each $ Q \in \Ls $, up to a rotation, $ E \cap Q^* \subset \set{(s,\phi(s)) : s \in \R} $, where $ \phi $ is as in Lemma \ref{lem.graph}. We enumerate 
	\begin{equation}
		E \cap Q^* = \set{(s_i, \phi(s_i)): i = 1, \cdots, N(Q)} 
		\enskip\text{ such that }
		s_1 < \cdots < s_{N(Q)}\,.
		\label{eq.enu}
	\end{equation}
	Let $ \Phi $ be as in Lemma \ref{lem.diffeo}. We also set
	\begin{equation}
	I_Q := \Phi(100Q)\big|_{\R\times\set{t = 0}}\,.
	\label{eq.IQ-def}
	\end{equation}
	For the rest of this section, whenever we consider $ Q \in \Ls $, we always assume that $ Q $ has been rotated so that enumeration of the form \eqref{eq.enu} holds.

	The next three definitions describe the objects of interest in this section. Definitions \ref{def.Sb} and \ref{def.SQ} concern the clusters, and Definition \ref{def.KQ} concerns the main polynomial convex sets.

	\begin{definition}\label{def.Sb}
		Let $ Q \in \Ls $. Let $ E \cap Q^* $ be enumerated as in \eqref{eq.enu}. 
		\begin{itemize}

			\item In the case $ 1 \leq N(Q) \leq 2 $, we set
			\begin{equation}
			\displayindent0pt
			\displaywidth\textwidth
			\Sb(Q,1) := E \cap Q^*
			\enskip
			\text{ and }
			\nu(Q):=1\,.
			\label{eq.Sb-def-2}
			\end{equation}

			\item Suppose $ N(Q) \geq 3 $, we set
			\begin{equation}
			\displayindent0pt
			\displaywidth\textwidth
			\Sb(Q,\nu) := \set{ (s_\nu,\phi(s_\nu)), (s_{\nu+1},\phi(s_{\nu+1})), (s_{\nu+2},\phi(s_{\nu+2})) },
			\label{eq.Sb-def}
			\end{equation}
			for 
			$ \nu = 1, \cdots, \nu(Q) $, 
			where
			$ \nu(Q) := N(Q)-2 $.

		\end{itemize}
	\end{definition}

	By the bounded intersection property of $ \Lz $ (see \eqref{def: BIP}), we have
	\begin{equation}
	\#\set{  \Sb (Q,\nu) : Q\in \Ls, \nu \in \set{1,\cdots,\nu(Q)} } \leq C\cdot \#(E)\,.
	\label{eq.Sb-bounded}
	\end{equation}

	\newcommand{\K}{\mathcal{K}}

	\begin{definition}\label{def.SQ}
		Let $ Q \in \Ls $. Let $ \Sb(Q,\nu) $ be as in Definition \ref{def.Sb}. Let $ \xqs $ be as in Lemma \ref{lem.rep}. Let $ \mu $ be the map in Lemma \ref{lem.Lsp}. Let $ S(\cdot) $ be as in \eqref{eq.Sx-def}. For each $ \nu ={1, \cdots, \nu(Q)} $, we set
		\begin{equation}
		S(Q,\nu) := \Sb(Q,\nu) \cup  \bigcup_{\substack{Q' \leftrightarrow Q\\Q'\in\Lz}}\brac{S(x_{Q'}^\sharp)\cup S(x_{\mu(Q')}^\sharp)}
		\label{eq.SQ-def}
		\end{equation}
	\end{definition}
	
	\begin{remark}
		The cluster $ S(Q,\nu) $ associated with each $ Q \in \Ls $ is the ``largest'' among all three types of clusters (the other two being $ \Sks(Q) $ in \eqref{eq.Sks} and $ \Ssp(Q) $ in \eqref{eq.Ssp}). This is as expected, since each $ Q \in \Ls $ satisfies $ E \cap Q^* \neq \void $, and must relay information to neighboring squares and their keystone representatives. 
	\end{remark}

	Thanks to \eqref{eq.Sxqs-bd}, the fact that $ \#(\Sb(Q,\nu)) \leq 3 $, and the bounded intersection property of $ \Lz $ (see \eqref{def: BIP}), we have
	\begin{equation}
	\#(S(Q,\nu)) \leq C
	\text{ for each }
	\nu = 1, \cdots, {\nu(Q)}\,.
	\label{eq.SQ bounded}
	\end{equation}
	Thanks to \eqref{eq.Sb-bounded}, we have
	\begin{equation}
	\#\set{  S (Q,\nu) : Q\in \Ls, \nu \in \set{1,\cdots,\nu(Q)} } \leq C\cdot \#(E)\,.
	\label{eq.Sb-B}
	\end{equation}

	To distinguish the roles of the clusters related to $ \Ls $, we make the following definition.
	
	\begin{definition}\label{def.KQ}
		Let $ Q \in \Ls $ and $ M \geq 0 $. Let $ \xqs $ be as in Lemma \ref{lem.rep}. Let $ S(Q,\nu) $ be as in Definition \ref{def.SQ}. We define
		\begin{equation}
		\K(Q,\nu,M) := \Gamma_+(\xqs,S(Q,\nu),M).
		\label{eq.KQ def}
		\end{equation}
	\end{definition}

	\subsection{Whitney compatibility}

	The next lemma is similar to Lemma \ref{lemma: neighboring polynomials}.

	\begin{lemma}\label{lem.S-S'}
		There exists a universal constant $ C $ such that the following holds.
		Let $ Q, Q' \in \Lz $. Let $ x_Q^\sharp, x_{Q'}^\sharp $ be as in Lemma \ref{lem.rep}. Let $ S(\xqs) $ be as in \eqref{eq.Sx-def}. Let $ S, S' \subset E $. Suppose
		\begin{equation}
		S(\xqs)\subset \brac{S \cap S'}.
		 \label{eq.ScupS'}
		\end{equation}
		Then given $ P \in \Gamma_+(\xqs,S,M) $ and $ P' \in \Gamma_+(x_{Q'}^\sharp,S',M) $, we have
		\begin{equation}
		\abs{\d^\alpha (P - P')(\xqs)} \leq CM\brac{ \dq + \abs{x_Q^\sharp - x_{Q'}^\sharp}  }^{2-\abs{\alpha}}
		\text{ for }\abs{\alpha} \leq 2.
		\label{eq.P-P'}
		\end{equation}
	\end{lemma}

	\begin{proof}
		\newcommand{\xqsp}{x_{Q'}^\sharp}
		Fix $ P $ and $ P' $ as in the hypothesis. By definition, there exist $ F, F' \in C^2_+(\R^2) $ such that the following hold.
		\begin{itemize}
			\item $ F|_S = f $ and $ F'|_{S'} = f $.
			\item $ \norm{F}_{C^2(\R^2)} \leq M $ and $ \norm{F'}_{C^2(\R^2)} \leq M $. 
			\item $ \jet{\xqs}{F} = P $ and $ \jet{x_{Q'}^\sharp}{F'} = P' $.
		\end{itemize}
		Thanks to \eqref{eq.ScupS'}, we see that
		\begin{equation*}
		F - F' = 0
		\text{ on }
		S(\xqs).
		\end{equation*}
		By the definition of $ \sigma $ in Section \ref{section: convex}, we see that
		\begin{equation*}
		\jet{\xqs}{(F-F')} = P - \jet{\xqs}{F'} \in 2M\cdot \sigma(\xqs,S(\xqs)).
		\end{equation*}
		By Lemma \ref{lem.sk-ball} and the definition of $ S(\xqs) $ in \eqref{eq.Sx-def}, we see that
		\begin{equation}
		\abs{\d^\alpha (P - \jet{\xqs}{F'})(\xqs) } \leq CM\dq^{2-\abs{\alpha}}.
		\label{eq.pj}
		\end{equation}
		By the triangle inequality, we have
		\begin{equation*}
		\abs{\d^\alpha (P - P')(\xqs)} \leq \abs{\d^\alpha(P - \jet{\xqs}{F'})(\xqs)} + \abs{\d^\alpha (\jet{\xqs}{F'} - \jet{\xqsp}{F'})(\xqs)}.
		\end{equation*}
		Using \eqref{eq.pj} to estimate the first term, and using Taylor's theorem to estimate the second, we see that \eqref{eq.P-P'} holds.

	\end{proof}

	\begin{remark}
		We note that Lemma \ref{lem.S-S'} is a one-sided estimate, in the sense that the right hand side of \eqref{eq.P-P'} does not contain the lengthscale $ \delta_{Q'} $. However, this is remedied once we know that $ Q \touch Q' $. This is further examined in the next corollary, which states that suitable choices of clusters give rise to Whitney compatible jets.
	\end{remark}

	\begin{corollary}\label{lem.K-K}
		\newcommand{\xqsp}{{x_{Q'}^\sharp}}
		There exists a universal constant $ C $ such that the following holds. Let $ \Lz,\Ls,\Lks $ be as in Definition \ref{def.CZ2}. Let $ \Lsp $ and $ \mu $ be as in Lemma \ref{lem.Lsp}. For $ Q \in \Lz $, let $ \xqs $ be as in Lemma \ref{lem.rep}, and let $ S(\xqs) $ be as in \eqref{eq.Sx-def}. Suppose $ Q, Q'\in \Lz $ with $ Q\leftrightarrow Q' $ and $ P, P' \in \P $ satisfy one of the following conditions.
		\begin{enumerate}[(A)]
			\item Suppose $ Q , Q' \in \Ls $. Let $ \nu ={1, \cdots, \nu(Q)} $ and $ \nu' = {1, \cdots, \nu(Q')} $ (Definition \ref{def.SQ}). Let $ \K(\cdot,\cdot,\cdot) $ be as in Definition \ref{def.KQ}. Suppose $ P \in \K(Q,\nu,M) $ and $ P' \in \K(Q',\nu',M) $.
			
			\item Suppose $ Q \in \Ls $ and $ Q' \in \Lz \setminus (\Ls \cup \Lsp) $. Suppose $ P \in \K(Q,\nu,M) $ (Definition \ref{def.KQ}) for some $ \nu \in \set{1, \cdots, \nu(Q)} $ (Definition \ref{def.SQ}). Suppose $ P' \in \Gamma_+(x_{\mu(Q')}^\sharp, \Sks(\mu(Q')),M) $, with $ \Sks(\mu(Q')) $ as in \eqref{eq.Sks}.
			
			\item Suppose $ Q \in \Ls $ and $ Q' \in \Lsp \setminus \Ls $.  Suppose $ P \in \K(Q,\nu,M) $ (Definition \ref{def.KQ}) for some $ \nu \in \set{1, \cdots, \nu(Q)} $ (Definition \ref{def.SQ}). Suppose $ P' \in \Gamma_+(x_{Q'}^\sharp, \Ssp(Q'),M) $, with $ \Ssp(Q') $ as in \eqref{eq.Ssp}. 
			
			\item Suppose $ Q, Q' \in \Lz \setminus (\Ls \cup \Lsp) $. Suppose $ P \in \Gamma_+(x_{\mu(Q)}^\sharp, \Sks(\mu(Q)),M) $ and suppose $ P' \in \Gamma_+(x_{\mu(Q')}^\sharp, \Sks(\mu(Q')),M) $, with $ \Sks(\mu(Q)) $ and $ \Sks(\mu(Q')) $ as in \eqref{eq.Sks}.
			
			\item Suppose $ Q \in \Lz \setminus (\Ls \cup \Lsp)  $ and $ Q' \in \Lsp\setminus \Ls $. Suppose $ P \in \Gamma_+(x_{\mu(Q)}^\sharp, \Sks(\mu(Q)),M) $, with $ \Sks(\mu(Q)) $ as in \eqref{eq.Sks}. Suppose $ P' \in \Gamma_+(x_{Q'}^\sharp,\Ssp(Q'),M) $, with $ \Ssp(Q') $ as in \eqref{eq.Ssp}. 
			
			\item Suppose $ Q, Q' \in \Lsp\setminus \Ls $. Suppose $ P \in \Gamma_+(x_{Q}^\sharp,\Ssp(Q),M) $ and $ P' \in \Gamma_+(x_{Q'}^\sharp,\Ssp(Q'),M) $, with $ \Ssp(Q) $ and $ \Ssp(Q') $ as in \eqref{eq.Ssp}.
		\end{enumerate}
	
		Then
		\begin{equation*}
		\abs{\d^\alpha (P - P')(\xqs)}, \,\,	\abs{\d^\alpha (P - P')(x_{Q'}^\sharp)} \leq CM\dq^{2-\abs{\alpha}}
		\text{ for }
		\abs{\alpha} \leq 2.
		\end{equation*}
		
	\end{corollary}

	\begin{proof}
		Thanks to Lemma \ref{lemma: CZ} and Lemma \ref{lem.Lsp}, we know that
		\begin{equation*}
		\abs{\xqs - x_{Q'}^\sharp},\, \abs{\xqs- x_{\mu(Q)}^\sharp},\, \abs{x_{\mu(Q)}^\sharp- x_{\mu(Q')}^\sharp} \leq C\dq
		\text{ for } Q, Q' \in \Lz
		\text{ with }
		Q \leftrightarrow Q'
		\end{equation*}
		Therefore, by Lemma \ref{lem.S-S'} and Taylor's theorem, it suffices to show that in (A)-(F), the sets $ S, S' $ in $ \Gamma_+(x_\star^\sharp,S,M) \ni P $, $ \Gamma_+(x_{\star'}^\sharp,S',M)\ni P' $ satisfy
		\begin{equation}
		S(x_\star^\sharp) \subset S\cap S',
		\text{ for some } x_{\star}^\sharp \in \set{\xqs,x_{Q'}^\sharp, x_{\mu(Q)}^\sharp, x_{\mu(Q')}^\sharp}.
		\label{eq.star}
		\end{equation}

		We analyze each scenario.
		
		\begin{enumerate}[(A)]
			
			\item Recall from \eqref{eq.KQ def} that $ \K(Q,\nu,M) = \Gamma_+(\xqs,S(Q,\nu),M) $ and $ \K(Q',\nu',M) = \Gamma_+(x_{Q'}^\sharp,S(Q',\nu'),M) $. In this scenario, $ S = S(Q,\nu) $ and $ S' = S(Q',\nu') $. We let $ x_\star^\sharp = x_Q^\sharp $.
			
			We see from \eqref{eq.SQ-def} that $ S(\xqs) \subset S(Q,\nu) $ and $ S(\xqs) \subset S(Q',\nu') $, since $ Q \touch Q' $. Therefore, $ S(\xqs) \subset S(Q,\nu) \cap S(Q',\nu') $. \eqref{eq.star} follows.
			
			\item Recall from \eqref{eq.KQ def} that $ \K(Q,\nu,M) = \Gamma_+(\xqs,S(Q,\nu),M) $. In this scenario, $ S = S(Q,\nu) $ and $ S' = \Sks(\mu(Q')) $. We let $ x_\star^\sharp = x_{\mu(Q')}^\sharp $.
			
			We see from \eqref{eq.SQ-def} that $ S(x_{\mu(Q')}^\sharp) \subset S(Q,\nu) $, since $ Q\touch Q' $. Recall from \eqref{eq.Sks} that $  \Sks(\mu(Q')) = S(x_{\mu(Q')}^\sharp)$. \eqref{eq.star} follows.
			
			\item Recall from \eqref{eq.KQ def} that $ \K(Q,\nu,M) = \Gamma_+(\xqs,S(Q,\nu),M) $. Thus, in this scenario, $ S = S(Q,\nu) $ and $ S' = \Ssp(Q') $. We let $ x_\star^\sharp = x_{\mu(Q)}^\sharp $. 
			
			We see from \eqref{eq.SQ-def} that $ S(x_{\mu(Q)}^\sharp) \subset S(Q,\nu) $, since $ Q \touch Q $ by definition. We see from \eqref{eq.Ssp} that $ S(x_{\mu(Q)}^\sharp) \subset \Ssp(Q') $ since $ Q\touch Q' $. \eqref{eq.star} follows.

			\item In the current scenario, $ S = \Sks(\mu(Q)) $ and $ S' = \Sks(\mu(Q')) $. 
			
			By Lemma \ref{lem.Lsp}, we have $ \mu(Q) = \mu(Q') $. Hence, $ S = S' $. Taking $ x_\star^\sharp = x_{\mu(Q)}^\sharp $, we see from \eqref{eq.Sks} that $ S(x_{\mu(Q)}^\sharp) = S \cap S' $. \eqref{eq.star} follows.

			\item In the current scenario, $ S = \Sks(\mu(Q)) $ and $ S' = \Ssp(Q') $. Let $ x_\star^\sharp = x_{\mu(Q)}^\sharp $. 
			
			Recall from \eqref{eq.Sks} that $ S(x_{\mu(Q)}^\sharp) = \Sks(\mu(Q)) $. From \eqref{eq.Ssp}, we see that $ S(x_{\mu(Q)}^\sharp) \subset \Ssp(Q') $, since $ Q'\touch Q $. \eqref{eq.star} follows. 
			
			\item In this scenario, $ S = \Ssp(Q) $ and $ S' = \Ssp(Q') $. We let $ x_\star^\sharp = x_{\mu(Q)}^\sharp $. 
			
			By \eqref{eq.Ssp}, $ S(x_{\mu(Q)}^\sharp) \subset \Ssp(Q) $, since $ Q \touch Q $ by definition. By \eqref{eq.Ssp} again, $ S(x_{\mu(Q)}^\sharp) \subset \Ssp(Q') $, since $ Q'\touch Q $. \eqref{eq.star} follows.

		\end{enumerate}
		
		We have exhausted all the cases. This concludes the proof of the corollary. 
		
	\end{proof}

	\subsection{Local extension problem}

	The next lemma states that on the correct local scale, the two-dimensional trace norm behaves in a similar way as the one-dimensional trace norm.

	\begin{lemma}\label{lem.1D-rep}
		Let $ Q \in \Ls $ and let $ \phi $ be as in Lemma \ref{lem.graph}. Let $ \Sb \subset E \cap Q^* $. Recall the definition of $ I_Q $ in \eqref{eq.IQ-def}. There exists a universal constant $ C $ such that the following hold.

		\begin{enumerate}[(A)]
			\item Let $ f : \Sb \to [0,\infty) $. Suppose there exists $ F \in C^2_+(100Q) $ with $ F = f $ on $\Sb $, and $ \abs{\d^\alpha F}\leq M\delta_Q^{2-\abs{\alpha}} $ on $ 100Q $ for $ \abs{\alpha} \leq 2 $. Then there exists $ \bar{F}_\nu \in C^2_+(I_Q) $ with $ \bar{F}_\nu(s) = f(s,\phi(s)) $ for each $ (s,\phi(s)) \in \Sb $, and $ \abs{\d_s^k \bar{F}_\nu} \leq CM\delta_Q^{2-k} $ on $ I_Q $ for $ k \leq 2 $. 
			
			\item Let $ g : \Sb \to \R $. Suppose there exists $ G \in C^2(100Q) $ with $ G = g $ on $ \Sb$, and $ \abs{\d^\alpha G}\leq M\delta_Q^{2-\abs{\alpha}} $ on $ 100Q $ for $ \abs{\alpha} \leq 2 $. Then there exists $ \bar{G} \in C^2(I_Q) $ with $ \bar{G}(s) = g(s,\phi(s)) $ for each $ (s,\phi(s)) \in \Sb $, and $ \abs{\d_s^k \bar{G}} \leq CM\delta_Q^{2-\abs{\alpha}} $ on $ I_Q $ for $ k \leq 2 $. 
		\end{enumerate}

	\end{lemma}

	\begin{proof}
		We only prove (A) here. The proof for (B) is identical.
		
		Let $ \Phi $ be as in Lemma \ref{lem.diffeo}, and let $ \Psi = (\Psi_1, \Psi_2) := \Phi^{-1} $. Let $ F $ be as in the hypothesis. Consider the function 
		\begin{equation*}
		\bar{F}(s) := F \circ \Psi (s,0)\,.
		\end{equation*}
		Since $ F \geq 0 $, we have $ \bar{F} \geq 0 $. By Lemma \ref{lem.diffeo}, we have $ \bar{F}(s) = f(s,\phi(s)) $ for each $ (s,\phi(s)) \in \Sb $. It remains to estimate the derivatives for $ \bar{F} $. Setting $ \d_1 = \d_s $ and $ \d_2 = \d_t $, we have
		\begin{equation*}
		\begin{split}
		\d_i(F \circ \Psi) &= \sum_{k = 1}^2 \d_i\Psi_k\cdot \d_kF\circ\Psi,\enskip\text{ and }\\
		\d_{ij}(F \circ \Psi) &= \sum_{k,l = 1}^{2} C_{kl} \d_i\Psi_k \cdot \d_j \Psi_l + \d_{kl} F \circ \Psi + \sum_{k = 1}^2\d_{ij}\Psi_k \cdot \d_k F\circ\Psi \,.
		\end{split}
		\end{equation*}
		Therefore, thanks to Lemma \ref{lem.diffeo} and the hypothesis $ \abs{\d^\alpha F} \leq M\delta_Q^{2-\abs{\alpha}} $, we can conclude that $ \abs{\d_s^k \bar{F}} \leq CM\delta_Q^{2-k} $ on $ I_Q $ for $ k \leq 2 $. This concludes the proof of the lemma.
	\end{proof}

	We can think of the next lemma as a re-scaled local finiteness principle (without a prescribed jet). It is essentially a consequence of Theorem \ref{thm.FP-1d}. 
	
	\renewcommand{\hat}[1]{\widehat{#1}}
	
	\begin{lemma}\label{lem.nojet}
		Let $ Q \in \Ls $. For each $ \nu = 1, \cdots, \nu(Q) $, let $ \Sb(Q,\nu) $ be as in Definition \ref{def.Sb}. 
		
		\begin{enumerate}[(A)]
			\item Let $ f : E\cap Q^* \to [0,\infty) $. Suppose for each $ \nu $, there exists $ F_\nu \in C^2_+(100Q) $ such that $ F_\nu = f $ on $ \Sb(Q,\nu) $, and $ \abs{\d^\alpha F_\nu} \leq M\delta_Q^{2-\abs{\alpha}} $. Then there exist a universal constant $ C $ and a function $ \hat{F}_Q \in C^2_+(\R^2) $ such that
			\begin{enumerate}[(i)]
				\item $ \hat{F}_Q\big|_{E \cap Q^*} = f  $, and
				\item $ \abs{\d^\alpha \hat{F}_Q} \leq CM\delta_Q^{2-\abs{\alpha}} $ on $ 100Q $,  $ \abs{\alpha} \leq 2 $. 
			\end{enumerate}
			
			\item Let $ g : E\cap Q^* \to \R $. Suppose for each $ \nu $, there exists $ G_\nu \in C^2(100Q) $ such that $ G_\nu = g $ on $ \Sb(Q,\nu) $, and $ \abs{\d^\alpha G_\nu} \leq M\delta_Q^{2-\abs{\alpha}} $. Then there exist a universal constant $ C $ and a function $ \hat{G}_Q \in C^2(\R^2) $ such that
			\begin{enumerate}[(i)]
				\item $ \hat{G}_Q\big|_{E \cap Q^*} = g  $, and
				\item $ \abs{\d^\alpha \hat{G}_Q} \leq CM\delta_Q^{2-\abs{\alpha}} $ on $ 100Q $,  $ \abs{\alpha} \leq 2 $. 
			\end{enumerate}
		\end{enumerate}

	\end{lemma}

	\renewcommand{\bar}[1]{\overline{#1}}
	
	\begin{proof}
		
		We only prove (A) here. The proof for (B) is identical.
		
		If $ \#(E \cap Q^*) \leq 3 $, then $ \nu(Q) = 1 $ and $ \Sb(Q,\nu(Q)) = E \cap Q^* $, and the conclusions follow directly from the definition of $ \norm{f}_{{C}^2_+(\Sb(Q,\nu(Q)))} $. For the rest of the proof, we assume $ \#(E \cap Q^*) > 3 $. 
		
		Up to a rotation, we know that $ E \cap Q^* \subset \set{(s,\phi(s)) : s \in \R} $, where $ \phi $ is as in Lemma \ref{lem.graph}. Enumerate $ E \cap Q^* $ as in \eqref{eq.enu}. For $ \nu = 1, \cdots, N(Q) -2 $, we set
		\begin{equation}
		I_\nu := [s_{\nu}, s_{\nu+2}]. 
		\label{eq.In}
		\end{equation}
		We also set 
		\begin{equation}
		I_0 := (-\infty,s_2]
		\enskip
		\text{ and }
		\enskip
		I_{N(Q)-1}:= [s_{N(Q)-1}, \infty)\,.
		\label{eq.In2}
		\end{equation}

		Let $ \set{\bar{\theta}_\nu}_{\nu = 1}^{N(Q)-1} $ be a partition of unity subordinate to the cover $ \set{I_\nu}_{\nu = 1}^{N(Q)-1} $, such that
		\begin{equation}
		\abs{\d_s^k\bar{\theta}_\nu(s)} \leq\begin{cases}
		C\abs{s_\nu - s_{\nu-1}}^{-k} &\text{ if } s \in [s_{\nu-1}, s_{\nu}]\\
		C\abs{s_{\nu+1}-s_\nu}^{-k} &\text{ if } s \in [s_\nu,s_{\nu+1}]
		\end{cases}
		\enskip
		\text{ for }
		k = 0,1,2\,.
		\label{eq.theta-n}
		\end{equation}
		Here it is convenient to use $ s_0 := -\infty $, $ s_{N(Q) + 1} = \infty $, and $ \infty^0 = 1 $. We set
		\begin{equation}
		\theta_\nu(s,t) := \bar{\theta}_\nu(s)
		\enskip
		\text{ for }
		\nu = 0, 1, \cdots, N(Q)-1 \,.
		\label{eq.theta-n2}
		\end{equation}

		Let $ \bar{F}_\nu $ be as in Lemma \ref{lem.1D-rep} with $ \Sb = \Sb(Q,\nu) $ for $ \nu = 1, \cdots, N(Q) - 2 $.  By Rolle's Theorem, we have
		\begin{equation}
		\abs{\d_s^k(\bar{F}_\nu - \bar{F}_{\nu+1})} \leq CM\delta_Q^{2-k}\abs{s_{\nu+1}-s_\nu}^{2-k}
		\label{eq.F-n}
		\end{equation}
		for $ \nu = 1, \cdots, N(Q) - 2 $, $ s \in I_\nu\cap I_{\nu+1} $, and $ k \leq 2 $.
		
		We also set 
		$ \bar{F}_0 := \bar{F}_1 $ and $ \bar{F}_{N(Q) - 1} := \bar{F}_{N(Q)-2} $.
		
		Define
		\begin{equation*}
		F_\nu(s,t):= \bar{F}_\nu(s)
		\enskip
		\text{ for }
		\nu = 0,1, \cdots, N(Q) - 1\,.
		\end{equation*}
		Finally, we set
		\begin{equation*}
		\hat{F}_Q(x) := \sum_{\nu = 0}^{N(Q) - 1} \theta_\nu(x) \cdot  F_\nu(x)\,.
		\end{equation*}		
		It is clear that $ \hat{F}_Q \geq 0 $ on $ \R^2 $ and $ \hat{F}_Q = f $ on $ E \cap Q^* $. By construction, $ \d_t \theta_\nu = \d_t F_\nu \equiv 0 $ for each $ \nu = 0, \cdots, N(Q) - 1 $. Then, using estimates \eqref{eq.theta-n} and \eqref{eq.F-n}, we can conclude that $ \abs{\d^\alpha \hat{F}_Q} \leq CM\delta_Q^{2-\abs{\alpha}}$ on $ 100Q $ for $ \abs{\alpha} \leq 2 $.
		
	\end{proof}

	Repeating the proof of Lemma \ref{lem.BS} and using Lemma \ref{lem.nojet}(A), we have the following result tailored for the matter at hand. 
	
	\begin{lemma}\label{lem.BS-new}
		For each $ \bmin > 0 $ sufficiently large, we can find $ \bmax$, depending only on $ \bmin $, such that the following holds. Let $ Q \in \Ls $, and let $ \K(Q,\nu,M) $ be as in Definition \ref{def.KQ}. Suppose for each $ \nu = 1, \cdots, \nu(Q) $, $ \K(Q,\nu,M) \neq \void $. Then at least one of the following holds.
		\begin{enumerate}[(A)]
			
			\item $ f(x) \geq \bmin M\delta_Q^2 $ for all $ x \in E \cap Q^* $. 
			
			\item $ f(x) \leq \bmax M\delta_Q^2 $ for all $ x \in E \cap Q^* $.
			
		\end{enumerate}
		
	\end{lemma}

	\begin{proof}
		Suppose (A) holds. There is nothing to prove.
		
		\newcommand{\hn}{\hat{\nu}}
		\newcommand{\hx}{\hat{x}}
		Suppose (A) fails. We write $ B_0 = \bmin $ and we fix the number $ B_0 $ throughout. 
		
		By assumption, there exists $ \hx \in E \cap Q^* $ with $ f(\hx) < \bmin \delta_Q^2 $. There exists $ \hn \in \set{1,\cdots, N(Q)-2} $ such that $ \hat{x} \in \Sb(Q,\hn) \subset S(Q,\hn) $. By assumption, $ K(Q,\hat{\nu},M) \neq \void $, so there exists $ \hat{F} \in C^2_+(\R^2) $ with $ \hat{F}= f $ on $ S(Q,\hn) $, $ \norm{\hat{F}}_{C^2(\R^2)} \leq M $, and $ \jet{\xqs}{\hat{F}} \in \K(Q,\hn,M) $.  By Lemma \ref{lem.TC}, we have
		\begin{equation*}
		\abs{\nabla F(\hx)} \leq C\sqrt{M F(\hx)} \leq CM\sqrt{B_0}\delta_Q\,.
		\end{equation*}
		Then Taylor's theorem implies
		\begin{equation*}
		F(\xqs) \leq CM\brac{B_0\delta_Q^2+\sqrt{B_0}\,\abs{\xqs - \hx}\delta_Q} \leq CM(\sqrt{B_0} + 1)^2\delta_Q^2\,.
		\end{equation*}
		
		Let $ x_0 \in E \cap Q^* $. Then there exists $ \nu(x_0) \in \set{1, \cdots N(Q)} $ such that $ x' \in S(Q,\nu(x_0)) $. 
		
		By assumption, $ \K(Q,\nu(x_0),M) \neq \void $. Pick $ P \in \K(Q,\nu(x_0),M) $. By Corollary \ref{lem.K-K}, we see that
		\begin{equation*}
		\abs{\jet{\xqs}{\hat{F}} - P(\xqs)} \leq CM\delta_Q^2\,.
		\end{equation*}
		Therefore, we have $ P(\xqs) \leq CM(\sqrt{B_0} + 1)^2\delta_Q^2 $. By the definition of $ \K(Q,\nu(x_0),M) $, there exists $ F \in C^2_+(\R^2) $ with $ F(x_0) = f(x_0) $ and $ \jet{\xqs}{F} = P $. In particular, by Lemma \ref{lem.TC} and Taylor's Theorem, we have
		\begin{equation*}
		\abs{\nabla F(x)} \leq CM(\sqrt{B_0} + 1)\delta_Q
		\text{ for all } x \in Q^*\,.
		\end{equation*}
		By Taylor's theorem again, we have
		\begin{equation*}
		F(x_0) \leq CM(\sqrt{B_0} + 1)\delta_Q^2\,.
		\end{equation*}
		Since $ x_0 \in E \cap Q^* $ was chosen arbitrarily, (B) follows.

	\end{proof}

	The next lemma mirrors Lemma \ref{lem.perturb}. It says the following. When the local data is big, $ \K $ can be viewed as a translate of $ \sk $. When the local data is small, $ \K $ contains not much more information than the zero jet.

	\begin{lemma}\label{lem.perturb-new}
		
		Let $ Q \in \Ls $. Let $ \K(Q,\nu,M) $ be as in Definition \ref{def.KQ}. Suppose $ \K(Q,\nu,M) \neq \void $ for each $ \nu = 1, \cdots, \nu(Q) $. 
		\begin{enumerate}[(A)]
			\item There exists a number $ B > 0 $ exceeding a universal constant such that the following holds. Suppose $ f(x) \geq B M\delta_Q^2 $ for all $ x \in E \cap Q^* $. Then $ \K(Q,\nu,M) + M\cdot\sk(\xqs,4k) \subset \K(Q,\nu,CM) $ for each $ \nu ={1, \cdots, \nu(Q)} $. Here, $ C $ is a universal constant.
			\item Let $ A > 0 $. Suppose $ f(x) \leq A M\delta_Q^2 $ for all $ x \in E \cap Q^* $. Then $ 0 \in \K(Q,\nu,A'M) $ for each $ \nu \in \set{1,\cdots,\nu(Q)} $. Here the number $ A' $ depends only on $ A $. 
		\end{enumerate}
		
	\end{lemma}

	\begin{proof}
		We adapt the proof of Lemma \ref{lem.perturb} with $ \K $ in place of $ \G $, and use the fact that $ S(Q,\nu) $ contains $ \Sb(Q,\nu) \subset E \cap Q^* $. We include the relevant steps here for completeness.
		
		Fix $ \nu \in \set{1,\cdots, \nu(Q)} $.
		
		We begin with (A). Let $ B > 0 $ be a sufficiently large number. 
		
		By \eqref{eq.Sx-def}, we have $ \Sb(Q,\nu(Q)) \subset S(Q,\nu) $. Let $ P \in \K(Q,\nu,M) $. Repeating the proof of Claim \ref{claim.222} in Lemma \ref{lem.perturb}, we see that $ P(\xqs) \geq C(\sqrt{B}-1/2)^2M\delta_Q^2 $.

		By \eqref{eq.Gamma+def}, Lemma \ref{lem.TC-W}, and Definition \ref{def.KQ}, there exists a Whitney field 
		\begin{equation*}
		 \vec{P} = (P,(P^x)_{x \in S(Q,\nu)}) \in W^2_+(\{\xqs\}\cup S(Q,\nu)) 
		\end{equation*}
		such that $ P^x = f(x) $ for all $ x \in S(Q,\nu) $, and $ \norm{\vec{P}}_{W^2_+(\{\xqs\}\cup S(Q,\nu))} \leq CM $. 
		
		Let $ \tilde{P} \in M\cdot\sk(\xqs,4k) $. By Lemma \ref{lem.sk-ball}, $ \tilde{P} \in CM\cdot \B(x,\delta_Q)$. 
		
		Consider the Whitney field 
		\begin{equation*}
		\vec{P}':= (P + \tilde{P}, (P^x)_{x \in S(Q,\nu)}) \,.
		\end{equation*}
		By the same argument in the proof of Lemma \ref{lem.perturb}, we can verify that 
		\begin{equation*}
		\vec{P}' \in W^2_+(\{\xqs\}\cup S(Q,\nu)) 
		\enskip
		\text{ and }
		\norm{\vec{P}'}_{W^2_+(\{\xqs\}\cup S(Q,\nu))} \leq CM \,.
		\end{equation*}
		Part (A) then follows from Lemma \ref{lem.WE-map}. 
		
		Now we turn to (B). 
		
		Let $ P \in \K(Q,\nu,M) $. Repeating the argument in Claim \ref{Claim.Psmall} in the proof of Lemma \ref{lem.perturb}, we have $ P(\xqs) \leq C(\sqrt{A}+1)^2M\delta_Q^2  $. By the definitions of $ \K $ and $ \Gamma_+ $ in \eqref{eq.KQ def} and \eqref{eq.Gamma+def}, there exists $ F \in C^2_+(\R^2) $ with $ F(x) =f(x) $ for all $ x \in S(Q,\nu) $, $ \norm{F}_{C^2_+(\R^2)} \leq CM $, and $ \jet{\xqs}{F} = P $. By Lemma \ref{lem.TC} and Taylor's theorem, we have
		\begin{equation}
		\abs{\d^\alpha F(x)} \leq CA''M\delta_Q^{2-\abs{\alpha}}
		\text{ for all } x \in Q^*
		\text{ and }\abs{\alpha} \leq 2\,.
		\label{eq.FQ-small}
		\end{equation}
		Here, $ A'' $ depends only on $ A $. 
		
		Let $ \psi\in C^2_+(\R^2) $ be a cutoff function such that $ \psi \equiv 1 $ near $ \xqs $, $ \psi \equiv 0 $ outside of $ B(\xqs,\frac{\crep\delta_Q}{100}) $, and $ \abs{\d^\alpha\psi} \leq C\delta_Q^{-\abs{\alpha}} $. 
		
		We set 
		\begin{equation*}
		\tilde{F} := (1-\psi)\cdot F\,.
		\end{equation*}
		It is clear that $ F \geq 0 $ on $ \R^2 $, $ \tilde{F} = f $ on $ S(Q,\nu) $, and $ \jet{\xqs}{\tilde{F}} \equiv 0 $. Using \eqref{eq.FQ-small}, we see that $ \norm{\tilde{F}}_{C^2(\R^2)} \leq A'M $.  This proves part (B) and concludes the proof of the lemma.
		
	\end{proof}

	The next lemma mirrors Lemma \ref{lem.FP-loc}. It solves the local interpolation with a prescribed jet in $ \K $, so that they can be patched together by a partition of unity.

	\begin{lemma}\label{lem.FP-loc-S}
		Let $ Q \in \Ls $. Let $ \K(Q,\nu,M) $ be as in Definition \ref{def.KQ}. Suppose $ \K(Q,\nu, M) \neq \void $ for each $ \nu = 1, \cdots, \nu(Q) $. Then there exist a universal constant $ C $ and a function $ F_Q \in C^2_+(100Q) $ such that
		\begin{enumerate}[(A)]
			\item $ F_Q\big|_{E \cap Q^*} = f $, 
			
			\item $ \norm{F_Q}_{C^2(\R^2)} \leq CM $, and
			
			\item $ \jet{\xqs}{F_Q} \in \K(Q,\nu(Q),CM) $. 
		\end{enumerate}
	\end{lemma}

	\begin{proof}[Proof]
		We adapt the proof of Lemma \ref{lem.FP-loc} with the following main adjustments:
		\begin{itemize}
			\item We use Lemma \ref{lem.perturb-new} in place of Lemma \ref{lem.perturb}. 
			
			\item We use Lemma \ref{lem.BS-new} in place of Lemme \ref{lem.BS}.
			
			\item We use $ \K $ in place of $ \G $, and the condition $ \G(\xqs,4k,M) \neq \void $ is replaced by $ \K(Q,\nu,M) \neq \void $ for each $ \nu = 1, \cdots, \nu(Q) $. See Lemma \ref{lem.nojet} and Lemma \ref{lem.FP-loc-S}. 
		\end{itemize}
		
		Here we present the relevant steps for completeness.
		
		Fix $ Q \in \Ls $.
		
		Suppose $ \#(E \cap Q^*) \leq 3  $. Recall Definitions \ref{def.Gamma}, \ref{def.Sb}, \ref{def.SQ}, and \ref{def.KQ}. By assumption, $ \K(Q,\nu(Q),M) = \Gamma_+(\xqs,S(Q,\nu(Q)),M) \neq \void $. Pick $ P \in \Gamma_+(\xqs,S(Q,\nu(Q)),M) $. By the definition of $ \Gamma_+ $, there exists $ F_{Q} \in C^2_+(\R^2) $ such that $ F_{Q}\big|_{S(Q,\nu(Q))} = f $, $ \norm{F_{Q}}_{C^2(\R^2)} \leq M $, and $ \jet{\xqs}{F_Q} \in \Gamma_+(\xqs,S(Q,\nu(Q)),M) $. Since $ S(Q,\nu) \supset \Sb(Q,\nu(Q)) $ and $ \Sb(Q,\nu(Q)) = \Sb(Q,1) = E \cap Q^* $ in this case, the conclusions follow.
		
		From now on, we assume $ \#(E \cap Q^*) > 3 $. 
		
		Let $ \bmin > 0 $ be sufficiently large, and in particular, $ \bmin > B $, where $ B $ is as in Lemma \ref{lem.perturb-new}. Let $ \bmax $ be given as in Lemma \ref{lem.BS} with such $ \bmin $. 
		
		Thanks to Lemma \ref{lem.BS-new}, each $ Q \in \Ls $ falls into at least one of the following cases. 
		
		\begin{enumerate}[(i)]
			\item $ f(x) \geq \bmin M\delta_Q^2 $ for all $ x \in E \cap Q^* $. 
			
			\item $ f(x) \leq \bmax M\delta_Q^2 $ for all $ x \in E \cap Q^* $. 
		\end{enumerate}
		
		We treat (i) first. 
		
		Since $ \#(E\cap Q^*) \geq 3 $, we may select distinct $ x_1, x_2 \in E \cap S(Q,\nu(Q)) \cap Q^*  $. Pick $ P \in \K(Q,\nu(Q),M) $. Let $ P^\sharp $ be the unique affine polynomial that interpolates the points $ (x_1, f(x_1)) $, $ (x_2, f(x_2)) $, and $ (\xqs,P(\xqs)) $. We may repeat the argument in Claim \ref{clam: Psharp} in the proof of Lemma \ref{lem.FP-loc big} and use Lemma \ref{lem.perturb-new} to show that 
		\begin{equation*}
		P^\sharp \in \K(Q,,\nu(Q),CM)
		\enskip
		\text{ and }
		\enskip
		P^\sharp(\xqs) \geq C\bmin M \delta_Q^2\,.
		\end{equation*}
		This, together with Lemma \ref{lem.TC}, implies that 
		\begin{equation*}
		\dist{\xqs}{\set{P^\sharp = 0}} \geq C\sqrt{\bmin}\delta_Q \,.
		\end{equation*}
		Therefore, we have
		\begin{equation}
		P^\sharp (x) \geq C\bmin M\delta_Q^2
		\text{ for all } x \in 100Q\,.
		\label{eq.Psh-big}
		\end{equation}
		
		Let $ g(x) := f(x) - P^\sharp(x) $ for each $ x \in E \cap Q^* $. Note that $ g $ is not necessarily nonnegative. Since $ P^\sharp \in \K(Q,\nu(Q), CM) $, there exists a function $ F \in C^2_+(\R^2) $ such that $ F|_{S(Q,\nu)} = f $, $ \norm{F}_{C^2(\R^2)} \leq CM $, and $ \jet{\xqs}{F} = P $. This, together with the assumption $ \K(Q,\nu,M) \neq \void $ and Rolle's theorem, implies that for each $ \nu \in 1, \cdots, \nu(Q) $, there exists $ G_\nu \in C^2(\R^2) $ such that 
		\begin{equation}
		G_\nu = g  \text{ on }  \Sb(Q,\nu) \enskip \text{ and } \enskip
		 \abs{\d^\alpha G(x)} \leq CM\delta_Q^{2-\abs{\alpha}}\enskip
		 \text{ for all } x \in 100Q\,,\, \abs{\alpha} \leq 2\,.
		 \label{eq.Gn-est}
		\end{equation}
		By Lemma \ref{lem.nojet}(B), there exists $ G \in C^2(100Q) $ such that
		\begin{equation}
		G|_{E \cap Q^*} = g
		\enskip\text{ and }\enskip \abs{\d^\alpha G(x)} \leq CM\delta_Q^{2-m}
		\enskip\text{ for all } x \in 100Q\,,\, \abs{\alpha} \leq 2\,.
		\label{eq.G-small}
		\end{equation}

		Let $ \psi \in C^2_+(\R^2) $ be a cutoff function such that 
		\begin{equation}
		\text{$ \psi \equiv 1 $ near $ \xqs $, $ \psi \equiv 0 $ outside of $ B(\xqs,\frac{\crep \delta_Q}{100}) $, and $ \abs{\d^\alpha\psi} \leq C\delta_Q^{2-\abs{\alpha}} $.}
		\label{eq.sigh}
		\end{equation}
		
		Consider the function 
		\begin{equation*}
		F_Q := P^\sharp + (1-\psi)\cdot G\,.
		\end{equation*}
		
		\begin{itemize}
			\item By \eqref{eq.Psh-big} and \eqref{eq.G-small}, we have $ F_Q \geq 0 $ on $ 100Q $. 
			\item Since $ \supp{\psi}$ is disjoint from $ E \cap Q^* $, we have $ F_Q(x) = P^\sharp(x)  +g(x) = f(x) $ for each $ x \in E \cap Q^* $. (A) is satisfied.

			\item Since $ P^\sharp \in \K(Q,\nu(Q), CM) $, we have
			\begin{equation*}
			\norm{P^\sharp}_{C^2(100Q)} \leq CM\,.
			\end{equation*}
			By \eqref{eq.G-small} and \eqref{eq.sigh}, we have
			\begin{equation*}
			\norm{(1-\psi)\cdot G}_{C^2(100Q)} \leq CM\,.
			\end{equation*}
			Conclusion (B) then follows from the triangle inequality.
			
			\item 	Since $ \psi \equiv 1 $ near $ \xqs $, we have $ \jet{\xqs}{F_Q} = \jet{\xqs}{P^\sharp} + 0 = P^\sharp \in \K(Q,\nu(Q),CM) $. (C) is satisfied. 
		\end{itemize}

		 This proves case (i).

		Now we turn to case (ii).

		Recall the hypothesis $ \norm{f}_{C^2_+(\Sb(Q_\nu))} \leq M $ for each $ \nu $. By definition, for each $ \nu = 1, \cdots, N(Q) $, there exists $ F_\nu \in C^2_+(\R^2) $ such that $ F_\nu = f $ on $ \Sb(Q,\nu) $ and $ \norm{F_\nu}_{C^2(\R^2)} \leq CM $. Since $ f(x) \leq \bmax\delta_Q^2 $ for all $ x \in E \cap Q^* $, by Lemma \ref{lem.TC}, we have
		\begin{equation*}
		\abs{\d^\alpha F_\nu(x)} \leq CM\delta_Q^{2-\abs{\alpha}}
		\enskip \text{ for all } x \in 100Q\,.
		\end{equation*}
		Therefore, the hypotheses of Lemma \ref{lem.nojet}(A) are satisfied, and there exists $ F \in C^2_+(\R^2) $ such that $ F|_{E \cap Q^*} = f $ and 
		\begin{equation}
		\abs{\d^\alpha F(x)} \leq CM\delta_Q^{2-\abs{\alpha}}
		\enskip \text{ for all } x \in 100Q\,.
		\label{eq.F-small}
		\end{equation}

		Let $ \psi $ satisfy \eqref{eq.sigh}. Consider the function
		\begin{equation*}
		F_Q := (1-\psi) \cdot F\,.
		\end{equation*}
		
		\begin{itemize}
			\item Since $ F \geq 0 $ and $ 0 \leq \psi \leq 1 $, we have $ F_Q \geq 0 $ on $ 100Q $. 
			\item Since $ \supp{\psi} $ is away from $ E \cap Q^* $, we have $ F_Q = f $ on $ E \cap Q^* $. (A) is satisfied.

			\item Thanks to \eqref{eq.sigh} and \eqref{eq.F-small}, we have $ \norm{F_Q}_{C^2(100Q)} \leq CM $. (B) is satisfied.
			
			\item 	Since $ \psi \equiv 1 $ near $ \xqs $, we have $ \jet{\xqs}{F_Q} \equiv 0 \in \K(Q,\nu(Q), CM) $, thanks to Lemma \ref{lem.perturb-new}. (C) is satisfied.
		\end{itemize}

		This concludes the treatment for case (ii) and the proof of the lemma.

	\end{proof}

	\subsection{Proof of Theorem \ref{thm.SFP}}

	Now we define the $ S_\ell $'s in the statement of Theorem \ref{thm.SFP}. 
	
	\newcommand{\Lss}{\mathcal{S}^{\sharp}}

	\begin{definition}
		Recall the definitions of $ \Lz,\Ls,\Lks,\Lsp $ in Definition \ref{def.CZ2} and Lemma \ref{lem.Lsp}. We set
		\begin{equation}
		\Lss := \Lss_1 \cup \Lss_2 \cup \Lss_3,
		\label{eq.Lss}
		\end{equation}
		where
		\begin{itemize}
			\item $ \Lss_1 := \set{S(Q,1), \cdots, S(Q,\nu(Q)) \enskip : \enskip Q \in \Ls
			}$ with $ S(Q,\nu) $ as in \eqref{eq.SQ-def},
			\item $ \Lss_2 := \set{\Sks(Q) : Q \in \Lks} $ with $ \Sks(Q) $ as in \eqref{eq.Sks}, and
			\item $ \Lss_3 := \set{\Ssp(Q) : Q \in \Lsp} $ with $ \Ssp(Q) $ as in \eqref{eq.Ssp}.
		\end{itemize}
		
	\end{definition}

	\renewcommand{\dq}{\delta_Q}

	\begin{proof}[Proof of Theorem \ref{thm.SFP}]

		Let $ \Lss $ be as in \eqref{eq.Lss}. We enumerate \begin{equation*}
		\Lss:= \set{S_\ell: \ell = 1, \cdots, L}.
		\end{equation*}
		We claim that the list $ S_1, \cdots, S_L $ satisfies the conclusions of Theorem \ref{thm.SFP}.
		
		We examine (A):
		\begin{itemize}
			\item Thanks to \eqref{eq.SQ bounded}, we have $ \#(S) \leq C $ for $ S \in \Lss_1 $.
			
			\item Thanks to \eqref{eq.Sks-bound}, we have $ \#(S) \leq C $ for $ S \in \Lss_2 $.
			
			\item Thanks to \eqref{eq.Ssp-bound}, we have $ \#(S) \leq C $ for $ S \in \Lss_3 $.
		\end{itemize}
		Therefore, conclusion (A) holds.
		
		Now we examine (B):
		\begin{itemize}
			\item Thanks to \eqref{eq.Sb-B}, $ \#(\Lss_1) \leq C\cdot \#(E) $.
			
			\item Thanks to Lemma \ref{lem.keystone}, $ \#(\Lss_2) \leq C\cdot \#(E) $.
			
			\item Thanks to Lemma \ref{lem.Lsp}, $ \#(\Lss_3) \leq C\cdot \#(E) $.
		\end{itemize}
		Therefore, conclusion (B) holds.

		We now turn to conclusion (C). Set
		\begin{equation*}
		M := \max\limits_{\ell = 1, \cdots, L} \norm{f}_{C^2_+(S_\ell)}
		\end{equation*}
		It suffices to show that there exists $ F \in C^2_+(\R^2) $ such that $ F|_E = f $ and $ \norm{F}_{C^2(\R^2)} \leq CM $. 
		
		By the definition of $ M $, we have $ \norm{f}_{C^2_+(S_\ell)} \leq M $ for all $ \ell = 1, \cdots, L $. This implies the following.
		\begin{itemize}
			\item Recall Definitions \ref{def.Sb}, \ref{def.SQ}, and \ref{def.KQ}. For each $ Q \in \Ls $, we have
			\begin{equation*}
			\eqindent
			\K(Q,\nu,CM) \neq \void
			\enskip
			\text{ and }
			\enskip
			\norm{f}_{C^2_+(\Sb(Q,\nu))} \leq M
			\text{ for } \nu = 1, \cdots, \nu(Q)\,.
			\end{equation*}
			This follows from the fact that $ S(Q,\nu) \in \Lss_1 \subset \Lss $ for $ \nu = 1, \cdots, \nu(Q) $ and $ \K(Q,\nu,CM) = \Gamma_+(\xqs,S(Q,\nu),CM) $ (see Definition \ref{def.KQ}). Therefore, the hypotheses of Lemma \ref{lem.FP-loc-S} are satisfied.
			
			\item For $ Q \in \Lks $ and $ \xqs $ as in Lemma \ref{lem.rep}, $ \Gamma_+(\xqs,\Sks(Q),CM) \neq \void $. This follows from the fact that $ \Sks(Q) \in \Lss_2 \subset \Lss $ for $ Q \in \Lks $.
			
			\item For $ Q \in \Lsp $ and $ \xqs $ as in Lemma \ref{lem.rep}, $ \Gamma_+(\xqs,\Ssp(Q),CM) \neq \void $ . This follows from the fact that $ \Ssp(Q) \in \Lss_3 \subset \Lss $ for $ Q \in \Lsp $.

		\end{itemize}

		We distinguish three types of squares $Q\in \Lz $.
		
		\newcommand{\pqs}{{P_Q^\sharp}}
		\newcommand{\fqs}{F_Q^\sharp}
		\newcommand{\xmqs}{{x_{\mu(Q)}^\sharp}}

		\begin{enumerate}[\text{Type }1]
			\item Suppose $ E \cap Q^* \neq \void $, that is, $ Q \in \Ls $. We set $ F_Q^\sharp:= F_Q $, where $ F_Q $ is as in Lemma \ref{lem.FP-loc-S}. In particular, we have
			\begin{equation}
			\eqindent
			P_Q^\sharp:= \jet{\xqs}{\fqs}  \in \K(Q,\nu(Q),CM) = \Gamma_+(\xqs,S(Q,\nu(Q)),CM),
			\label{eq.T1}
			\end{equation}
			with $ \xqs $ as in Lemma \ref{lem.rep} and $ S(Q,\nu(Q)) $ as in \eqref{eq.SQ-def}.
			
			\item Suppose $ E \cap Q^* = \void $ but $ \delta_Q < 1 $. Let $\Lsp, \mu $ be as in Lemma \ref{lem.Lsp}. 
			\begin{itemize}
				\item Suppose $ Q \notin \Lsp $. Pick 
				\begin{equation}
				\eqindent
				P_Q^\sharp \in \Gamma_+(\xmqs, \Sks(\mu(Q)), CM),
				\label{eq.T2-1}
				\end{equation}
				with $ x_{\mu(Q)}^\sharp $ as in Lemma \ref{lem.rep} and $ \Sks(\mu(Q)) $ as in \eqref{eq.Sks}.
				We set $ \fqs := \ew^{\xmqs}(\pqs) $, where $ \ew^{\xmqs} $ is as in Lemma \ref{lem.WE-map}.
				
				\item Suppose $ Q \in \Lsp $. Pick
				\begin{equation}
				\eqindent
				\pqs \in \Gamma_+(\xqs, \Ssp(Q),CM),
				\label{eq.T2-2}
				\end{equation}
				with $ \xqs $ as in Lemma \ref{lem.rep} and $ \Ssp(Q) $ as in \eqref{eq.Ssp}.
				We set $ \fqs := \ew^{\xqs}(\pqs) $, where $ \ew^{\xqs} $ is as in Lemma \ref{lem.WE-map}.
				
			\end{itemize}

			\item Suppose $ E \cap Q^* = \void $ and $ \delta_Q = 1 $. We set $ F_Q^\sharp:\equiv 0 $. 
		\end{enumerate}
	
		To wit, we associate Type 1 squares with clusters in $ \Lss_1 $, Type 2 non-special squares with clusters in $ \Lss_2 $, and Type 2 special squares with clusters in $ \Lss_3 $.

		Let $ \set{\theta_Q : Q \in \knice} $ be a $ C^2 $ partition of unity that is CZ compatible with $ \knice $. 
	
		We set
		\begin{equation*}
		F(x) := \sum_{Q \in \knice}\theta_Q(x) \cdot F_Q^\sharp(x)\,.
		\end{equation*}

		By construction, $ F_Q^\sharp \geq 0 $ on $ 100Q $ and $ F_Q^\sharp \big|_{E \cap Q^*} = f $ for each $ Q \in \Lz $. Therefore, $ F(x) \geq 0 $ and $ F = f $ on $ E $. 
		
		Now we estimate the derivatives of $ F $.

		Let $ x \in \mathbb{R}^2 $. Then there exists $ Q \in \knice $ such that $ Q \ni x $. We have
		\begin{equation}
		\d^\alpha F(x) = \sum_{Q' \leftrightarrow Q}\d^\alpha F_{Q'}^\sharp (x) \cdot \theta_{Q'}(x)
		+
		\sum_{Q'\leftrightarrow Q}\sum_{0 < \beta \leq \alpha}\binom{\alpha}{\beta}\d^{\alpha - \beta}(F_{Q'}^\sharp - F_Q^\sharp)(x)\cdot \d^\beta \theta_{Q'}(x)\,.
		\label{eq.last}
		\end{equation}

		\newcommand{\fqsp}{F_{Q'}^\sharp}
		\begin{claim}\label{last claim}
			Fix $ x \in \mathbb{R}^2 $. Let $ Q \ni x $, and let $ Q' \in \Lz $ with $ Q' \leftrightarrow Q $. Then
			\begin{equation}
			\abs{\d^\alpha (\fqs - \fqsp)(x)} \leq CM\dq^{2-\abs{\alpha}}
			\text{ for }
			\abs{\alpha} \leq 2.
			\label{eq.fqs}
			\end{equation}
		\end{claim}

		Suppose the claim is true. Then applying Lemma \ref{lem.WE-map}, Lemma \ref{lem.FP-loc-S}, and \eqref{eq.fqs} to estimate \eqref{eq.last}, we can conclude that $ \norm{F}_{C^2(\R^2)} \leq CM $. 
		
		Therefore, it suffices to prove Claim \ref{last claim}.

		\begin{proof}[Proof of Claim \ref{last claim}]
			
			By the triangle inequality, we can write 
			\begin{equation}
			\begin{split}
			\abs{\d^\alpha (F_Q^\sharp - F_{Q'}^\sharp)(x)} &\leq
			\abs{\d^\alpha (F_Q^\sharp - \jet{\xqs}{F_Q^\sharp})(x) } + \abs{\d^\alpha (F_{Q'}^\sharp - \jet{x_{Q'}^\sharp}{F_{Q'}^\sharp})(x) } \\&\,\,\,\,\,\,\,\,\,\,\,\,\,\,\,\,\,\,\,\,\,\,\,\,\,\,\,\,\,\,\,\,\,\,\,\,\,\,\,\,\,\,+ \abs{\d^\alpha (\jet{\xqs}{F_Q^\sharp} - \jet{x_{Q'}^\sharp}{F_{Q'}^\sharp})(x) }
			\\ 
			&=:\eta_1 + \eta_2 + \eta_3.
			\end{split}
			\end{equation}
			
			By Lemma \ref{lem.WE-map}, Lemma \ref{lem.FP-loc-S}, and Taylor's theorem, we have
			\begin{equation}
			\eta_1 + \eta_2 \leq CM\delta_Q^{2-\abs{\alpha}}.
			\label{eq.last-1}
			\end{equation}
			
			We want to show that
			\begin{equation}
			\eta_3 = \abs{\d^\alpha (\jet{x_Q^\sharp}{F_Q^\sharp - \jet{x_{Q'}^\sharp}{F_{Q'}^\sharp}})(x)} = \abs{\d^\alpha(P_Q^\sharp - P_{Q'}^\sharp)(x)} \leq CM\dq^{2-\abs{\alpha}}.
			\label{eq.eta-3}
			\end{equation}

			We consider the following cases. 
			
			\begin{enumerate}[\text{Case} 1]
				
				\item Suppose either $ Q $ or $ Q' $ is of Type 3. Then \eqref{eq.eta-3} follows from Lemma \ref{lem.WE-map}, Lemma \ref{lem.FP-loc-S}, and Taylor's theorem.
				
				
				\item Suppose both $ Q $ and $ Q' $ are of Type 1, that is, $ Q, Q' \in \Ls $. Then \eqref{eq.eta-3} follows from \eqref{eq.T1}, scenario (A) of Corollary \ref{lem.K-K}, and Taylor's theorem.

				\item Suppose one of $ Q, Q' $ is of Type 1 and the other is of Type 2. Without loss of generality, we may assume $ Q \in \Ls $ and $ Q' \in \Lz $. Recall $ \Lsp $ from Lemma \ref{lem.Lsp}.
				
				\begin{enumerate}[\text{Case 3-}a]
					\item Suppose $ Q' \notin \Lsp $. Then \eqref{eq.eta-3} follows from \eqref{eq.T1}, \eqref{eq.T2-1}, scenario (B) of Corollary \ref{lem.K-K}, and Taylor's theorem.
					
					\item Suppose $ Q' \in \Lsp $. Then \eqref{eq.eta-3} follows from \eqref{eq.T1}, \eqref{eq.T2-2}, scenario (C) of Corollary \ref{lem.K-K}, and Taylor's theorem.
				\end{enumerate}

				\item Suppose both $ Q, Q' $ are of Type 2. 
				
				\begin{enumerate}[\text{Case 4-}a]
					\item Suppose $ Q, Q' \notin \Lsp $. Then \eqref{eq.eta-3} follows from \eqref{eq.T2-1},  scenario (D) of Corollary \ref{lem.K-K}, and Taylor's theorem.

					\item Suppose $ Q \in \Lsp $ and $ Q' \notin \Lsp $. Then \eqref{eq.eta-3} follows from \eqref{eq.T2-1}, \eqref{eq.T2-2}, scenario (E) of Corollary \ref{lem.K-K}, and Taylor's theorem.
					
					\item Suppose $ Q, Q' \in \Lsp $. Then \eqref{eq.eta-3} follows from \eqref{eq.T2-2}, scenario (F) of Corollary \ref{lem.K-K}, and Taylor's theorem.

				\end{enumerate}
				
			\end{enumerate}
			
		This proves the claim.	
			
		\end{proof}
		
		The theorem is proved.
		
	\end{proof}

	\appendix

	\bigskip
	\bigskip

	\bibliographystyle{amsplain}

\end{document}